\documentclass[preprint]{elsarticle}

\usepackage{amsmath}
\usepackage{amsthm}
\usepackage{amssymb}		
\usepackage{eucal}		

\renewcommand\ge\geqslant
\renewcommand\le\leqslant

\newtheorem{theorem}{Theorem}[section]
\newtheorem{corollary}[theorem]{Corollary}
\newtheorem{proposition}[theorem]{Proposition}
\newtheorem{lemma}[theorem]{Lemma}

\theoremstyle{definition}
\newtheorem{definition}[theorem]{Definition}
\newtheorem{example}[theorem]{Example}

\theoremstyle{remark}

\newtheorem{remarkaftertheorem}{Remark}[theorem]

\numberwithin{equation}{section}

\newcommand{\divides}{\mid}
\newcommand{\rsmcoloneq}{:=}

\makeatletter
\newenvironment{sizeequation}[1]{%
  \skip@=\baselineskip
  #1%
  \baselineskip=\skip@
  \equation
}{\endequation \ignorespacesafterend}
\newenvironment{sizedisplaymath}[1]{%
  \skip@=\baselineskip
  #1%
  \baselineskip=\skip@
  \displaymath
}{\enddisplaymath \ignorespacesafterend}
\newenvironment{sizegather}[1]{%
  \skip@=\baselineskip
  #1%
  \baselineskip=\skip@
  \gather
}{\endgather \ignorespacesafterend}
\newenvironment{sizealign}[1]{%
  \skip@=\baselineskip
  #1%
  \baselineskip=\skip@
  \align
}{\endalign \ignorespacesafterend}  
\newenvironment{sizemultline}[1]{%
  \skip@=\baselineskip
  #1%
  \baselineskip=\skip@
  \multline
}{\endmultline \ignorespacesafterend}  
\newenvironment{sizealignat}[2]{%
  \skip@=\baselineskip
  #1%
  \baselineskip=\skip@
  \alignat
  #2%
}{\endalignat \ignorespacesafterend}  
\makeatother

\begin{document}

\begin{frontmatter}
\title{The Uniformization of Certain Algebraic\\ Hypergeometric Functions}
\author{Robert S. Maier}
\ead{rsm@math.arizona.edu}
\ead[url]{http://math.arizona.edu/\~{}rsm}
\address{Depts.\ of Mathematics and Physics, University of Arizona, Tucson,
AZ 85721, USA}
\begin{abstract}
The hypergeometric functions ${}_nF_{n-1}$ are higher transcendental
functions, but for certain parameter values they become algebraic, because
the monodromy of the defining hypergeometric differential equation becomes
finite.  It is shown that many algebraic ${}_nF_{n-1}$'s, for which the
finite monodromy is irreducible but imprimitive, can be represented as
combinations of certain explicitly algebraic functions of a single
variable; namely, the roots of trinomials.  This generalizes a result of
Birkeland, and is derived as a corollary of a family of binomial
coefficient identities that is of independent interest.  Any tuple of roots
of a trinomial traces~out a projective algebraic curve, and it is also
determined when this so-called Schwarz curve is of genus zero and can be
rationally parametrized.  Any such parametrization yields a hypergeometric
identity that explicitly uniformizes a family of algebraic ${}_nF_{n-1}$'s.
Many examples of such uniformizations are worked~out explicitly.  Even when
the governing Schwarz curve is of positive genus, it is shown how it is
sometimes possible to construct explicit single-valued or multivalued
parametrizations of individual algebraic ${}_nF_{n-1}$'s, by parametrizing
a quotiented Schwarz curve.  The parametrization requires computations in
rings of symmetric polynomials.
\end{abstract}
\end{frontmatter}

\section{Introduction}
\label{sec:intro}

The hypergeometric functions ${}_{n}F_{n-1}(\zeta)$, $n\ge1$, are
parametrized special functions of fundamental importance.  Each
${}_{n}F_{n-1}(\zeta)$ is a function of a single complex variable, and
in~general is a higher transcendental function.  It is parametrized by
complex numbers $a_1,\dots,a_n;\allowbreak b_1,\dots,b_{n-1}$, and
is written as ${}_{n}F_{n-1}(a_1,\dots,a_n;b_1,\dots,b_{n-1};\zeta)$.  It
is analytic on $\left|\zeta\right|<1$, with definition
\begin{equation}
\label{eq:Fdef}
{}_nF_{n-1}\left(
\begin{array}{l}
a_1,\dots,a_n\\ b_1,\dots,b_{n-1}
\end{array}
\!\!\biggm|\zeta
\right)=\sum_{k=0}^\infty \frac{(a_1)_k\dotsm (a_n)_k}{(b_1)_k\dotsm
(b_{n-1})_k(1)_k}\,\zeta^k.
\end{equation}
Here $(c)_k \rsmcoloneq c(c+1)\cdots(c+k-1)$, and the lower parameters
$b_1,\dots,b_{n-1}$ may not be non-positive integers.  The $n=1$ function
${}_1F_0(a_1;\text{--};\zeta)$ equals $(1-\nobreak\zeta)^{-a_1}$, and the
$n=2$ function ${}_2F_1(a_1,a_2;b_1;\zeta)$ is the Gauss hypergeometric
function.

If its parameters are suitably chosen, ${}_{n}F_{n-1}(\zeta)$ will become
an \emph{algebraic} function of~$\zeta$.  Equivalently, if one regards
${}_{n}F_{n-1}$ as a single-valued function on a certain Riemann surface,
defined by continuation from the disk $|\zeta|<1$, then the surface will
become compact.  If $n=1$, this occurs when $a_1\in\mathbb{Q}$.  In the
first nontrivial case $n=2$, the characterization of the triples
$(a_1,a_2;b_1)$ for which ${}_2F_1(a_1,\nobreak a_2;\allowbreak b_1;\zeta)$
is algebraic is a classical result of Schwarz.  There is a finite list of
possible \emph{normalized} triples (the famous `Schwarz list'), and
${}_2F_1$ is algebraic iff $(a_1,a_2;b_1)$ is a denormalized version of a
triple on the list.  Denormalization involves integer displacements of the
parameters.  For specifics, see~\cite[\S\,2.7.2]{Erdelyi53}.

More recently, Beukers and Heckman~\cite{Beukers89} treated $n\ge3$, and
obtained a complete characterization of the parameters
$(a_1,\dots,a_n;\allowbreak b_1,\dots,b_{n-1})$ that yield algebraicity.
Like the $n=2$ result of Schwarz, their result was based on the fact that
the function ${}_{n}F_{n-1}(a_1,\dots,a_n;\allowbreak
b_1,\dots,b_{n-1};\zeta)$ satisfies an order-$n$ differential equation on
the Riemann sphere~$\mathbb{P}^1_\zeta$, called
$E_n(a_1,\dots,a_n;\allowbreak b_1,\dots,b_{n-1},1)_\zeta$ below.  In
modern language, $E_n$~specifies a flat connection on an
$n$\nobreakdash-dimensional vector bundle over~$\mathbb{P}_\zeta^1$.  It
has singular points at $\zeta=0,1,\infty$, and its (projective) monodromy
group is generated by loops about these three points.  The monodromy,
resp.\ projective monodromy group is a subgroup of ${\it
  GL}_n(\mathbb{C})$, resp.\ ${\it PGL}_n(\mathbb{C})$, and algebraicity
occurs iff the monodromy is finite.  Schwarz exploited the classification
of the finite subgroups of ${\it PGL}_2(\mathbb{C})$.  In a \emph{tour
  de~force}, Beukers and Heckman handled the $n\ge3$ case by exploiting the
Shephard--Todd classification of the finite subgroups of~${\it
  GL}_n(\mathbb{C})$ generated by complex reflections.  Their
characterization result, however, is non-constructive: it supplies an
algorithm for determining whether a given function
${}_{n}F_{n-1}(a_1,\dots,a_n;\allowbreak b_1,\dots,b_{n-1};\zeta)$ is
algebraic in~$\zeta$, but it does not yield a polynomial equation (with
coefficients polynomial in~$\zeta$) satisfied by the function.

In this paper, a new approach is taken to the problem of constructing
equations satisfied by algebraic ${}_{n}F_{n-1}$'s.  Several classes of
hypergeometric function, known to be algebraic by Beukers--Heckman, are
made explicit by being \emph{uniformized}.  Recall that an algebraic
function~$F$ may be of genus zero, i.e., may have a `uniformization,' or
parametrization, by rational functions.  That~is, one may have
$F(R_1(t))=R_2(t)$ for certain rational functions $R_1,R_2$, so that
formally, $F=\allowbreak R_2\circ\nobreak R_1^{-1}$.  In this case the
Riemann surface of $F=F(\zeta)$, on which $\zeta$ and~$F$ are single-valued
meromorphic functions, is isomorphic to the Riemann
sphere~$\mathbb{P}^1_t$.  A~sample result of this nature, obtained below,
is the following.  Let $n=p+\nobreak2$, where $p\ge1$~is odd.  Then
\begin{multline}
\label{eq:sample}
{}_nF_{n-1}\left(
\begin{array}{l}
  \tfrac{a}{n},\dots,\tfrac{a+(n-1)}{n}\\
  \tfrac{a+1}p,\dots,\tfrac{a+p}p,\tfrac12
\end{array}
\!\!\biggm| \frac{4n^n}{p^p}\, \frac{t^2(1-t^2)^{2p}\left[(1+t)^p+(1-t)^p\right]^p}{\left[(1+t)^n+(1-t)^n\right]^n}
\right)\\
= \tfrac12\left[(1+t)^{-2a}+(1-t)^{-2a}\right] 
\left[\frac{(1+t)^n+(1-t)^n}{(1+t)^p+(1-t)^p}\right]^a.
\end{multline}
Equation~(\ref{eq:sample}) holds as an equality for all $t$ in a
neighborhood of $t=0$, or equivalently if one expands each side about
$t=0$, as an equality between power series.  It is a \emph{hypergeometric
  identity}, which is more general than a uniformization of a single~$F$.
This is because $a\in\mathbb{C}$ is arbitrary.  (If~any lower
hypergeometric parameter is a non-positive integer, the identity must be
interpreted in a limiting sense.)  If~$a\in\mathbb{Q}$, it follows that
${}_nF_{n-1}\bigl(\tfrac{a}{n},\dots,\tfrac{a+(n-1)}{n};\allowbreak\tfrac{a+1}p,\dots,\tfrac{a+p}p,\tfrac12;\zeta\bigr)$
is algebraic in~$\zeta$.  If moreover $a\in\mathbb{Z}$, in which case
(\ref{eq:sample})~is a uniformization in the strict sense, this algebraic
${}_{n}F_{n-1}$ must be of genus zero.

Most of the algebraic ${}_nF_{n-1}$'s parametrized below, or more
specifically uniformized, are of the following type.  The ${}_nF_{n-1}$
in~(\ref{eq:sample}) is the $q=2$ case of
\begin{equation}
\label{eq:sample2}
  {}_nF_{n-1}\left(
\begin{array}{l}
  \tfrac{a}{n},\dots,\tfrac{a+(n-1)}{n}\\
  \tfrac{a+1}p,\dots,\tfrac{a+p}p,\tfrac1q,\dots,\tfrac{q-1}q
\end{array}
\!\!\biggm|\zeta
\right),
\end{equation}
where $n=p+q$.  The order-$n$ differential equation~$E_n$
on~$\mathbb{P}^1_\zeta$ associated to~(\ref{eq:sample2}),
\begin{equation}
\label{eq:impriv}
E_n\bigl(\tfrac{a}{n},\dots,\tfrac{a+(n-1)}{n};\,
\tfrac{a+1}p,\dots,\tfrac{a+p}p,\tfrac1q,\dots,\tfrac{q-1}q,1\bigr)_\zeta\,,
\end{equation}
plays a role in the Beukers--Heckman analysis of algebraicity
(see~\S\,\ref{sec:prelims}).  The treatment of (\ref{eq:sample2})
and~(\ref{eq:impriv}) relies on the following fact.  If $\gcd(p,\nobreak
q)=1$, the solution space of~(\ref{eq:impriv}) is spanned by
$x_1^\gamma(\zeta),\dots,x_n^\gamma(\zeta)$, where $\gamma\rsmcoloneq -q a$
and the $n$~algebraic functions $x_1(\zeta),\dots,x_n(\zeta)$ are the roots
of the \emph{trinomial equation}
\begin{equation}
\label{eq:trinomial}
  x^n - g\,x^p - \beta = 0.
\end{equation}
Here $g\neq0$ is arbitrary but fixed, and $\beta$ is determined implicitly
by
\begin{equation}
\label{eq:zetadef0}
\zeta=(-)^q \frac{n^n}{p^pq^q}\, \frac{\beta^q}{g^n}.
\end{equation}
(This statement assumes $\gamma\notin\mathbb{Z}$.)  As defined
by~(\ref{eq:zetadef0}), $\zeta$~is a projectivized and normalized version
of the \emph{discriminant} of~(\ref{eq:trinomial}): if $\beta\neq0$, the
trinomial has coincident roots if and only if $\zeta=1$.

Formulas relating such hypergeometric functions as~(\ref{eq:sample2}) to
trinomial roots can be traced to the 1758 work of Lambert, who expressed
such roots as finite sums of hypergeometric series.
(See~\cite[p.~72]{BerndtI}.)  Once trinomial roots have been expressed
in~terms of hypergeometric series, typically by Lagrange inversion,
formulas for hypergeometric functions in~terms of trinomial roots can be
derived.  Two especially useful reversed formulas of this sort were
obtained in the 1920s by Birkeland~\cite{Birkeland27}.
(See~\cite{Belardinelli59,Passare2004} for general reviews; the classical
literature is cited in~\cite{Hall41}, and some explicit series for the
roots of~(\ref{eq:trinomial}) are given in~\cite{Eagle39,Glasser2000}.)
For example, Birkeland expressed the algebraic function~(\ref{eq:sample2})
as a combination of only $q$~of the $n$~functions
$x_1^\gamma(\zeta),\dots,x_n^\gamma(\zeta)$.

The first result of this paper is the extension of Birkeland's two formulas
to the case when the ${}_nF_{n-1}$ being represented in~terms of trinomial
roots has a so-called parametric excess $S=\allowbreak
\sum_{i=1}^{n-1}b_i-\nobreak\sum_{i=1}^na_i$ that does not equal~$\frac12$,
as it does in~(\ref{eq:sample2}), but rather is an arbitrary
half-odd-integer.  There are in~fact two families of formulas indexed by
$\ell\in\mathbb{Z}$, with $S=\frac12-\nobreak\ell$, that subsume
Birkeland's `$\ell=0$' formulas.  The new families are shown to follow from
a family of binomial coefficient identities, which is derived first and is
of independent interest: it subsumes several combinatorial identities that
are listed individually in~\cite{Gould72}.  The derivation of the family of
binomial coefficient identities resembles that of Chu~\cite{Chu95}.  The
technique used is not Lagrange inversion (explicitly) but rather the
Vandermonde convolution transform of Gould~\cite{Gould61}.

The heart of this paper is a study of the trinomial
equation~(\ref{eq:trinomial}), focused on basic algebraic geometry rather
than on controlling or bounding its roots.  (For the latter,
see~\cite{Fell80} and papers cited therein.)  The key concept is that of a
\emph{Schwarz curve}.  For each coprime ${p,q\ge1}$ with
$n\rsmcoloneq\allowbreak p+\nobreak q$, this is defined following Kato and
Noumi~\cite{Kato2003} to be a complex projective curve
$\mathcal{C}^{(n)}_{p,q}\subset\mathbb{P}^{n-1}$ comprising all
${[x_1:\ldots: x_n]}\in\mathbb{P}^{n-1}$ such that $x_1,\dots, x_n$ are the
roots of \emph{some} equation of the form~(\ref{eq:trinomial}).  That~is,
$\mathcal{C}^{(n)}_{p,q}$~is the common zero set of the $n-2$ elementary
symmetric polynomials $\sigma_1,\dots,\sigma_{q-1}$ and
$\sigma_{q+1},\dots,\sigma_{n-1}$ in $x_1,\dots,x_n$.  Clearly
$g=(-)^{q-1}\sigma_q$ and $\beta=(-)^{n-1}\sigma_n$, giving a formula
for~$\zeta$; so there is a degree\nobreakdash-$n!$ covering
$\mathcal{C}^{(n)}_{p,q}\to\mathbb{P}^1_\zeta$.  The curve
$\mathcal{C}^{(n)}_{p,q}$ is irreducible~\cite[Cor.~4.7]{Kato2003}.

For each~$k$, $n-1\ge k\ge2$, a subsidiary Schwarz curve
$\mathcal{C}^{(k)}_{p,q}\subset\mathbb{P}^{k-1}$, also irreducible, is
defined to include all $[x_1:\ldots: x_k]\in\mathbb{P}^{k-1}$ such that
$x_1,\dots, x_k$ are $k$~of the $n$~roots of some equation of the
form~(\ref{eq:trinomial}).  The curve $\mathcal{C}^{(2)}_{p,q}$ is
$\mathbb{P}^1$ itself, and one can introduce a final curve
$\mathcal{C}^{(1)}_{p,q}$, also isomorphic to~$\mathbb{P}^1$.  These curves
are joined by maps $\phi_{p,q}^{(k)}\colon
\mathcal{C}_{p,q}^{(k)}\to\mathcal{C}_{p,q}^{(k-1)}$ with
$\deg\phi_{p,q}^{(k)}=n-k+1$, i.e.,
\begin{equation}
  \mathcal{C}^{(n)}_{p,q}\stackrel{\phi_{p,q}^{(n)}}{\longrightarrow}\mathcal{C}^{(n-1)}_{p,q}\stackrel{\phi_{p,q}^{(n-1)}}{\longrightarrow}\dots
  \stackrel{\phi_{p,q}^{(3)}}{\longrightarrow}\mathcal{C}^{(2)}_{p,q}\stackrel{\phi_{p,q}^{(2)}}{\longrightarrow}\mathcal{C}^{(1)}_{p,q}\stackrel{\phi_{p,q}^{(1)}}{\longrightarrow}\mathbb{P}^1_\zeta.
\end{equation}
Each covering $\mathcal{C}^{(k)}_{p,q}\to\mathbb{P}^1_\zeta$ is a Bely\u\i\
map, i.e., is ramified only over $\zeta=0,1,\infty$.

The Schwarz curves are not in general of genus zero, though
$\mathcal{C}^{(2)}_{p,q}$ and~$\mathcal{C}^{(1)}_{p,q}$ are.  One can write
$\mathcal{C}^{(2)}_{p,q}\cong\mathbb{P}^1_t$ and
$\mathcal{C}^{(1)}_{p,q}\cong\mathbb{P}^1_s$, where $t$ and~$s$ are
rational parameters, and make the concrete choice $t=\allowbreak
(x_1+\nobreak x_2)/\allowbreak (x_1-\nobreak x_2)$.  It is a consequence of
the rational parametrizations of $\mathcal{C}^{(2)}_{p,q}$
and~$\mathcal{C}^{(1)}_{p,q}$, substituted into the appropriate
Birkeland-style representation formulas, that when $q=2$ or ${q=1}$, the
function ${}_nF_{n-1}$ of~(\ref{eq:sample2}) can be parametrized by $t$
or~$s$, with its argument~$\zeta$ being respectively a
degree-$[n(n-\nobreak1)]$ rational function of~$t$, and a
degree\nobreakdash-$n$ rational function of~$s$.  This turns~out to be the
origin of the sample identity~(\ref{eq:sample}).  Its right side is the sum
of two terms, proportional to $(1+\nobreak t)^{-2a},\allowbreak(1-\nobreak
t)^{-2a}$; and these come from the exponentiated roots
$x_1^\gamma,x_2^\gamma$.

The fundamental fact is that if a hypergeometric function
${}_nF_{n-1}(\zeta)$ can be represented in~terms of a $k$\nobreakdash-tuple
of roots of the trinomial equation~(\ref{eq:trinomial}) with $n=p+\nobreak
q$, any parametrization of the curve $\mathcal{C}^{(k)}_{p,q}$ will yield a
parametrization of ${}_nF_{n-1}(\zeta)$; and if this curve is of genus
zero, with a rational parameter, the parametrization of
${}_nF_{n-1}(\zeta)$ will be by rational functions of the parameter.  Many
interesting examples of this are worked~out below, going well beyond the
sample identity~(\ref{eq:sample}).  The curves employed include `top'
Schwarz curves $\mathcal{C}^{(n)}_{p,q}$, as well as subsidiary curves
$\mathcal{C}^{(k)}_{p,q}$ with $k=1,2,3$.

To derive parametrizations of certain algebraic ${}_nF_{n-1}(\zeta)$'s, a
more sophisticated technique is useful.  If there is a representation of an
${}_nF_{n-1}(\zeta)$ in~terms of $k$~trinomial roots that is symmetric in
the roots, the parametrization is really by a \emph{quotiented} Schwarz
curve, obtained from~$\mathcal{C}^{(k)}_{p,q}$ by quotienting~out the
symmetric group~$\mathfrak{S}_k$.  Even if $\mathcal{C}^{(k)}_{p,q}$ is of
positive genus, the quotiented curve may be of genus zero.  Quotienting~out
the cyclic group~$\mathfrak{C}_k$ or the dihedral group~$\mathfrak{D}_k$
may also be useful.  Several examples of explicit parametrizations of
algebraic ${}_nF_{n-1}$'s that can be obtained by quotienting are
worked~out, by computation in rings of symmetric polynomials.

\smallskip
The paper is organized as follows.  The first half focuses on
hypergeometric equations and series, and series identities.  In
\S\,\ref{sec:prelims} the equations~$E_n$ are introduced, and their
monodromy is discussed.  In~\S\,\ref{sec:binomial}, the family of binomial
coefficient identities is derived with the aid of the Vandermonde
convolution transform.  In~\S\,\ref{sec:alghyper}, families of formulas of
Birkeland type for various hypergeometric functions, including algebraic
ones, are derived from these identities.  The main results are
Thm.~\ref{thm:invrep}, which extends Birkeland's two representations and is
parametrized by $\ell\in\mathbb{Z}$, and Thm.~\ref{thm:invrepg0}.  The
former expresses ${}_nF_{n-1}$'s of a type that
generalizes~(\ref{eq:sample2}) in~terms of trinomial roots.  The latter
expresses certain non-algebraic ${}_{n+1}F_n$'s in~terms of ${}_2F_1$'s
and~${}_3F_2$'s.

The second half is explicitly algebraic-geometric.
In~\S\S\,\ref{sec:curves} and~\ref{sec:genera}, the Schwarz curves are
introduced and studied.  Their ramification structures and genera are
computed, and the curves of genus zero are determined.
(Theorem~\ref{thm:genus}, the genus formula, generalizes a theorem
of~\cite{Kato2003}.)  In~\S\,7, many hypergeometric identities with a free
parameter $a\in\mathbb{C}$ are derived, including uniformizations of
algebraic ${}_nF_{n-1}$'s of genus zero.  A~typical result is
Thm.~\ref{thm:73}, of which Eq.~(\ref{eq:sample}) above is a special case;
it comes from~$\mathcal{C}^{(2)}_{p,2}$.  In~\S\,8, the quotiented Schwarz
curves are defined and in some cases parametrized, yielding
parametrizations in which $a\in\mathbb{Q}$ is fixed.  Some of these are
multivalued parametrizations with radicals.

\section{Monodromy and degeneracy}
\label{sec:prelims}

The parametrized functions
${}_{n}F_{n-1}(a_1,\dots,a_n;b_1,\dots,b_{n-1};\zeta)$ are locally defined
by Eq.~(\ref{eq:Fdef}).  Two sorts of degeneracy are of interest, for the
functions themselves, and for the equations that they satisfy.  The first
is when an upper and a lower parameter equal each other.  If~so, they may
be `cancelled', reducing the ${}_nF_{n-1}$ to an~${}_{n-1}F_{n-2}$.  The
second arises when one of the lower parameters is taken to an integer~$-m$,
and one of the upper ones to an integer~$-m'$, with $-m\le -m'\le0$.  Let
$[\zeta^{>m'}]F$ denote the sum of the terms proportional to~$\zeta^k$,
$k>m'$, in the series for~$F$; and let $(\mathfrak{a})$
and~$(\mathfrak{b})$ denote $(a_1,\dots,a_{n-1})$ and
$(b_1,\dots,b_{n-2})$.  The following auxiliary lemma permits such
hypergeometric identities as~(\ref{eq:sample}) to be interpreted in a
limiting sense; it will be used in~\S\,\ref{sec:nofreeparams}.

\begin{lemma}
\label{lem:1}
  Let\/ $a_n\to-m'$, $b_{n-1}\to-m$, with\/
  $(a_n+m')/(b_{n-1}+m)\to\alpha$.  Then
\begin{multline*}
    [\zeta^{>m'}]\,
  {}_nF_{n-1}\left(
\begin{array}{l}
  (\mathfrak{a}),\,a_n\\
  (\mathfrak{b}),\,b_{n-1}
\end{array}
\!\!\biggm|\zeta
\right)\longrightarrow \\
\alpha \, (-)^{m-m'}{\binom{m}{m'}}^{-1}
\frac{(\mathfrak{a})_{m+1}}{(\mathfrak{b})_{m+1}(1)_{m+1}}\,\zeta^{m+1}\, 
{}_nF_{n-1}\left(
\begin{array}{l}
  (\mathfrak{a})+m+1,\,m-m'+1\\
  (\mathfrak{b})+m+1,\,m+2
\end{array}
\!\!\biggm|\zeta
\right).
\end{multline*}
If\/ $m=m'$, this limit simplifies to
\begin{equation*}
\alpha\,\cdot\,[\zeta^{>m'}]\,
{}_{n-1}F_{n-2}\left(
\begin{array}{l}
  (\mathfrak{a})\\
  (\mathfrak{b})
\end{array}
\!\!\biggm|\zeta
\right).
\end{equation*}
\end{lemma}
\begin{proof}
Take the term-by-term limit of the series in~(\ref{eq:Fdef}).
\end{proof}

\smallskip
It is well known that the function ${}_nF_{n-1}(a_1,\dots,a_n;b_1,\dots
b_{n-1};\zeta)$, defined on the disk $\left|\zeta\right|<1$, satisfies the
$b_n=1$ case of the order-$n$ equation $\mathcal{D}_nF=0$, where
(with $D_\zeta \rsmcoloneq {\rm d}/{\rm d}\zeta$)
\begin{align}
\label{eq:En}
\mathcal{D}_n&=\mathcal{D}_n(a_1,\dots,a_n;\,b_1,\dots,b_n)\\
&=
(\zeta D_\zeta+b_1-1)\cdots(\zeta D_\zeta+b_n-1)
-
\zeta\,(\zeta D_\zeta+a_1)\cdots(\zeta D_\zeta+a_n).\notag
\end{align}
The equation~$\mathcal{D}_nF=0$ is denoted
by~$E_n(a_1,\dots,a_n;b_1,\dots,b_n)$ here; a subscript~$\zeta$ will be
added as appropriate to indicate the independent variable.  This equation
has three regular singular points on~$\mathbb{P}^1_\zeta$, namely
$\zeta=0,1,\infty$, with respective characteristic exponents $1-\nobreak
b_1,\dots,1-\nobreak b_n$; and $0,1,\dots,n-\nobreak 2,S$; and
$a_1,\dots,a_n$.  In~this,
\begin{equation}
  S = \left(\sum_{i=1}^n b_i\right) - \left(\sum_{i=1}^n a_i\right) -1
\end{equation}
is the `parametric excess,' which reduces to $\sum_{i=1}^{n-1} b_i
-\nobreak \sum_{i=1}^n a_n$ if $b_n=1$.  The canonical solution
${}_nF_{n-1}$, defined if none of $b_1,\dots,b_{n-1}$ is a non-positive
integer, is analytic at $\zeta=0$ and belongs to the exponent $1-\nobreak
b_n=0$.  Allowing $b_n$ to differ from unity in~$E_n$ is not a major
matter, since adding~$\delta$ to each parameter of~$E_n$ merely multiplies
all solutions by~$\zeta^{-\delta}$.  The following is a standard fact.

\begin{lemma}
\label{lem:2}
If $b_j-b_{j'}\notin\mathbb{Z}$, $\forall j,j'$, so that the singular
point\/ $\zeta=0$ is non-logarithmic, the solution space of\/ $E_n$
at\/ $\zeta=0$ is spanned by the functions
\begin{equation*}
\zeta^{1-b_j}\,{}_nF_{n-1}
\left(
\begin{array}{l}
a_1-b_j+1,\dots,a_n-b_j+1\\
b_1-b_j+1,\dots,\widehat{b_j-b_j+1},\dots, b_n-b_j+1  
\end{array}
\!\!\biggm|\zeta
\right), \quad j=1,\dots,n.
\end{equation*}
\end{lemma}

\smallskip
The differential equation~$E_n$ on~$\mathbb{P}^1_\zeta$ is the canonical
order-$n$ one with three regular singular points, the monodromy at one of
them being special.  This is seen as follows.  Associated to any~$E_n$ is a
monodromy representation of $\pi_1(X,\zeta_0)$, the fundamental group of
the triply punctured sphere $X=\mathbb{P}^1_\zeta\setminus\{0,1,\infty\}$.
(The base point~$\zeta_0$ plays no role, up~to a similarity
transformation.)  The image~$H$ of $\pi_1(X,\zeta_0)$ in~${\it
  GL}_n(\mathbb{C})$ is the monodromy group of the specified~$E_n$.  It~is
generated by $h_0,h_1,h_\infty$, the monodromy matrices around
$\zeta=0,1,\infty$, which satisfy $h_\infty h_1h_0=I$.  Their eigenvalues
are exponentiated characteristic exponents, so they have characteristic
polynomials
\begin{equation}
\begin{split}
  \det(\lambda I-h_\infty)=\prod_{j=1}^n (\lambda -\alpha_j),\qquad
  \det(\lambda  I-h_0^{-1})=\prod_{j=1}^n (\lambda -\beta_j),\\
  \det(\lambda I-h_1)=(\lambda-1)^{n-1}(\lambda -e^{2\pi{\rm i}S}),\hskip1in
\end{split}
\end{equation}
where $\alpha_j=e^{2\pi{\rm i}a_j}$, $\beta_j=e^{2\pi{\rm i}b_j}$.
Moreover, if $S\notin\mathbb{Z}$ then $h_1$~is
diagonalizable~\cite{Beukers89}, i.e., the singular point $\zeta=1$ is not
logarithmic.  Hence if $S\notin\mathbb{Z}$, $h_1$~will act
on~$\mathbb{C}^n$ as multiplication by $e^{2\pi{\rm i}S}$ on a
1\nobreakdash-dimensional subspace, and as the identity on a complementary
$(n-\nobreak1)$-dimensional subspace.  It will be a \emph{complex
  reflection}.

The monodromy representation is irreducible iff no upper parameter~$a_j$
and lower one~$b_{j'}$ differ by an integer, i.e.,
$\alpha_j\neq\beta_{j'}$, $\forall j,j'$; and it is unaffected (up~to
similarity transformations) by integer displacements of
parameters~\cite{Beukers89}.  If exactly one upper and one lower parameter
of ${}_nF_{n-1}(a_1,\dots,a_n;b_1,\dots,b_{n-1};\zeta)$ are equal,
permitting cancellation, then the corresponding~$E_n$ will have reducible
monodromy.  One can show that in this case, the differential
operator~$\mathcal{D}_n$ will have an order-$(n-1)$ right factor, which
comes from the ${}_{n-1}F_{n-2}$ to which the ${}_{n}F_{n-1}$ reduces.  In
the theory of the Gauss function~${}_2F_1$, the reducible case, when $a_1$
or~$a_2$ differs by an integer from $b_1$ or~$b_2=1$, is often called the
`degenerate' case~\cite[\S\,2.2]{Erdelyi53}.  Reducibility of~$E_2$
facilitates the explicit construction of solutions~\cite{Kimura91}.

A monodromy group~$H<{\it GL}_n(\mathbb{C})$ of an equation~$E_n$, if
irreducible, is said to be \emph{imprimitive} if there is a direct sum
decomposition $V_1\oplus\dots\oplus V_d$ of~$\mathbb{C}^n$ with $d\ge2$ and
$\dim V_j\ge1$, such that each element of $H$~permutes the spaces~$V_j$.
The following is a key theorem of Beukers and Heckman.  It refers to the
complex reflection subgroup~$H_r$ of~$H$, which is generated by the
elements $h_\infty^k\cdot h_1\cdot h_\infty^{-k}$, $k\in\mathbb{Z}$.
A~`scalar shift' by $\delta\in\mathbb{C}$ means the addition of~$\delta$ to
each parameter.

\begin{theorem}[cf.~\cite{Beukers89},\ Thm.~5.8]
\label{thm:BH0}
Suppose that\/ $H,H_r$, computed from\/ $E_n$, are irreducible in\/ ${\it
  GL}_n(\mathbb{C})$.  Then\/ $H$ will be imprimitive if and only if there
exist relatively prime\/ $p,q\ge1$ with\/ $p+\nobreak q=n$, and\/
$a\in\mathbb{C}$, such that\/ $E_n$ takes on one of the two forms
\begin{equation*}
  E_n\left(
  \begin{array}{l}
  \tfrac{a}{n},\dots,\tfrac{a+(n-1)}{n}\\
  \tfrac{a+1}p,\dots,\tfrac{a+p}p;\,\tfrac1q,\dots,\tfrac{q}q
  \end{array}
  \right),
  \quad
  E_n\left(
  \begin{array}{l}
  \tfrac{-a}p,\dots,\tfrac{-a+(p-1)}p;\,\tfrac{a}q,\dots,\tfrac{a+(q-1)}q\\
  \tfrac{1}{n},\dots,\tfrac{n}{n}
  \end{array}
  \right),
\end{equation*}
up to integer displacements of parameters, and up to a scalar shift by
some\/ $\delta\in\mathbb{C}$.  {\rm(}For irreducibility, one must have
that\/ $qa\notin\mathbb{Z}$, resp.\ $na\notin\mathbb{Z}$; this condition
ensures that the upper and lower parameters satisfy\/ $a_j-\nobreak
b_{j'}\notin\mathbb{Z}$, $\forall j,j'$.{\rm)} Moreover, if\/
$a\in\mathbb{Q}$ then\/ $H$ will be finite.
\end{theorem}

This theorem yields a class of algebraic ${}_nF_{n-1}$'s, including
(\ref{eq:sample2}).  One must choose the shift~$\delta$ so that one of the
lower parameters equals~$1$; e.g., $\delta=0$.  Beukers and Heckman do not
compute the monodromy group~$H$ of the two $E_n$'s in the theorem, but it
follows readily from their proof that if, e.g., $a=-1/mq$,
resp.\ $a=-1/mn$, with $m\ge2$, then $H$~is of order $m^{n-1}n!$, and is
isomorphic to the `symmetric' index-$m$ subgroup of the wreath product
${\mathfrak{C}}_m \wr \mathfrak{S}_n$, where ${\mathfrak{C}}_m$
and~$\mathfrak{S}_n$ are the usual cyclic and symmetric groups.  So is the
corresponding projective monodromy group $\overline H < {\it
  PGL}_n(\mathbb{C})$ (see~\cite[Cor.~4.6]{Kato2003}).

If $a\in\mathbb{Z}$ then Thm.~\ref{thm:BH0} does not apply, due to
reducibility: in each~$E_n$ an upper and a lower parameter will differ by
an integer.  But, one has the following.

\begin{theorem}[cf.~\cite{Beukers89},\ Prop.~5.9]
\label{thm:BH1}
If there are \emph{equal} upper and lower parameters in either\/ $E_n$ of
Theorem\/~{\rm\ref{thm:BH0}}, and\/ $a\in\mathbb{Z}$ {\rm(}e.g., if\/
$a=\pm1${\rm)}, then the\/ $E_{n-1}$ obtained by cancelling them will have
monodromy group\/ $H<{\it GL}_{n-1}(\mathbb{C})$ isomorphic
to\/~$\mathfrak{S}_n$.
\end{theorem}

This yields a class of algebraic ${}_{n-1}F_{n-2}$'s.  The hypergeometric
functions appearing in the following sections will include algebraic
${}_nF_{n-1}$'s of the type arising from Thm.~\ref{thm:BH0}, but for
certain choices of parameter, they will reduce to ${}_{n-1}F_{n-2}$'s.  On
the power series level, Lemma~\ref{lem:1} will perform some of these
reductions.

\section{Sequence transformations and coefficient identities}
\label{sec:binomial}

In this section and~\S\,\ref{sec:alghyper}, it is shown that the solutions
of hypergeometric equations~$E_n$ of the types introduced in
Thm.~\ref{thm:BH0}, with monodromy that is irreducible but imprimitive,
include algebraic functions that are solutions of trinomial equations.
Explicit expressions for the associated ${}_nF_{n-1}$'s as combinations of
algebraic functions will be derived, which extend those of
Birkeland~\cite{Birkeland27}.

The chief result of this section is Theorem~\ref{thm:cor1}, which provides
a family of binomial coefficient identities, the family being indexed by
$\ell\in\mathbb{Z}$.  They are derived from the standardized trinomial
equation $y-\nobreak 1-\nobreak zy^B=0$, and will be used
in~\S\,\ref{sec:alghyper}.  Theorem~\ref{thm:g0} provides two more
identities, which are more sophisticated in the sense that they have an
additional free parameter.  Some of these identities were found by Lambert,
Ramanujan, P\'olya, and Gould, and the family indexed by
$\ell\in\mathbb{Z}$ was explored by Chu~\cite{Chu95} in a little-noticed
paper.  The tool employed in deriving the identities of Theorems
\ref{thm:cor1} and~\ref{thm:g0} is the Vandermonde convolution transform of
Gould~\cite{Gould56,Gould61}, which is useful in deriving identities
relating coefficients of related series, such as binomial coefficient
identities.

In what follows,
\begin{equation}
  \binom{a}{r} \rsmcoloneq  (a-r+1)_r / r!
\end{equation}
extends the binomial coefficient to the case of arbitrary upper
parameter, and
\begin{align}
  (a)_r \rsmcoloneq  
  \begin{cases}
    (a)\dotsb(a+r-1),&r\ge0;\\
    \left[(a-s)\dotsb(a-1)\right]^{-1},&r=-s\le0,
  \end{cases}
\end{align}
extends the usual rising factorial, so that $(a)_r =
\left[(a+r)_{-r}\right]^{-1}$ for all $r\in\mathbb{Z}$.  The standard
forward finite difference operator $\Delta_{A,B}$, which acts on functions
of~$A$, is defined by $\Delta_{A,B}[h](A) =\allowbreak h(A+\nobreak
B)-\nobreak h(A)$.  Its $n$'th power $\Delta_{A,B}^n$ satisfies
\begin{equation}
\label{eq:333}
  \Delta_{A,B}^n[h](A)=\sum_{k=0}^n (-)^{n-k}\binom{n}{k}\,h(A+kB).
\end{equation}

\begin{theorem}[\cite{Gould61}]
\label{thm:mainGould}
Let\/ $A,B\in\mathbb{C}$ be given, with\/ $B\neq0,1$.  Let\/ $f(k)$,
$k\ge0$, be an infinite sequence of complex numbers, and let the sequence
$\hat f(n)$, $n\ge0$, be defined by
\begin{equation}
\label{eq:vctransform}
  \hat f(n) = \sum_{k=0}^n (-)^k\binom{n}{k}\binom{A+Bk}{n}f(k).
\end{equation}
In a neighborhood of\/ $z=0$, let\/ $y$ near\/ $1$ be defined implicitly
by the standardized trinomial equation\/ $y -\nobreak 1 -\nobreak z\,y^B=0$.  Then
\begin{equation}
\label{eq:extraeqn2}
\sum_{k=0}^\infty \binom{A+Bk}{k}f(k)\,z^k
=
y^A
\sum_{n=0}^\infty \hat f(n)  \left(\frac{1-y}{y}\right)^n
\end{equation}
holds as an equality between power series in\/ $z$; and hence, if either
has a positive radius of convergence, as an equality in a neighborhood of\/
$z=0$.
\end{theorem}
If the theorem applies, the sequence $\hat f(n)$, $n\ge0$, will be called
the $(A,B)$-Vandermonde convolution transform of the sequence $f(k)$,
$k\ge0$.  If $f(0)=1$, then $\hat f(0)=1$ also.  For the inverse transform,
see~\cite{Gould61}.

\begin{definition}
\label{def:fdef}
For each $\ell\in\mathbb{Z}$, the sequence $f_\ell(A,B;k)$, $k\ge0$, is
defined by
\begin{equation}
\label{eq:fdef}
f_\ell(A,B;\,k) \rsmcoloneq  \frac{(A+Bk+1)_{\ell-1}}{(A+1)_{\ell-1}}
= \frac{(A+\ell)_{1-\ell}}{(A+Bk+\ell)_{1-\ell}}.
\end{equation}
It is normalized so that $f_\ell(A,B;0)=1$.  Note that $f_1(A,B;k)\equiv1$.
\end{definition}

It is an immediate consequence of~(\ref{eq:333}) that the
$(A,B)$-Vandermonde convolution transform of $f_\ell(A,B;k)$, defined by
Eq.~(\ref{eq:vctransform}), has the representation
\begin{equation}
\label{eq:defconvol}
\begin{split}
\hat f_\ell(A,B;\,n) &= 
\frac{(-)^n}{(A+1)_{\ell-1}}
\,\Delta_{A,B}^n
\left[\binom{A}{n} (A+1)_{\ell-1}\right] \\
&=\frac{(-)^n}{n!\,(A+1)_{\ell-1}}
\,\Delta_{A,B}^n
\left[(A-n+1)_{n+\ell-1}\right].
\end{split}
\end{equation}
Here and in subsequent equations and identities, parameters such as $A,B$
are constrained, sometimes tacitly, so that division by zero occurs on
neither side.

\begin{theorem}
\label{thm:ratchar}
  {\rm(i)} For each\/ $\ell=-m\le0$, $\hat f_\ell(A,B;n)$ is nonzero only
  if\/ $n=0,\dots,m$, and is rational in\/ $A,B$ for each\/~$n$.
  {\rm(ii)} For each\/ $\ell\ge1$, $\hat
  f_\ell(A,B;n)$ equals a certain degree-$(\ell-\nobreak1)$ polynomial in\/
  $n$ with coefficients polynomial in\/ $A,B$, multiplied by the function\/
  $(n+\nobreak1)_{\ell-1}(-B)^n/(A+\nobreak1)_{\ell-1}$ of\/~$n$.
\end{theorem}

\begin{proof}
  Part~(i) follows from the fact that if $n\ge m+1$, then $(A-\nobreak
  n+\nobreak1)_{n-m-1}$, appearing in~(\ref{eq:defconvol}), is a polynomial
  in~$A$ of degree $n-\nobreak m-\nobreak1$; hence its $n$'th finite
  difference must equal zero.  Part~(ii) comes from a finite difference
  computation (cf.~\cite[\S\,76]{Jordan65}).  Explicitly, $\hat
  f_\ell(A,B;n)$ equals $(-B)^n/(A+\nobreak1)_{\ell-1}$ times
  \begin{equation}
    \sum_{i=0}^{\ell-1}\sum_{j=0}^i s(n+\ell-1,n+i)\,S(n+i,n+j)\,(n+1)_j\:B^i\binom{(A+\ell-1)/B}{j},
  \end{equation}
  where $s,S$ are the Stirling numbers of the first and second kinds.  But
  $s(\nu,\nu-\nobreak\rho)$, $S(\nu,\nu-\nobreak\rho)$ are degree-$2\rho$
  polynomials in~$\nu$ that are divisible by $(\nu-\rho+1)_\rho$, so each
  summand is of degree at~most $2(\ell-\nobreak1)$ in~$n$ and is divisible
  by~$(n+\nobreak1)_{\ell-1}$.
\end{proof}

\begin{remarkaftertheorem}
  The formula given in the preceding proof, valid when $\ell\ge1$, can be
  written as
  \begin{equation}
    \label{eq:rep2}
    \hat f_\ell(A,B;n) = \frac{(-)^n}{(A+1)_{\ell-1}}
    \sum_{j=0}^{\ell-1} C(n+\ell-1, n+j;\,B)\, (n+1)_j 
    \binom{(A+\ell-1)/B}{j}
    ,
  \end{equation}
  where the coefficients $C(n,k;B)\rsmcoloneq\sum_{i=k}^n s(n,i)B^i S(i,k)$
  are `C\nobreakdash-numbers'~\cite{Charalambides77}.  Closed-form
  expressions for $C(n,k;B)$ when $B=-1,\frac12,2$ are
  known~\cite[p.~158]{Comtet74}.  However, these values of~$B$ play no
  special role in the present analysis.
\end{remarkaftertheorem}

Examples of transformed sequences~$\hat f_\ell(A,B;n)$, $n\ge0$, include
\begin{subequations}
\label{eq:subeqs}
  \begin{align}
  \hat f_{-1}(A,B;\,n) &= \delta_{n,0} + \frac{AB}{A+B-1}\,\delta_{n,1},\\
  \hat f_0(A,B;\,n) &= \delta_{n,0},\\
  \hat f_1(A,B;\,n) &= (-B)^n,\\
  \hat f_2(A,B;\,n) &= \left[(A+1) + \tfrac12(B-1)n \right]\frac{(n+1)(-B)^n}{A+1}.
  \label{eq:subeqsd}
  \end{align}
\end{subequations}
These follow readily from the representation~(\ref{eq:defconvol}), or with
more effort from~(\ref{eq:rep2}).  For later reference, note that
\begin{equation}
\label{eq:inserta}
  \hat f_2(A,B;\;n) + B\,\hat f_2(A,B;\;n-1)  = \left[A+1+(B-1)n\right]\frac{(-B)^n}{A+1}.
\end{equation}

\begin{theorem}
\label{thm:cor1}
For each\/ $\ell\in\mathbb{Z}$, there is a certain rational function of\/
$y$, $F_\ell(A,B;y)$, with coefficients polynomial in\/ $A,B$, such that
\begin{equation}
\label{eq:fund}
y^A\,F_\ell(A,B;\,y)=
  1 + \sum_{k=1}^\infty
  \frac{(A+Bk+1)_{\ell-1}}{(A+1)_{\ell-1}}
  \binom{A+Bk}{k}z^k
\end{equation}
holds in a neighborhood of\/ $z=0$, if\/ $y$ near\/ $1$ is the solution of
the standardized trinomial equation\/ $y -\nobreak 1 -\nobreak z\,y^B=0$.
{\rm(}By assumption, $B\neq0,1$.{\rm)} Specifically,
\begin{equation}
\label{eq:insertform}
  F_\ell(A,B;\,y)=\sum_{n=0}^\infty \hat f_\ell(A,B;\,n)\left(\frac{1-y}{y}\right)^n,
\end{equation}
where the sequence $\hat f_\ell$ was defined in\/ {\rm(\ref{eq:defconvol})}
and characterized in Theorem\/~{\rm\ref{thm:ratchar}}.

The rational function $F_\ell$ has at most one pole on\/ $\mathbb{P}^1_y$.
If\/ $\ell=-m<0$, the pole is at\/ $y=0$ and is of order\/~$m$; and if\/
$\ell>0$, it is at\/ $y=B/(B-\nobreak 1)$ and is of order\/ $2\ell-1$.
If\/ $\ell>0$, $F_\ell$ can be written as\/ $(y^2D_y)^{\ell-1}\tilde
F_\ell$, where\/ $D_y={\rm d}/{\rm d}y$ and\/ $\tilde F_\ell$ is rational
with a pole at\/ $y=B/(B-\nobreak1)$ of order\/~$\ell$.
\end{theorem}

\begin{proof}
  The identity~(\ref{eq:fund}) is the specialization
  of~(\ref{eq:extraeqn2}) to the case $f=f_\ell$, $\hat f=\hat f_\ell$.  By
  Thm.~\ref{thm:ratchar}, $F_\ell$~is a degree\nobreakdash-$m$ polynomial
  in $v\rsmcoloneq (1-\nobreak y)/y$ if $\ell=-m\le0$, and if $\ell>0$, it
  is rational on~$\mathbb{P}^1_v$ with a unique pole (of order
  $2\ell-\nobreak1$) at $v=(-B)^{-1}$, i.e., at $y=B/(B-\nobreak 1)$.
  Since $(n+\nobreak1)_{\ell-1}\divides \hat f_\ell$, the function $F_\ell$
  on~$\mathbb{P}_v^1$ is the $(\ell-\nobreak1)$'th derivative of some
  rational function with a unique finite\nobreakdash-$v$ pole (of
  order~$\ell$) at $v=(-B)^{-1}$.  The final claim thus follows from $D_v =
  -y^2D_y$.
\end{proof}

As examples of the use of formula~(\ref{eq:insertform}), four of the
rational functions~$F_\ell$ are listed in Table~\ref{tab:1}.  They come
from the sequences $\hat f_\ell$ given in~(\ref{eq:subeqs}).

\begin{table}
  \caption{The first few rational functions $F_\ell$, with
  $S\rsmcoloneq\frac12-\ell$.}
\hfil
  \begin{tabular}{rrl}
    \hline
    $\ell$ & $S$ & $F_\ell(A,B;y)$\\
    \hline
    $-1$ & $\frac32$ & $\frac{AB-(A-1)(B-1)y}{(A+B-1)\,y}$ \\
    $0$ & $\frac12$ & $1$ \\
    $1$ & $-\frac12$ & $\frac{y}{(1-B)y+B}$ \\
    $2$ & $-\frac32$ & $\frac{y^2\left[(B-1)(B-A-1)y-B(B-A-2)\right]}{(A+1)\left[(1-B)y+B\right]^3}$ \\
         &  &\quad $=\left({y^2}D_y\right)\left\{
\tfrac{2(B-1)(B-A-1)y - B(2B-2A-3)}{2(A+1)(B-1)\left[(1-B)y + B\right]^2}
\right\}$ \\
    \hline
  \end{tabular}
\label{tab:1}
\end{table}

Theorem~\ref{thm:cor1}, with the aid of Table~\ref{tab:1}, yields the
following binomial coefficient identities, which are indexed by
$\ell=-1,0,1,2$, respectively.  They hold near $z=0$, if $y$ near~$1$ is
defined by $y-\nobreak1-\nobreak z\,y^B=0$.
\begin{subequations}
  \begin{gather}
    y^A\left[\frac{AB-(A-1)(B-1)y}{(A+B-1)\,y}\right] =
    1+\sum_{k=1}^\infty \frac{A-1}{A+Bk-1}\,\frac{A}{A+Bk}\binom{A+Bk}{k}z^k,\label{eq:newident1}\\
    y^A = 1 + \sum_{k=1}^\infty \frac{A}{A+Bk}\binom{A+Bk}{k} z^k,\label{eq:LambertRam}\\
    y^A\left[\frac{y}{(1-B)y+B}\right] = 1 + \sum_{k=1}^\infty \binom{A+Bk}{k} z^k,\label{eq:Polya}\\
y^A\left\{\frac{y^2\left[(B-1)(B-A-1)y-B(B-A-2)\right]}{(A+1)\left[(1-B)y+B\right]^3}\right\} \hskip1.25in\label{eq:newident2}\\
    \hskip2.0in
    =1+\sum_{k=1}^\infty \frac{A+Bk+1}{A+1}\binom{A+Bk}{k}z^k.\notag
  \end{gather}
\end{subequations}
The $\ell=0$ identity~(\ref{eq:LambertRam}) and the $\ell=1$
identity~(\ref{eq:Polya}) are well known.  The former was derived by
Lambert in 1758 (and by Ramanujan~\cite[p.~72]{BerndtI}), and the latter by
P\'olya~\cite{Polya21}.  They can be proved by Lagrange
inversion~\cite[Thm.~2.1]{Merlini2006}, and imply each
other~\cite[Lem.~1]{Zeitlin70}.  But the identities
(\ref{eq:newident1}),(\ref{eq:newident2}) are less familiar.  The
$\ell=3,4,\dotsc$ and $\ell=-2,-3,\dots$ identities can also be worked~out.

Each identity in this family indexed by $\ell\in\mathbb{Z}$ can be obtained
by differentiating the preceding.  This is formalized in the following
recurrence, which is equivalent to one obtained by Chu~\cite{Chu95}.  By
starting with $F_{0}\equiv1$ and iterating, one can generate all $F_\ell$,
$\ell>0$.  By doing the reverse, i.e., by integrating (and exploiting the
fact that $F_\ell(1)=1$, $\forall\ell\in\mathbb{Z}$, which determines each
constant of integration), one can generate all $F_\ell$, $\ell<0$.

\begin{theorem}
\label{thm:iteratable}
The rational functions\/ $F_\ell(A,B;y)$, $\ell\in\mathbb{Z}$, satisfy
  \begin{equation*}
    F_{\ell+1}(A,B;\,y) = \frac{(A-B+2)_{\ell-1}}{(A+1)_{\ell}}\,
\frac{y^{-A+B+1}}{(1-B)y + B} \, D_y\left[y^{A-B+1}F_\ell(A-B+1,B;\,y)\right].
  \end{equation*}
\end{theorem}
\begin{proof}
Apply $D_z$ to both sides of~(\ref{eq:fund}), using $D_z=\frac{y^{2B}}{y^B
  - B(y-1)y^{B-1}}D_y$, which follows from $y -\nobreak 1 -\nobreak
z\,y^B=0$.
\end{proof}

In the binomial coefficient identities of the form~(\ref{eq:fund}) with
$\ell=3,4,\dots$ and $\ell=-2,-3,\dots$, not given here, the left-hand
functions $F_\ell(A,B;y)$ become increasingly complicated.  But for all
$\ell\in\mathbb{Z}$, the right-hand power series in~$z$ has radius of
convergence $\left|(B-\nobreak1)^{B-1}/B^B\right|$ if $B\neq1$, and unity
if ${B=1}$.  This is consistent with the presence (when $\ell>0$, at~least)
of a pole at $y=B/(B-\nobreak 1)$ on the left-hand side, since $y=B/(B-1)$
corresponds to $z=(B-\nobreak1)^{B-1}/B^B$.  As $B\to1$, the pole moves to
$y=\infty$, i.e., to $z=1$.

The following theorem reveals when the evaluation of the rational function
$F_\ell(A,B;y)$, in any identity in this family, may lead to a division by
zero.

\begin{theorem}
  The rational function\/ $F_\ell(A,B;y)$ in\/ {\rm(\ref{eq:fund})} is of
  the form
  \begin{subequations}
  \begin{sizealignat}{\small}{2}
    &\frac{\tilde\Pi_m(A,B;\,y)}{(A+B-m)_m(A+2B-m)_{m-1}\dotsb(A+mB-m)_1\,y^m},&\qquad&\ell=-m\le0{\rm;}\label{eq:nonremovablepoles}\\
    &\frac{y^\ell\,\Pi_{\ell-1}(A,B;\,y)}{(A+1)_{\ell-1}\left[(1-B)y+B\right]^{2\ell-1}},&\qquad&\ell\ge1.\label{eq:removablepoles}
  \end{sizealignat}
  \end{subequations}
where $\tilde\Pi_m$, $\Pi_{\ell-1}$ are polynomials in~$y$, the subscripts
indicating their degrees.
\end{theorem}
\begin{proof}
  Iterate Thm.~\ref{thm:iteratable} toward negative and
  positive~$\ell$.
\end{proof}

It should be mentioned that in the identities indexed by $\ell\ge1$, the
poles in~$A$ that are evident in~(\ref{eq:removablepoles}), located at
$A=-1,\dots,-\ell+1$, are removable.  The denominator factor
$(A+1)_{\ell-1}$ is present on the right as~well as the left side
of~(\ref{eq:fund}), and can simply be cancelled.  But the poles in~$A$
present in~$F_\ell$, $\ell\le-1$, which are evident
in~(\ref{eq:nonremovablepoles}), are less easily removed.

\smallskip
This family of identities can be generalized by modifying the initial
sequence $f_\ell(A,B,k)$ of~(\ref{eq:fdef}) to include one or more free
parameters.  A~pair of such generalizations, involving the functions
${}_2F_1$ and~${}_3F_2$, will prove useful.  Let
\begin{equation}
g_\ell(A,B,C;\,k) \rsmcoloneq  \frac{A+C+\ell}{A+C+Bk+\ell}\, f_{\ell+1}(A,B;\,k),\qquad
k\ge0, 
\end{equation}
so that $g_\ell$ is `interpolating': it reduces to~$f_{\ell}$ if~$C=0$ and
to~$f_{\ell+1}$ as $C\to\infty$.  It follows from~(\ref{eq:333}), much as
in~(\ref{eq:defconvol}), that the $(A,B)$-transforms of~$g_{0},g_{1}$ are
\begin{subequations}
\begin{align}
\hat g_{0}(A,B,C;\,n) &= (-)^n(A+C)\:
\mathop{\Delta_{A,B}^n}\!\left[\binom{A}{n}\left(\frac1{A+C}\right)\right],\\
\hat g_{1}(A,B,C;\,n) &= (-)^n\left(\frac{A+C+1}{A+1}\right)\:
\mathop{\Delta_{A,B}^n}\!\left[\binom{A}{n}\left(\frac{A+1}{A+C+1}\right)\right].
\end{align}
\end{subequations}
By a finite difference computation or an expansion in partial fractions
(cf.~\cite[\S\,6]{Gould56}), these transformed sequences satisfy
\begin{subequations}
\label{eq:hatg0}
\begin{align}
\label{eq:hatg0a}
    \hat g_0(A,B,C;\,n) &= (-)^n \frac{(C)_n}{\left(\tfrac{A+C}B+1\right)_n},\\
    \hat g_{1}(A,B,C;\,n) + B\,\hat g_{1}(A,B,C;\,n-1) &= (-)^n \frac{(C)_n\bigl(\tfrac{A+1}{B-1}+1\bigr)_n}{\left(\tfrac{A+C+1}{B}+1\right)_n\bigl(\tfrac{A+1}{B-1}\bigr)_n},
\label{eq:hatg0b}
\end{align}
\end{subequations}
for each $n\ge0$, resp.\ $n\ge1$, showing by comparison with
(\ref{eq:subeqs}) and~(\ref{eq:inserta}) that $\hat g_0(A,B,C;n)$ reduces
to $\hat f_0(A,B;n)=(-B)^n$ if~$C=0$ and to $\hat f_1(A,B;n)=\delta_{n,0}$
as $C\to\infty$; and that $\hat g_{1}(A,B,C;n)$ reduces to $\hat
f_1(A,B;n)=1$ if~$C=0$, and to $\hat f_{2}(A,B;n)$, given
in~(\ref{eq:subeqsd}), as $C\to\infty$.  To avoid division by zero
in~(\ref{eq:hatg0}), one may assume that $B\neq0,1$ and that $(A+\nobreak
C)/B$~is not a negative integer.

The right sides of (\ref{eq:hatg0a}),(\ref{eq:hatg0b}) have the form of
coefficients of ${}_2F_1$ and~${}_3F_2$ series.  Applying
Thm.~\ref{thm:mainGould} therefore yields a pair of hypergeometric
identities.
\begin{theorem}
\label{thm:g0}
Under the preceding assumptions, one has the \emph{interpolating}
identities
  \begin{align*}
    &y^A \:\,
    {}_2F_1\left(
    \begin{array}{c}
      C,\, 1\\
      \tfrac{A+C}B+1
    \end{array}
    \!\!\biggm|
    \frac{y-1}y
    \right)
    = 1 + \sum_{k=1}^\infty \frac{A+C}{A+C+Bk}\binom{A+Bk}{k} z^k,
    \\
    &y^A\:\, \left[\frac{y}{(1-B)y+B}\right]\,
    {}_3F_2\left(
    \begin{array}{c}
      C,\, \tfrac{A+1}{B-1}+1,\, 1\\
      \tfrac{A+C+1}{B}+1,\, \tfrac{A+1}{B-1}
    \end{array}
    \!\!\biggm|
    \frac{y-1}y
    \right)
    \\
    &\qquad\qquad\qquad\qquad\qquad{}= 1 + \sum_{k=1}^\infty \frac{A+C+1}{A+C+Bk+1}\,\frac{A+Bk+1}{A+1}\binom{A+Bk}{k} z^k,
  \end{align*}
which are valid near\/ $z=0$, if\/ $y$ near\/ $1$ is defined by\/ $y - 1 -
z\,y^B=0$.
\end{theorem}

Of these two identities the first was found by Gould~\cite[\S\,6]{Gould56},
but the more complicated second one is new.  The parametrized functions
multiplied by~$y^A$ on their left sides will be denoted by
$G_\ell(A,B,C;y)$, $\ell=0,1$, respectively.

The first identity reduces if $C=0$ to~(\ref{eq:LambertRam}), the $\ell=0$
identity of Lambert and Ramanujan, and as~$C\to\infty$ to~(\ref{eq:Polya}),
the $\ell=1$ identity of P\'olya.  The function $G_0(A,B,C;y)$ reduces to
$F_0(A,B;y)\equiv1$ and $F_1(A,B;y)=\allowbreak y/[(1-\nobreak B)y+\nobreak
  B]$, respectively.  Similarly, the second identity reduces if $C=0$ to
P\'olya's identity, and as $C\to\infty$ to the $\ell=2$
identity~(\ref{eq:newident2}).  The function $G_{1}(A,B,C;y)$ reduces
respectively to $F_{1}(A,B;y)$ and~$F_2(A,B;y)$.

\section{Algebraic and hypergeometric functions}
\label{sec:alghyper}

In this section certain algebraic functions, namely the solutions of
trinomial equations, are expressed in~terms of ${}_nF_{n-1}$'s that are
solutions of $E_n$'s with monodromy that is irreducible but imprimitive.
In Thm.~\ref{thm:invrep}, which is the main result of this section,
the~${}_nF_{n-1}$'s are expressed in~terms of algebraic functions by
`inverting' the just-mentioned representations.  If the parameter
$\ell\in\mathbb{Z}$ in Thm.~\ref{thm:invrep} is set to zero, the
expressions reduce to those of Birkeland~\cite{Birkeland27}.
Theorem~\ref{thm:invrepg0} is a partial extension of Thm.~\ref{thm:invrep}
that has an additional free parameter, and is less algebraic.

The chief result of~\S\,\ref{sec:binomial} was Theorem~\ref{thm:cor1},
which provides a family of binomial coefficient identities indexed
by~$\ell$.  It was based on the \emph{standardized} trinomial equation
$y-1-z\,y^B=0$, and the solution~$y$ that is near~$1$ for $z$ near~$0$.
The theorem is now restated in~terms of the \emph{general} trinomial
equation
\begin{equation}
\label{eq:trinomial1}
  x^n - g\,x^p - \beta =0.
\end{equation}
Here $n=p+q$, for integers $p,q\ge1$, and $g,\beta\in\mathbb{C}$ with
at~most one of $g,\beta$ equal to zero.  The condition $\gcd(p,q)=1$ will
be added later.

Let the $n$~solutions of~(\ref{eq:trinomial1}), with multiplicity, be
denoted $x_1,\dots,x_n$.  In the limit $\beta\to0$ with fixed $g>0$, one
may choose
\begin{equation}
\label{eq:reduction1}
  x_j=\left\{
  \begin{alignedat}{2}
     &\varepsilon_q^{-(j-1)}g^{1/q},&\qquad& j=1,\dots,q,\\
     &0,&\qquad& j=q+1,\dots,n.
  \end{alignedat}
\right.
\end{equation}
In the limit $g\to0$ with fixed $\beta>0$, one may choose
\begin{equation}
\label{eq:reduction2}
  x_j = \varepsilon_n^{-(j-1)}\beta^{1/n},\qquad j=1,\dots,n.
\end{equation}
Here and below, $\varepsilon_r\rsmcoloneq \exp(2\pi{\rm i}/r)$ signifies a primitive
$r$'th root of unity.

To reduce the first case of~(\ref{eq:trinomial1}), i.e., that of $\beta$
near zero, to $y-1-z\,y^B=0$, let
\begin{alignat}{3}
\label{eq:scalings1}
  y&=g^{-1}\,x^q,&\qquad&z=\varepsilon_q^{(j-1)n}g^{-n/q}\,\beta,&\qquad&B=-p/q.\\
\intertext{To reduce the second case, i.e., that of $g$ near zero, let}
y&=\beta^{-1}\,x^n,&\qquad&z=\varepsilon_n^{(j-1)q}\beta^{-q/n}\,g,&\qquad&B=p/n.\label{eq:scalings2}
\end{alignat}
By undoing these two reductions, and letting $A=\gamma/q$,
resp.\ $A=\gamma/n$, one obtains the following `de-standardized' version of
Thm.~\ref{thm:cor1}, in which the rational functions $F_\ell(A,B;y)$,
$\ell\in\mathbb{Z}$, were defined.

\begin{theorem}
\label{thm:41}
The following hold for\/ $\ell\in\mathbb{Z}$ and\/
$\gamma\in\mathbb{C}$.\hfil\break {\rm(i)} In a neighborhood of\/ $\beta=0$
with fixed\/ $g>0$, and with the trinomial root\/ $x_j$ near\/
$\varepsilon_q^{-(j-1)}g^{1/q}$ defined by\/
{\rm(\ref{eq:trinomial1}),(\ref{eq:reduction1}),}
  \begin{multline*}
    [\varepsilon_q^{(j-1)} g^{-1/q}\,
    x_j]^\gamma\,F_\ell(\gamma/q,-p/q;\,g^{-1}x_j^q)
    \\
    {}= 1 + \sum_{k=1}^\infty
    \frac{(\gamma/q - pk/q + 1)_{\ell-1}}
	 {(\gamma/q + 1)_{\ell-1}}
\binom{\gamma/q-pk/q}{k}
    \left[\varepsilon_q^{(j-1)n}g^{-n/q}\,\beta\right]^k,
  \end{multline*}
  for $j=1,\dots,q$.\hfil\break {\rm(ii)} In a neighborhood of\/ $g=0$ with
  fixed\/ $\beta>0$, and with the trinomial root\/ $x_j$ near\/
  $\varepsilon_n^{-(j-1)}\beta^{1/n}$ defined by\/
  {\rm(\ref{eq:trinomial1}),(\ref{eq:reduction2}),}
  \begin{multline*}
    [\varepsilon_n^{(j-1)} \beta^{-1/n}\,
    x_j]^\gamma\,F_\ell(\gamma/n,p/n;\,\beta^{-1}x_j^n)
    \\
    {}= 1 + \sum_{k=1}^\infty
    \frac{(\gamma/n+pk/n+1)_{\ell-1}}
    {(\gamma/n+1)_{\ell-1}}
    \binom{\gamma/n+pk/n}{k}
    \left[\varepsilon_n^{(j-1)q}\beta^{-q/n}\,g\right]^k,
  \end{multline*}
  for $j=1,\dots,n$.\hfil\break In {\rm(i)} and {\rm(ii),}
$\gamma\in\mathbb{C}$ is constrained so that no~division by zero occurs in
the evaluation of either the rational function\/ $F_\ell$ or the first
factor of the summand.
\end{theorem}

The connection to hypergeometric equations~$E_n$ with imprimitive
monodromy, and their canonical solutions ${}_nF_{n-1}$, can now be made.
Henceforth, let a Riemann sphere~$\mathbb{P}^1_\zeta$ be parametrized by
\begin{equation}
\label{eq:zetadef}
  \zeta \rsmcoloneq  (-)^q \frac{n^n}{p^pq^q}\,\frac{\beta^q}{g^n}.
\end{equation}
The following formulas extend those of
Birkeland~\cite[\S\,3]{Birkeland27}, who considered only the case $\ell=0$.
Since $F_0\equiv1$, his formulas contain no left-hand $F_\ell$ factor.

\begin{theorem}
\label{thm:preinv}
The following hold for\/ $\ell\in\mathbb{Z}$ and\/ $a\in\mathbb{C}$, with\/
$\zeta$ defined by\/ {\rm(\ref{eq:zetadef})}.\hfil\break {\rm(i)} In a
neighborhood of\/ $\beta=0$ with fixed\/ $g>0$, and with the trinomial root\/
$x_j$ near\/ $\varepsilon_q^{-(j-1)}g^{1/q}$ defined by\/
{\rm(\ref{eq:trinomial1}),(\ref{eq:reduction1}),}
\begin{sizealign}
{\small}
 & [\varepsilon_q^{(j-1)} g^{-1/q}\, x_j]^{-qa}\,F_\ell(-a,-p/q;\,g^{-1}x_j^q)  \notag\\
&\qquad=\sum_{\kappa=0}^{q-1} \left(\varepsilon_q^{n}\right)^{(j-1)\kappa} \frac{(-)^\kappa
    (a)_{1-\ell}}{(a+n\kappa/q)_{1-\ell-\kappa}(1)_\kappa}\left[\frac{(-)^qp^pq^q}{n^n}\,\cdot\zeta\,\right]^{\kappa/q}\notag\\
&\qquad\qquad\qquad{}_nF_{n-1}\left(
\begin{array}{l}
\frac{a}n+\frac{\kappa}q,\dots,\frac{a+(n-1)}n+\frac{\kappa}q\\
\frac{a-\ell+1}p+\frac{\kappa}q, \dots, \frac{a-\ell+p}p+\frac{\kappa}q;\,
\frac1q+\frac{\kappa}q,\dots,\widehat{\frac{q-\kappa}q+\frac{\kappa}q},\dots,\frac{q}q+\frac{\kappa}q
\end{array}
\!\!\biggm|\zeta\right),\notag
\end{sizealign}
for\/ $j=1,\dots,q$.\hfil\break {\rm(ii)} In a neighborhood of\/ $g=0$ with
fixed\/ $\beta>0$, and with the trinomial root\/ $x_j$ near\/
$\varepsilon_n^{-(j-1)}\beta^{1/n}$ defined by\/
{\rm(\ref{eq:trinomial1}),(\ref{eq:reduction2}),}
\begin{sizealign}
{\small}
 & [\varepsilon_n^{(j-1)} \beta^{-1/n}\,x_j]^{-na}\,F_\ell(-a,p/n;\,\beta^{-1}x_j^n) \notag\\
&\qquad=\sum_{\kappa=0}^{n-1}\left(\varepsilon_n^{q}\right)^{(j-1)\kappa} \frac{(-)^\kappa
    (a)_{1-\ell}}{(a+q\kappa/n)_{1-\ell-\kappa}(1)_\kappa}\left[\frac{(-)^qp^pq^q}{n^n}\,\cdot\zeta\,\right]^{-\kappa/n}\notag\\
&\qquad\qquad\qquad{}_nF_{n-1}\left(
\begin{array}{l}
\frac{-a+\ell}p+\frac{\kappa}n,\dots,
\frac{-a+\ell+(p-1)}p+\frac{\kappa}n;\,
\frac{a}q+\frac{\kappa}n,\dots,\frac{a+(q-1)}q+\frac{\kappa}n\\
\frac1n+\frac{\kappa}n,\dots,\widehat{\frac{n-\kappa}n+\frac{\kappa}n},\dots,\frac{n}n+\frac{\kappa}n
\end{array}
\!\!\biggm|\zeta^{-1}\right),\notag
\end{sizealign}
for\/ $j=1,\dots,n$.\hfil\break In {\rm(i)} and {\rm(ii),} $a\in\mathbb{C}$
is constrained so that no~division by zero occurs on either side, and\/
{\rm(}in\/ {\rm(i))} so that no~lower parameter of any\/ ${}_nF_{n-1}$ is a
non-positive integer.
\end{theorem}
\begin{remarkaftertheorem}
The quantity in square brackets on the right sides of (i),(ii) equals
$g^{-n}\beta^q$.  It is raised to the power $\kappa/q$, resp.\ power
$-\kappa/n$.  The branch should be chosen so that the result equals
$g^{-\kappa n/q}\beta^\kappa$, resp.\ $\beta^{-\kappa q/n} g^{\kappa}$.
\end{remarkaftertheorem}

\begin{proof}
  Partition the series of Thm.~\ref{thm:41} into residue classes of
  $k\textrm{ mod }q$, resp.\ $k\textrm{ mod }n$.  That~is, let
  $k=\allowbreak\kappa+\nobreak i\cdot q$,
  resp.\ $k=\allowbreak\kappa+\nobreak i\cdot n$, and for each~$\kappa$,
  sum over $i=0,1,\dots$, using the identity
  \begin{equation}
    (t)_{i\cdot r} = r^r \left(\tfrac{t}r\right)_i
    \left(\tfrac{t+1}r\right)_i\dots
    \left(\tfrac{t+r-1}r\right)_i.
  \end{equation}
  Also, let $\gamma=-qa$, resp.\ $\gamma=-na$.  (That~is, $A=-a$.)
\end{proof}

\begin{corollary}
\label{cor:short}
The following hold for\/ $\ell\in\mathbb{Z}$ and\/
$a\in\mathbb{C}$.\hfil\break{\rm(i)} Near\/ $\zeta=0$ on\/
$\mathbb{P}^1_\zeta$, with the trinomial root\/ $x_j$ near\/
$\varepsilon_q^{-(j-1)}g^{1/q}$ defined by\/
{\rm(\ref{eq:trinomial1}),(\ref{eq:reduction1}),} the coefficient\/ $g>0$
being fixed and arbitrary and the coefficient\/ $\beta$ being defined in
terms of\/ $\zeta$ by\/ {\rm(\ref{eq:zetadef}),} the quantities
\begin{equation*}
x_j^{-qa}\,F_\ell(-a,-p/q;\,g^{-1}x_j^q),\qquad j=1,\dots,q,
\end{equation*}
regarded as functions of\/ $\zeta$, lie in the solution space of
the differential equation
\begin{equation*}
  E_n\left(
  \begin{array}{l}
  \tfrac{a}{n},\dots,\tfrac{a+(n-1)}{n}\\
  \tfrac{a-\ell+1}p,\dots,\tfrac{a-\ell+p}p;\,\tfrac1q,\dots,\tfrac{q}q
  \end{array}
  \right)_\zeta.
\end{equation*}
{\rm(ii)} Near\/ $\zeta=\infty$ on\/ $\mathbb{P}^1_\zeta$, with the
trinomial root\/ $x_j$ near\/ $\varepsilon_n^{-(j-1)}\beta^{1/n}$ defined
by\/ {\rm(\ref{eq:trinomial1}),(\ref{eq:reduction2}),} the coefficient\/
$\beta>0$ being fixed and arbitrary and the coefficient\/ $g$ being defined
in terms of\/ $\zeta$ by\/ {\rm(\ref{eq:zetadef}),} the quantities
\begin{equation*}
x_j^{-na}\,F_\ell(-a,p/n;\,\beta^{-1}x_j^n),\qquad j=1,\dots,n,
\end{equation*}
regarded as functions of the reciprocal variable\/
$\tilde\zeta\rsmcoloneq \zeta^{-1}$, lie in the solution space near\/ $\tilde\zeta=0$
of the differential equation
\begin{equation*}
  E_n\left(
  \begin{array}{l}
  \tfrac{-a+\ell}p,\dots,\tfrac{-a+\ell+(p-1)}p;\,\tfrac{a}q,\dots,\tfrac{a+(q-1)}q\\
  \tfrac{1}{n},\dots,\tfrac{n}{n}
  \end{array}
  \right)_{\tilde\zeta}.
\end{equation*}
In\/ {\rm(i),} it is assumed that\/
$a\in\mathbb{C}$ is such that  no~pair of
lower parameters differs by an integer; i.e., that the singular point\/
$\zeta=0$ of the\/ $E_n$ is non-logarithmic.
\end{corollary}
\begin{proof}
Immediate, by Lemma~\ref{lem:2} applied to Thm.~\ref{thm:preinv}.
\end{proof}

If in either $E_n$ of Corollary~\ref{cor:short}, $a$~is chosen so that
no~upper parameter differs by an integer from a lower one, then the
monodromy group of the~$E_n$ will be irreducible; and moreover, if
$\gcd(p,q)=1$ then the~$E_n$ will be of the imprimitive irreducible type
characterized in Thm.~\ref{thm:BH0}.  The latter fact is obvious
if~$\ell=0$, though if $\ell\in\mathbb{Z}\setminus\{0\}$, integer
displacements of parameters are needed.

The representations of Theorem~\ref{thm:preinv} can now be inverted, to
express the ${}_nF_{n-1}$'s satisfying these~$E_n$'s in~terms of the
trinomial roots $x_1,\dots,x_n$.  In part~(i) of the theorem below, only
$x_1,\dots,x_q$ are needed; but in part~(ii), all are needed.

\begin{theorem}
\label{thm:invrep}
If\/ $\gcd(p,q)=1$, the following hold for each\/ $\ell\in\mathbb{Z}$ and\/
$a\in\mathbb{C}$.\hfil\break {\rm(i)} Near\/ $\zeta=0$ on\/
$\mathbb{P}^1_\zeta$, for arbitrary fixed\/ $g>0$ and for\/
$\kappa=0,\dots,q-1$,
\begin{align}
&\zeta^{\kappa/q}\,{}_nF_{n-1}\left(
\begin{array}{l}
\frac{a}n+\frac{\kappa}q,\dots,\frac{a+(n-1)}n+\frac{\kappa}q\\
\frac{a-\ell+1}p+\frac{\kappa}q,\dots,
\frac{a-\ell+p}p+\frac{\kappa}q;\,
\frac1q+\frac{\kappa}q,\dots,\widehat{\frac{q-\kappa}q + \frac{\kappa}q},\dots,\frac{q}q+\frac{\kappa}q
\end{array}
\!\!\biggm|\zeta\right)\notag\\
&\qquad=
\frac{(-)^\kappa(1)_\kappa(a+n\kappa/q)_{1-\ell-\kappa}}{(a)_{1-\ell}}\left[\frac{(-)^qn^n}{p^pq^q}\right]^{\kappa/q} \notag\\
&\quad\qquad\qquad{}\times q^{-1}\sum_{j=1}^q
\left(\varepsilon_q^{-n}\right)^{(j-1)\kappa}
[\varepsilon_q^{(j-1)}g^{-1/q}
\,x_j]^{-qa}\,F_\ell(-a,-p/q;\,g^{-1}x_j^q), \notag
\end{align}
in which the trinomial root\/ $x_j$ near\/ $\varepsilon_q^{-(j-1)}g^{1/q}$,
$j=1,\dots,q$, is defined by\/
{\rm(\ref{eq:trinomial1}),(\ref{eq:reduction1}),} with the coefficient\/
$\beta$ determined by\/ $\zeta$ according to\/
{\rm(\ref{eq:zetadef}).}\hfil\break {\rm(ii)} Near\/ $\zeta=\infty$ on\/
$\mathbb{P}^1_\zeta$, for arbitrary fixed\/ $\beta>0$ and for\/
$\kappa=0,\dots,{n-1}$,
\begin{align}
&\zeta^{-\kappa/n}\,{}_nF_{n-1}\left(
\begin{array}{l}
\frac{-a+\ell}p+\frac{\kappa}n,\dots,
\frac{-a+\ell+(p-1)}p+\frac{\kappa}n;\,
\frac{a}q+\frac{\kappa}n,\dots,\frac{a+(q-1)}q+\frac{\kappa}n\\
\frac1n+\frac{\kappa}n,\dots,\widehat{\frac{n-\kappa}n + \frac{\kappa}n},\dots,\frac{n}n+\frac{\kappa}n
\end{array}
\!\!\biggm|\zeta^{-1}\right)\notag\\
&\qquad= \frac{(-)^\kappa(1)_\kappa(a+q\kappa/n)_{1-\ell-\kappa}}{(a)_{1-\ell}}\left[\frac{(-)^qn^n}{p^pq^q}\right]^{-\kappa/n}\notag\\
&\quad\qquad\qquad{}\times n^{-1}\sum_{j=1}^n 
\left(\varepsilon_n^{-q}\right)^{(j-1)\kappa}
[\varepsilon_n^{(j-1)} \beta^{-1/n}
x_j]^{-na}\,F_\ell(-a,p/n;\,\beta^{-1}\,x_j^n),
\notag
\end{align}
in which the trinomial root\/ $x_j$ near\/
$\varepsilon_n^{-(j-1)}\beta^{1/n}$, $j=1,\dots,n$, is defined by\/
{\rm(\ref{eq:trinomial1}),(\ref{eq:reduction2}),} with the coefficient\/
$g$ determined by\/ $\zeta$ according to\/ {\rm(\ref{eq:zetadef}).}
\hfil\break In {\rm(i)} and {\rm(ii),} $a\in\mathbb{C}$ is constrained so
that no~division by zero occurs, and\/ {\rm(}in\/ {\rm(i))} so that
no~lower parameter of any\/ ${}_nF_{n-1}$, for any\/ $\kappa$, is a
non-positive integer.
\end{theorem}
\begin{remarkaftertheorem}
A convention for choosing branches of fractional powers must be adhered~to,
in interpreting these formulas.  In~(i), there are $q$~choices for each of
\begin{equation*}
  \zeta^{1/q},\qquad \left[\frac{(-)^qn^n}{p^pq^q}\right]^{1/q}\!\!,\qquad
  \beta.
\end{equation*}
(The last is evident from~(\ref{eq:zetadef}).)  Any choices will work, so
long as the formal $q$'th root of Eq.~(\ref{eq:zetadef}), i.e.,
\begin{equation}
  \zeta^{1/q} = \left[\frac{(-)^qn^n}{p^pq^q}\right]^{1/q}\!\cdot\, g^{-n/q}\cdot\beta,
\end{equation}
is satisfied.  This removes one degree of freedom.  Part~(ii) is
similarly interpreted.
\end{remarkaftertheorem}

\begin{proof}
  Each of the two formulas in Thm.~\ref{thm:preinv} has the form of a
  linear transformation expressed as a matrix--vector product, i.e.,
  $u_j=\sum_{\kappa}M_{j\kappa}v_\kappa$.  Here,
  $\mathsf{M}=(M_{j\kappa})=(\varepsilon^{(j-1)\kappa})$ is a $q\times q$,
  resp.\ $n\times n$ matrix, with $\varepsilon\rsmcoloneq \varepsilon_q^{n}$,
  resp.~$\varepsilon_n^{q}$.  If $\gcd(p,q)=1$, or equivalently
  $\gcd(q,n)=1$, then the root of unity $\varepsilon$~will be a
  \emph{primitive} $q$'th, resp.~$n$'th, root, and $\mathsf{M}$~will be
  nonsingular.  The $v_\kappa$ are expressed in~terms of the~$u_j$ by
  inverting~$\mathsf{M}$.  But $\mathsf{M}^{-1}=q^{-1}\mathsf{M}^*$, resp.\
  $\mathsf{M}^{-1}=n^{-1}\mathsf{M}^*$.  Elementwise multiplication
  by~$\mathsf{M}^*$ yields the claimed summation formulas.
\end{proof}

\begin{corollary}
  If\/ $\gcd(p,q)=1$, then the following holds for\/ $a\in\mathbb{C}$ with
  $na\notin\mathbb{Z}$.  Near\/ $\zeta=\infty$ on\/ $\mathbb{P}^1_\zeta$,
  if\/ $x_j$ near\/ $\varepsilon_n^{-(j-1)}\beta^{1/n}$, $j=1,\dots,n$, is
  defined by\/ {\rm(\ref{eq:trinomial1}),(\ref{eq:reduction2}),} with the
  coefficient\/~$g$ determined by\/~$\zeta$ according to\/
  {\rm(\ref{eq:zetadef}),} then the exponentiated roots
  \begin{equation*}
    x_j^{-na},\qquad j=1,\dots,n,
  \end{equation*}
  regarded as functions of\/ $\tilde\zeta\rsmcoloneq \zeta^{-1}$, span the solution
  space near\/ $\tilde\zeta=0$ of the differential equation
  \begin{equation*}
  E_n\left(
  \begin{array}{l}
  \tfrac{-a}p,\dots,\tfrac{-a+(p-1)}p;\,\tfrac{a}q,\dots,\tfrac{a+(q-1)}q\\
  \tfrac{1}{n},\dots,\tfrac{n}{n}
  \end{array}
  \right)_{\tilde\zeta}.
  \end{equation*}
\end{corollary}
\begin{proof}
Immediate by Lemma~\ref{lem:2}, applied to the $\ell=0$ case of
Thm.~\ref{thm:invrep}(ii).  (The condition $na\notin\mathbb{Z}$
ensures linear independence of the exponentiated roots.)
\end{proof}

The representations of certain ${}_nF_{n-1}$'s given in
Thm.~\ref{thm:invrep}, in~terms of the solutions of trinomial equations,
are the point of departure for the rest of this paper.  In~each
${}_nF_{n-1}$, the parametric excess~$S$ (the sum of the lower parameters,
minus the sum of the upper ones) equals $\frac12-\ell$, by examination.  As
was mentioned, the $\ell=0$ case of the theorem, when the $F_\ell$ factor
on each right-hand side degenerates to unity, was previously obtained by
Birkeland.

Many classical hypergeometric identities, such as the cubic transformations
of ${}_3F_2$ found by Bailey, are restricted to ${}_nF_{n-1}$'s with
$S=\frac12$.  Sometimes there are `companion' identities that are satisfied
by hypergeometric functions with $S=-\frac12$.  (For Bailey's identities,
see \cite[(5.3,5.4) and (5.6,5.7)]{Gessel82}.)  It is remarkable that in
the present context, the cases $S=\frac12$ and~$-\frac12$, i.e., $\ell=0$
and~$1$, are merely the most easily treated steps on an infinite ladder of
possibilities.

\smallskip
The preceding theorems relating hypergeometric and algebraic functions,
including Thm.~\ref{thm:invrep}, came ultimately from Thm.~\ref{thm:cor1},
which applied the Vandermonde convolution transform to the
sequences~$f_\ell$ given in Defn.~\ref{def:fdef}.  One could start instead
with Thm.~\ref{thm:g0}, i.e., with the `interpolating' sequences
$g_0,g_{1}$ parametrized by~$C$, which reduce to~$f_0,f_1$ if~$C=0$ and
to~$f_1,f_{2}$ as~$C\to\infty$.  Deriving the following theorem, which
interpolates between the $\ell=0,1$ and $1,2$ cases of
Thm.~\ref{thm:invrep}, is straightforward.  (In~it, $c$~stands for~$-C$.)

\begin{theorem}
\label{thm:invrepg0}
If\/ $\gcd(p,q)=1$, the following hold for\/ $\ell=0,1$ and\/
$a\in\mathbb{C}$, with
\begin{align*}
G_0(A,B,C;\,y)&\rsmcoloneq{}_2F_1\left(
\begin{array}{l}
C,\,1\\
\frac{A+C}B+1  
\end{array}
\!\!\biggm|\frac{y-1}y
\right),\\
G_{1}(A,B,C;\,y)&\rsmcoloneq
\left[\frac{y}{(1-B)y+B}\right]\,
    {}_3F_2\left(
    \begin{array}{c}
      C,\, \tfrac{A+1}{B-1}+1,\, 1\\
      \tfrac{A+C+1}{B}+1,\, \tfrac{A+1}{B-1}
    \end{array}
    \!\!\biggm|
    \frac{y-1}y
    \right).
\end{align*}
\hfil\break {\rm(i)} Near\/ $\zeta=0$ on\/ $\mathbb{P}^1_\zeta$, for
arbitrary fixed\/ $g>0$ and for\/ $\kappa=0,\dots,q-1$,
\begin{sizealign}
{\footnotesize}
&\zeta^{\kappa/q}\,{}_{n+1}F_{n}\left(
\begin{array}{ll}
\frac{a}n+\frac{\kappa}q,\dots,\frac{a+(n-1)}n+\frac{\kappa}q; & \frac{a+c-\ell}p+\frac\kappa{q}\\
\frac{a-\ell}p+\frac{\kappa}q,\dots,
\frac{a-\ell+(p-1)}p+\frac{\kappa}q;\,
\frac1q+\frac{\kappa}q,\dots,\widehat{\frac{q-\kappa}{q}+\frac{\kappa}q},\dots,\frac{q}q+\frac{\kappa}q;&\frac{a+c-\ell+p}{p}+\frac\kappa{q}
\end{array}
\!\!\biggm|\zeta\right)\notag\\
&\qquad=
\frac{(-)^\kappa(1)_\kappa(a+n\kappa/q)_{-\ell-\kappa}(a+c-\ell+p\kappa/q)}{(a)_{-\ell}(a+c-\ell)}\left[\frac{(-)^qn^n}{p^pq^q}\right]^{\kappa/q} \notag\\
&\quad\qquad\qquad{}\times q^{-1}\sum_{j=1}^q
\left(\varepsilon_q^{-n}\right)^{(j-1)\kappa}
[\varepsilon_q^{(j-1)}g^{-1/q}
\,x_j]^{-qa}\,G_\ell(-a,-p/q,-c;\,g^{-1}x_j^q), \notag
\end{sizealign}
in which the trinomial root\/ $x_j$ near\/ $\varepsilon_q^{-(j-1)}g^{1/q}$,
$j=1,\dots,q$, is defined by\/
{\rm(\ref{eq:trinomial1}),(\ref{eq:reduction1}),} with the coefficient\/
$\beta$ determined by\/ $\zeta$ according to\/ {\rm(\ref{eq:zetadef}).}
\hfil\break {\rm(ii)} Near\/ $\zeta=\infty$ on\/ $\mathbb{P}^1_\zeta$, for
arbitrary fixed\/ $\beta>0$ and for\/ $\kappa=0,\dots,n-1$,
\begin{sizealign}
{\footnotesize}
&\zeta^{-\kappa/n}\,{}_{n+1}F_{n}\left(
\begin{array}{ll}
\frac{-a+\ell+1}p+\frac{\kappa}n,\dots,
\frac{-a+\ell+p}p+\frac{\kappa}n;\,
\frac{a}q+\frac{\kappa}n,\dots,\frac{a+(q-1)}q+\frac{\kappa}n; & \frac{-a-c+\ell}p+\frac\kappa{n}\\
\frac1n+\frac{\kappa}n,\dots,\widehat{\frac{n-\kappa}n+\frac\kappa{n}},\dots,\frac{n}n+\frac{\kappa}n;
& \frac{-a-c+\ell+p}p + \frac{\kappa}n
\end{array}
\!\!\biggm|\zeta^{-1}\right)\notag\\
&\qquad= \frac{(-)^\kappa(1)_\kappa(a+q\kappa/n)_{-\ell-\kappa}(a+c-\ell-p\kappa/n)}{(a)_{-\ell}(a+c-\ell)}\left[\frac{(-)^qn^n}{p^pq^q}\right]^{-\kappa/n}\notag\\
&\quad\qquad\qquad{}\times n^{-1}\sum_{j=1}^n 
\left(\varepsilon_n^{-q}\right)^{(j-1)\kappa}
[\varepsilon_n^{(j-1)} \beta^{-1/n}
x_j]^{-na}\,G_\ell(-a,p/n,-c;\,\beta^{-1}\,x_j^n),
\notag
\end{sizealign}
in which the trinomial root\/ $x_j$ near\/
$\varepsilon_n^{-(j-1)}\beta^{1/n}$, $j=1,\dots,n$, is defined by\/
{\rm(\ref{eq:trinomial1}),(\ref{eq:reduction2}),} with the coefficient\/
$g$ determined by\/ $\zeta$ according to\/ {\rm(\ref{eq:zetadef}).}
\hfil\break In {\rm(i)} and {\rm(ii),} $a,c\in\mathbb{C}$ are constrained
so that no~division by zero occurs, and\/ so that no~lower parameter of
any\/ ${}_{n+1}F_{n}$, for any\/ $\kappa$, is a non-positive integer.
\end{theorem}

The $\ell=0,1$ representations of ${}_{n+1}F_{n}$'s given in
Thm.~\ref{thm:invrepg0} reduce if~$c=0$ to the $\ell=0,1$ representations
of ${}_{n}F_{n-1}$'s given in Thm.~\ref{thm:invrep}, by cancelling equal
upper and lower parameters; and in the limit~$c\to\infty$, to the
$\ell=1,2$ ones.

It should be noted that unlike Thm.~\ref{thm:invrep}, which involves
rational right-hand functions $F_\ell$, $\ell\in\mathbb{Z}$,
Thm.~\ref{thm:invrepg0} is essentially non-algebraic: the right-hand
functions $G_\ell$, $\ell=0,1$, which depend on the additional
parameter~$c$, are hypergeometric and are generically transcendental.

\section{Schwarz curves: Generalities}
\label{sec:curves}

This section introduces what will be called Schwarz curves, which are
projective algebraic curves that parametrize ordered $k$-tuples of roots
(with multiplicity) of the general degree\nobreakdash-$n$ trinomial
equation
\begin{equation}
\label{eq:trinomial2}
  x^n - g\,x^p - \beta=0.
\end{equation}
As in~\S\,\ref{sec:alghyper}, it is
assumed that $n=\allowbreak p+\nobreak q$ for relatively prime integers
$p,q\ge1$, and that $g,\beta\in\mathbb{C}$ with at~most one of $g,\beta$
equaling zero.  For each~$k$, any ordered $k$-tuple of roots will trace~out
a Schwarz curve as
\begin{equation}
\label{eq:zetadef2}
  \zeta = (-)^q\frac{n^n}{p^pq^q}\,\frac{\beta^q}{g^n}
\end{equation}
varies.  Any Schwarz curve is a projective curve, since multiplying each
root by any $c\neq0$ will multiply $g,\beta$ by~$c^q,c^n\!$, but leave
$\zeta$ invariant.  The dependence on~$\zeta$, which is really a
projectivized (and normalized) version of the discriminant
of~(\ref{eq:trinomial2}), will be interpreted as specifying a covering
of~$\mathbb{P}^1_\zeta$ by the Schwarz curve, the genus of which will be
calculated.  Determining the dependence of a point on the curve on the base
point~$\zeta$ will be of value when uniformizing the solutions of $E_n$'s
with imprimitive monodromy, since the formulas of
Thm.~\ref{thm:invrep}(i),(ii) express many solutions in~terms of
$q$-tuples, resp.\ $n$-tuples of roots of~(\ref{eq:trinomial2}).

Let $\sigma_l=\sigma_l(x_1,\dots,x_n)$ denote the $l$'th elementary
symmetric polynomial in the roots $x_1,\dots,x_n$, so that $\sigma_0=1$,
$\sigma_1=\sum_{i=1}^n x_i$, etc.  From~(\ref{eq:trinomial2}),
\begin{equation}
\label{eq:basecoin0}
  \beta=(-)^{n-1}\sigma_n,\quad\qquad g=(-)^{q-1}\sigma_q;
\end{equation}
and $\sigma_l=0$ for $l=1,\dots,q-1$ and $q+1,\dots,n-1$.  Parametrizing
tuples of roots, given $\zeta\in\mathbb{P}^1$, means parametrizing them
given the ratio $[{\sigma_n}^q:{\sigma_q}^n]$.

\begin{definition}
\label{def:15}
  The top Schwarz curve $\mathcal{C}_{p,q}^{(n)}\subset\mathbb{P}^{n-1}$ is
  the algebraic curve comprising all points
  $[x_1:\ldots:x_n]\in\mathbb{P}^{n-1}$ such that $x_1,\dots,x_n$ satisfy
  the $n-2$ homogeneous equations
  \begin{equation}
    \label{eq:insyst}
    \sigma_l(x_1,\dots,x_n)=0,\qquad 
    l=1,\dots,q-1\text{ and }q+1,\dots,n-1.
  \end{equation}
  That is, $\mathcal{C}_{p,q}^{(n)}$ comprises all
  $[x_1:\ldots:x_n]\in\mathbb{P}^{n-1}$ such that $x_1,\dots,x_n$ are the
  $n$~roots (with multiplicity) of {some} trinomial equation of the
  form~(\ref{eq:trinomial2}).  The curve $\mathcal{C}_{p,q}^{(n)}$ is
  stable under the action of the symmetric group $\mathfrak{S}_n$ on
  $[x_1:\ldots:x_n]$.  An associated degree\nobreakdash-$n!$ covering map
  $\pi_{p,q}^{(n)}\colon \mathcal{C}_{p,q}^{(n)}\to\mathbb{P}^1_\zeta$ is
  defined by
  \begin{equation}
  \zeta = (-)^n\frac{n^n}{p^pq^q}\,\frac{{\sigma_n}^q}{{\sigma_q}^n}.    
  \end{equation}
\end{definition}

\begin{proposition}[\cite{Kato2003}, Cor.~4.7]
  The curve\/ $\mathcal{C}_{p,q}^{(n)}$ is irreducible.
\end{proposition}

\begin{definition}
  For each $k$ with $n> k\ge2$, the subsidiary Schwarz curve
  $\mathcal{C}_{p,q}^{(k)}\subset\mathbb{P}^{k-1}$, also irreducible, is
  the image of~$\mathcal{C}_{p,q}^{(n)}$ under the map
  $[x_1:\ldots:x_n]\mapsto[x_1:\ldots:x_k]$, which
  on~$\mathcal{C}_{p,q}^{(n)}$ is $(n-k)!$-to-$1$.  A~more concrete
  definition of $\mathcal{C}_{p,q}^{(k)}$ is the following.  It is the
  closure in~$\mathbb{P}^{k-1}$ of the set of all points
  $[x_1:\ldots:x_k]\in\mathbb{P}^{k-1}$ such that $x_1,\dots,x_k$ are
  $k$~of the $n$ roots (with multiplicity) of {some} trinomial equation of
  the form~(\ref{eq:trinomial2}).  The curve $\mathcal{C}_{p,q}^{(k)}$ is
  stable under the action of $\mathfrak{S}_k$ on $[x_1:\ldots:x_k]$.

  For each~$k$ with $n\ge k>2$, let a projection
  $\phi_{p,q}^{(k)}\colon\mathcal{C}_{p,q}^{(k)}\to\mathcal{C}_{p,q}^{(k-1)}$
  be defined by $\phi_{p,q}^{(k)}([x_1:\ldots:x_k]) =
  [x_1:\ldots:x_{k-1}]$, so that
  $\phi_{p,q}^{(k+1)}\circ\dots\circ\phi_{p,q}^{(n)}$ is the projection
  $\mathcal{C}_{p,q}^{(n)}\to\mathcal{C}_{p,q}^{(k)}$.  Since
  $\mathcal{C}_{p,q}^{(n)}\to\mathcal{C}_{p,q}^{(k)}$ is $(n-k)!$-to-$1$
  for each~$k$, it follows that for each~$k$, $\phi_{p,q}^{(k)}$~is
  $(n-k+1)$-to-$1$, i.e., is a degree-$(n-\nobreak k+\nobreak1)$ map.
\end{definition}

\begin{remarkaftertheorem}
  Each subsidiary curve $\mathcal{C}^{(k)}_{p,q}\subset\mathbb{P}^{k-1}$,
  $n>k\ge2$, is the solution set of a system of $k-\nobreak2$ homogeneous
  equations in~$x_1,\dots,x_k$, obtained by eliminating $x_{k+1},\dots,x_n$
  from the system~(\ref{eq:insyst}).  The details of this will be given
  shortly.
\end{remarkaftertheorem}

\begin{remarkaftertheorem}
\label{rem:51}
  As will be explained in~\S\,\ref{subsec:61}, defining
  $\mathcal{C}^{(k)}_{p,q}$ as a closure in~$\mathbb{P}^{k-1}$ appends
  at~most a finite number of points to~it.  Specifically, if $k\le p$, it
  appends each point $[x_1:\ldots:x_k]\in\mathbb{P}^{k-1}$ in which
  $x_1,\dots,x_k$ are distinct $p$'th roots of unity.  There are
  $(p-\nobreak k+\nobreak 1)_{k-1} = \allowbreak (p-\nobreak k+\nobreak
  1)_{k}/p$ such points in~$\mathbb{P}^{k-1}$, none of which comes directly
  (via a $k$-tuple of roots) from a trinomial equation of the
  form~(\ref{eq:trinomial2}).  Rather, each is a limit
  in~$\mathbb{P}^{k-1}$ of points that~do.
\end{remarkaftertheorem}

\begin{remarkaftertheorem}
  The projections $\mathcal{C}_{p,q}^{(n)}\to\mathcal{C}_{p,q}^{(k)}$ and
  $\phi_{p,q}^{(k)}$ are really only \emph{partial} maps, being undefined
  at a finite number of singular points, as is typical of maps between
  algebraic curves.  (The symbol~$\dashrightarrow$ could be used instead
  of~$\to$.)  Specifically, if $x_1=\ldots =x_{k-1}=0$ then
  $\phi_{p,q}^{(k)}([x_1:\ldots:x_k])$ is undefined.  The problems with
  zero tuples, associated to trinomials with~$\beta=0$, will be dealt with
  in~\S\,\ref{sec:genera}, where each curve $\mathcal{C}_{p,q}^{(k)}$ will
  be lifted to a desingularized curve~$\tilde{\mathcal{C}}_{p,q}^{(k)}$.
\end{remarkaftertheorem}

\smallskip
The reason for the term `Schwarz curve' is this.  The ratios of any
$n$~independent solutions $y_1,\dots,y_n$ of an order-$n$ differential
equation on~$\mathbb{P}^1_\zeta$ define a (multivalued) \emph{Schwarz map}
from $\mathbb{P}^1_\zeta$ to~$\mathbb{P}^{n-1}$.  Its image is a curve
in~$\mathbb{P}^{n-1}$, and in some cases the inverse map from the image is
single-valued, i.e., supplies a covering of~$\mathbb{P}_\zeta^1$ by the
curve.  Studying Schwarz maps is a standard way of computing the monodromy
groups of differential equations~\cite{Kato2006,Kato2003}, and goes back to
Schwarz's classical work on~$E_2$, the Gauss hypergeometric equation.  In
the present paper, solutions of $E_n$'s with imprimitive monodromy have
been expressed in~terms of tuples of solutions of trinomial equations,
which are algebraic in~$\zeta$; so a purely algebraic use of the Schwarz
map and curve concepts seems warranted.  The top Schwarz
curve~$\mathcal{C}_{p,q}^{(n)}$ was introduced and used by Kato and
Noumi~\cite{Kato2003}, though not under that name; the subsidiary Schwarz
curves seem not to have been treated or exploited before.

It is evident that as defined,
$\mathcal{C}_{p,q}^{(n)}\cong\mathcal{C}_{p,q}^{(n-1)}\!$ and
$\mathcal{C}_{p,q}^{(2)}\cong\mathbb{P}^1_t$, where birational equivalence
is meant.  Here, $t$~is any homogeneous degree\nobreakdash-$1$ rational
function of~$x_1,x_2$.  Henceforth the choice $t=\allowbreak (x_1+\nobreak
x_2)/\allowbreak(x_1-\nobreak x_2)$ will be made, so $[x_1:x_2]=\allowbreak
[t+\nobreak1:t-\nobreak1]$ will be a uniformization
of~$\mathcal{C}^{(2)}_{p,q}$ by the rational parameter~$t$.  It will also
prove useful to define $\mathcal{C}_{p,q}^{(1)}\rsmcoloneq\mathbb{P}^1$.
This will make it possible to define
$\phi_{p,q}^{(2)}\colon\mathcal{C}^{(2)}_{p,q}\to\mathcal{C}^{(1)}_{p,q}$
and $\phi_{p,q}^{(1)}\colon\mathcal{C}^{(2)}_{p,q}\to\mathbb{P}^1_\zeta$ in
such a way that
\begin{equation}
\pi_{p,q}^{(n)}=\phi_{p,q}^{(1)}\circ\dots\circ\phi_{p,q}^{(n)},
\end{equation}
in the sense of being an equality on a cofinite domain.  (See
(\ref{eq:basercoin}) and~(\ref{eq:basecoin}) below.)  For each~$k$ with $n>
k\ge1$,
\begin{equation}
\pi_{p,q}^{(k)}=\phi_{p,q}^{(1)}\circ\dots\circ\phi_{p,q}^{(k)}
\end{equation}
will then define a subsidiary (partial) map
$\pi_{p,q}^{(k)}\colon\mathcal{C}^{(k)}_{p,q}\to\mathbb{P}_\zeta^1$ of
degree $(n-\nobreak k+\nobreak1)_k$.

One can derive a system of polynomial equations for each
$\mathcal{C}^{(k)}_{p,q}$, $n>k>2$, in~terms of $x_1,\dots,x_k$, by using
resultants to eliminate $x_{k+1},\dots,x_n$.  But one can also eliminate
them by~hand in the following way, which will incidentally indicate how
best to define the curve $\mathcal{C}^{(1)}_{p,q}$.  For any~$k$, $0< k<
n$, let $\bar\sigma_m,\hat\sigma_m$ denote the $m$'th elementary symmetric
polynomial in $x_1,\dots,x_k$, resp.\ $x_{k+1},\dots,x_n$.  Then
\begin{equation}
  \sigma_l = \sum_{m=0}^l\bar\sigma_m\hat\sigma_{l-m},\qquad
  0\le l\le n,
\end{equation}
it being understood that $\bar\sigma_m=0$ if $m\notin\{0,\dots,k\}$, and
similarly, that $\hat\sigma_m=0$ if $m\notin\{0,\dots,n-k\}$, with
$\bar\sigma_0=\hat\sigma_0=1$.  The defining equations of the top
curve~$\mathcal{C}^{(n)}_{p,q}$, by Defn.~\ref{def:15}, are
\begin{equation}
\left\{
\begin{array}{lcl}
1=\sigma_{0} = \sum_{m=0}^n \bar\sigma_m\hat\sigma_{0-m}\,;&& \\[\jot]
0=\sigma_1 = \sum_{m=0}^n \bar\sigma_m\hat\sigma_{1-m}, &\ldots,& 0=\sigma_{q-1} = \sum_{m=0}^n \bar\sigma_m\hat\sigma_{(q-1)-m}\,; \\[\jot]
\sigma_{q} = \sum_{m=0}^n \bar\sigma_m\hat\sigma_{q-m}\,;&& \\[\jot]
0=\sigma_{q+1} = \sum_{m=0}^n \bar\sigma_m\hat\sigma_{(q+1)-m},&\ldots,&0=\sigma_{n-1} = \sum_{m=0}^n \bar\sigma_m\hat\sigma_{(n-1)-m}\,; \\[\jot]
\sigma_{n} = \sum_{m=0}^n \bar\sigma_m\hat\sigma_{n-m}.&&
\end{array}
\right.
\label{eq:eqnarray}
\end{equation}

\begin{lemma}
\label{lem:dual}
$\{\hat\sigma_m\}_{m=0}^{n-k}$ can be expressed in terms of\/
  $x_1,\dots,x_k$ {\rm(}and\/ $\sigma_n${\rm)} by
\begin{displaymath}
\hat\sigma_m = \left\{
\begin{aligned}
& (-)^m\! \sum_{\substack{m_1+\dots+m_k=m\\ (\forall j)\:0\le m_j\le m}} x_1^{m_1}\dotsm x_k^{m_k},\\
& \quad\qquad\qquad\qquad\qquad\qquad\qquad\qquad m=0,\dots,\min(q-1,n-k);\\
& \frac{(-)^{n-k-m}\,\sigma_n}{(x_1\dotsm x_k)^{n-k-m+1}} \sum_{\substack{m_1+\dots+m_k=(k-1)(n-k-m)\\ (\forall j)\:0\le m_j\le n-k-m}} x_1^{m_1}\dotsm x_k^{m_k},\\
& \quad\qquad\qquad\qquad\qquad\qquad\qquad\qquad m=\max(q-k+1,0),\dots,n-k.
\end{aligned}
\right.
\end{displaymath}
Note that\/ $0\le m\le\min(q-1,n-k)$ and\/ $\max(q-k+1,0)\le m\le n-k$, the
$m$-ranges of validity of these two formulas, may overlap.
\end{lemma}
\begin{proof}
  The first formula is proved by induction.  The idea is that one solves
  the zeroth equation in~(\ref{eq:eqnarray}) for~$\hat\sigma_0$, then the
  first equation for~$\hat\sigma_1$, etc.; stopping with the
  $\min({q-1},\allowbreak{n-k})$'th equation, since the next one may
  involve~$\sigma_q$, which is not known.  But since the $l$'th equation,
  for $1\le l\le \min(q-1,n-k)$, says that
  \begin{equation}
    \hat\sigma_l = -(\bar\sigma_1\hat\sigma_{l-1}+\dots+\bar\sigma_k\hat\sigma_{l-k}),
  \end{equation}
  the inductive step amounts to verifying that
  \begin{gather}
    \label{eq:foo}
    P(m) = \bar\sigma_1\,P(m-1) - \bar\sigma_2\,P(m-2) + \dots
    + (-)^{k-1}\bar\sigma_k\,P(m-k),\\
    P(m) := \! \sum_{\substack{m_1+\dots+m_k=m\\ (\forall j)\:0\le m_j\le m}} x_1^{m_1}\dotsm x_k^{m_k}.\notag
  \end{gather}
  This is a well-known identity.

  The second formula in the lemma is dual to the first and is proved by a
  similar induction, `downward.'  One solves the $n$'th equation
  in~(\ref{eq:eqnarray}) for~$\hat\sigma_{n-k}$, then the $(n-1)$'st
  equation for~$\hat\sigma_{n-k-1}$, etc.; stopping with the
  $\max(q-k+1,0)$'th equation, since the next one may involve~$\sigma_q$,
  which is not known.  As with the first formula, the inductive step
  reduces to a verification of Eq.~(\ref{eq:foo}).
\end{proof}

The formulas of the lemma are accompanied by \emph{constraints}.  First,
there are the implicit constraints coming from the equivalence between the
two formulas for~$\hat\sigma_m$, for each~$m$ in the doubly covered range
\begin{equation}
  \max(q-k+1,0)\le m\le \min(q-1,n-k).
\end{equation}
Second, there are the equations $\sigma_l=0$ in~(\ref{eq:eqnarray}),
for each~$l$ in the ranges
\begin{equation}
  \min(n-k+1,q)\le l\le q-1,\qquad
  q+1\le l\le\max(k-1,q),
\end{equation}
which were not exploited in the proof of the lemma.  In~all, there are
$k-\nobreak1$ constraint equations.  They serve (if~$k\ge2$) to do two
things: (i)~they yield an expression for $\sigma_n$ as a rational function
of~$\hat\sigma_1,\dots,\hat\sigma_{n-k}$, and hence as a rational function
of~$x_1,\dots,x_k$, homogeneous of degree~$n$; and, (ii)~they impose
$k-\nobreak2$ homogeneous conditions on~$x_1,\dots,x_k$, which are the
desired defining equations of the curve
$\mathcal{C}^{(k)}_{p,q}\subset\mathbb{P}^{k-1}$.

The formula for~$\sigma_q$ in~(\ref{eq:eqnarray}), as yet unused, yields an
expression for~$\sigma_q$ as a rational function of~$x_1,\dots,x_k$,
homogeneous of degree~$q$.  Therefore,
$\zeta\propto{\sigma_n}^q/{\sigma_q}^n$ is homogeneous of degree zero and
is in the function field of~$\mathcal{C}^{(k)}_{p,q}$, and
\begin{equation}
\mathbb{P}^{k-1}\supset\mathcal{C}^{(k)}_{p,q}\ni[x_1:\ldots:x_k]\mapsto \zeta\in\mathbb{P}^1  
\end{equation}
yields a formula for the degree-$(n-\nobreak k+\nobreak1)_k$ (partial) map
$\pi_{p,q}^{(k)}\colon \mathcal{C}^{(k)}_{p,q}\to\mathbb{P}^1_\zeta$, i.e.,
an explicit formula $\zeta=\zeta(x_1,\dots,x_k)$.

The just-sketched elimination procedure can be carried~out by~hand for the
cases $k=2$ (see Lemma~\ref{lem:54} below) and $k=3$ (see
Thm.~\ref{thm:long3}).  When $k$~increases further, carrying it out by~hand
becomes increasingly difficult.

\begin{lemma}
\label{lem:54}
  Applying the elimination procedure to the case $k=2$ yields the formulas
  \begin{gather*}
    \sigma_n=(-)^n (x_1x_2)^p\,\frac{x_1^q - x_2^q}{x_1^p-x_2^p},
    \qquad\quad
    \sigma_q=(-)^{q-1}\,\frac{x_1^n - x_2^n}{x_1^p - x_2^p},\\
    \zeta=(-)^n\frac{n^n}{p^pq^q}\,\frac{\sigma_n^q}{\sigma_q^n} =
    \frac{n^n}{p^pq^q}(x_1x_2)^{pq}\,\frac{(x_1^p-x_2^p)^p\, (x_1^q-x_2^q)^q}
    {(x_1^n-x_2^n)^n},
  \end{gather*}
the last of which defines a degree-$[n(n-1)]$ covering map\/
$\pi_{p,q}^{(2)}\colon\mathcal{C}^{(2)}_{p,q}\to\mathbb{P}^1_\zeta$.
\end{lemma}

\begin{proof}
  By elementary algebra.
\end{proof}

\begin{remarkaftertheorem}
If one simply defines $\pi_{p,q}^{(k)}\colon
\mathcal{C}^{(k)}_{p,q}\to\mathbb{P}^1_\zeta$ as the composition of the two
rational maps $[x_1:\ldots:x_k]\mapsto [x_1:x_2]$ and~$\pi_{p,q}^{(2)}$,
then it will follow immediately from Lemma~\ref{lem:54} that $\zeta$~is in
the function field of~$\mathcal{C}^{(k)}_{p,q}$, without actually employing
the just-sketched elimination procedure to derive an explicit rational
formula $\zeta=\zeta(x_1,\dots,x_k)$.
\end{remarkaftertheorem}

\smallskip
The case $k=1$ obviously requires special treatment, since one cannot
express $\sigma_n,\sigma_q$ in~terms of $x_1$ alone.  If~$k=1$, the
system~(\ref{eq:eqnarray}) reduces to
\begin{equation}
\left\{
\begin{array}{lcl}
0=\sigma_1 = x_1+\hat\sigma_1, &\ldots,& 0=\sigma_{q-1} = x_1\hat\sigma_{q-2}+\hat\sigma_{q-1}\,; \\[\jot]
\sigma_{q} = x_1\hat\sigma_{q-1}+\hat\sigma_q\,;&& \\[\jot]
0=\sigma_{q+1} = x_1\hat\sigma_q + \hat \sigma_{q+1},&\ldots,&0=\sigma_{n-1} = x_1\hat\sigma_{n-2}+\hat\sigma_{n-1}\,; \\[\jot]
\sigma_{n} = x_1\hat\sigma_{n-1},&&
\end{array}
\right.
\end{equation}
the solution of which can be written as
\begin{equation}
\sigma_n=(-)^{n-q-1} x_1^{n-q}\hat\sigma_q,\qquad\quad
  \sigma_q = \hat\sigma_q + (-)^{q-1} x_1^q .
\end{equation}
This suggests focusing on~$s$, the element of the function field
of~$\mathcal{C}^{(n)}_{p,q}$ defined by
\begin{alignat}{2}
\label{eq:focusedalignat}
  s&\rsmcoloneq(-)^{n-1}\sigma_n / x_1^n,&\qquad\quad&1-s:=(-)^{q-1}\sigma_q/x_1^q\\
   &\hphantom{:}=\beta / x_1^n,       &\qquad\quad&\hphantom{1-s:}\!=g/x_1^q\notag
\end{alignat}
(the two definitions being equivalent), and defining the special $k=1$
Schwarz curve $\mathcal{C}^{(1)}_{p,q}$ to be~$\mathbb{P}^1_s$.  

There is then a degree-$(n-\nobreak1)$ map
$\phi_{p,q}^{(2)}\colon\mathcal{C}^{(2)}_{p,q}\to \mathcal{C}^{(1)}_{p,q}$
given by the map $t=\allowbreak (x_1+\nobreak x_2)/\allowbreak
(x_1-\nobreak x_2)\mapsto s$, since both $\sigma_n,\sigma_q$ are rational
functions of~$x_1,x_2$.  By comparing the expressions
for~$\sigma_n,\sigma_q$ in Lemma~\ref{lem:54} with those for
$\sigma_n/x_1^n,\sigma_q/x_1^q$ in~(\ref{eq:focusedalignat}), one finds
that this degree-$(n-\nobreak1)$ map $\phi^{(2)}_{p,q}$ is given by
\begin{equation}
\label{eq:basercoin}
  \begin{alignedat}{2}
    s&=-\,\frac{x_2^p}{x_1^q}\,\frac{x_1^q-x_2^q}{x_1^p-x_2^p}&&=1-\frac{1}{x_1^q}\,\frac{x_1^n-x_2^n}{x_1^p-x_2^p}\\
     &=-\,\frac{(t-1)^p}{(t+1)^q}\,\frac{(t+1)^q-(t-1)^q}{(t+1)^p-(t-1)^p}&&=1-\frac1{(t+1)^q}\,\frac{(t+1)^n-(t-1)^n}{(t+1)^p-(t-1)^p}.
  \end{alignedat}
\end{equation}
The degree-$(n-1)!$ composition $\phi^{(2)}_{p,q} \circ \dots \circ
\phi^{(n)}_{p,q} \colon \mathcal{C}^{(n)}_{p,q} \to
\mathcal{C}^{(1)}_{p,q}\cong \mathbb{P}^1_s$ is made explicit by the ratios
$z_j=x_j/x_1$, $j=1,\dots,n$, being the solutions~$z$ of
\begin{equation}
\label{eq:znp0}
  z^n - (1-s)z^p - s = 0,
\end{equation}
so that $z_2,\dots,z_n$ are the inverse images (with multiplicity) of~$s$
under the corresponding degree-$(n-\nobreak1)$ rational function $s=s(z)$,
which is
\begin{equation}
  s = -\,\frac{z^p(1-z^q)}{1-z^p}
  = 
-\,\frac{z^p(1+z+\dots+z^{q-1})}{1+z+\dots+z^{p-1}}.
\end{equation}
There is a final degree-$n$ map
$\phi_{p,q}^{(1)}\colon\mathcal{C}^{(1)}_{p,q}\cong\mathbb{P}^1_s\to\mathbb{P}^1_\zeta$,
given by
\begin{equation}
\label{eq:basecoin}
\zeta=(-)^n\frac{n^n}{p^pq^q}\,\frac{{\sigma_n}^q}{{\sigma_q}^n} =
(-)^q\frac{n^n}{p^pq^q}\,\frac{s^q}{(1-s)^n}.
\end{equation}
It completes the sequence of maps leading from the top curve
$\mathcal{C}^{(n)}_{p,q}$ down to~$\mathbb{P}^1_\zeta$.

Pulling everything together yields the following theorem, which
summarizes the results of this section.

\begin{theorem}
\label{thm:55}
  For any pair of relatively prime integers\/ $p,q\ge1$ with\/ $n=p+q$,
there is a sequence of algebraic curves and\/ {\rm(}partial\/{\rm)} maps
\begin{equation}
  \mathcal{C}^{(n)}_{p,q}\stackrel{\phi_{p,q}^{(n)}}{\longrightarrow}\mathcal{C}^{(n-1)}_{p,q}\stackrel{\phi_{p,q}^{(n-1)}}{\longrightarrow}\dots
  \stackrel{\phi_{p,q}^{(3)}}{\longrightarrow}\mathcal{C}^{(2)}_{p,q}\stackrel{\phi_{p,q}^{(2)}}{\longrightarrow}\mathcal{C}^{(1)}_{p,q}\stackrel{\phi_{p,q}^{(1)}}{\longrightarrow}\mathbb{P}^1_\zeta,
\end{equation}
in which\/ $\deg\phi_{p,q}^{(k)}=n-k+1$.  For\/ $k=n,n-1,\dots,2$, the
curve\/ $\mathcal{C}^{(k)}_{p,q}\subset\mathbb{P}^{k-1}$, prior to closure,
comprises all\/ $[x_1:\ldots:x_k]$ in which\/ $x_1,\dots,x_k$ is a nonzero
ordered\/ $k$-tuple of roots\/ {\rm(}with multiplicity\/{\rm)} of some
trinomial equation of the form\/ {\rm(\ref{eq:trinomial2})}.  Each partial
map\/ $\phi_{p,q}^{(k)}$, $n\ge k>2$, takes\/ $[x_1:\ldots:x_k]$, where at
least one of $x_1,\dots,x_{k-1}$ is nonzero, to\/ $[x_1:\ldots:x_{k-1}]$.
One writes\/ $\pi_{p,q}^{(k)}$ for\/
${\phi_{p,q}^{(1)}\circ\dots\circ\phi_{p,q}^{(k)}}$, which is of degree
$(n-k+1)_k$.  Any\/ $[x_1:\ldots:x_k]$ in the domain of\/ $\pi_{p,q}^{(k)}$
is taken by it to\/ $\zeta$, computed from the associated trinomial
equation by Eq.\/~{\rm(\ref{eq:zetadef2}).}

The final two curves\/ $\mathcal{C}^{(k)}_{p,q}$, $k=2,1$, are of genus
zero; i.e., $\mathcal{C}^{(2)}_{p,q}\cong\mathbb{P}^1_t$ and\/
$\mathcal{C}^{(1)}_{p,q}\cong\mathbb{P}^1_s$, where\/ $t:=\allowbreak
(x_1+\nobreak x_2)/\allowbreak(x_1-\nobreak x_2)$ and\/ $s$ are rational
parameters.  The final two maps\/ $\phi_{p,q}^{(2)},\phi_{p,q}^{(1)}$ are
given by
\begin{align*}
s &=  \phi_{p,q}^{(2)}(t) = -\,\frac{(t-1)^p}{(t+1)^q}\,\frac{(t+1)^q-(t-1)^q}{(t+1)^p-(t-1)^p},\\
\zeta &=  \phi_{p,q}^{(1)}(s) = (-)^q\frac{n^n}{p^pq^q}\,\frac{s^q}{(1-s)^n},
\end{align*}
and their composition\/
$\pi_{p,q}^{(2)}=\phi_{p,q}^{(1)}\circ\phi_{p,q}^{(2)}$ by
\begin{sizedisplaymath}{\small}
\zeta = \pi_{p,q}^{(2)}(t) =
\phi_{p,q}^{(1)}\left(\phi_{p,q}^{(2)}(t)\right) =
\frac{n^n}{p^pq^q}(t^2-1)^{pq} \frac{\left[(t+1)^p -
(t-1)^p\right]^p\left[(t+1)^q - (t-1)^q\right]^q}{\left[(t+1)^n -
(t-1)^n\right]^n}.
\end{sizedisplaymath}
\end{theorem}

\smallskip
\section{Schwarz curves: Ramifications and genera}
\label{sec:genera}

The Schwarz curves $\mathcal{C}^{(k)}_{p,q}$ introduced
in~\S\,\ref{sec:curves}, in particular
$\mathcal{C}^{(q)}_{p,q},\mathcal{C}^{(n)}_{p,q}$, will be used
in~\S\,\ref{sec:freeparams} to parametrize the solutions of hypergeometric
differential equations ($E_n$'s) with imprimitive monodromy.  The explicit
formulas of Thm.~\ref{thm:55} will be especially useful.  They exploit the
parametrizations of $\mathcal{C}^{(2)}_{p,q},\mathcal{C}^{(1)}_{p,q}$ by
respective rational (i.e., $\mathbb{P}^1$-valued) parameters $t,s$, which
exist because $\mathcal{C}^{(2)}_{p,q},\mathcal{C}^{(1)}_{p,q}$ are of
genus zero.

The question arises whether `higher'~$\mathcal{C}^{(k)}_{p,q}$, such as the
family of projective plane curves
$\mathcal{C}^{(3)}_{p,q}\subset\mathbb{P}^2$, are ever of genus zero.  In
principle, the genus of each
$\mathcal{C}^{(k)}_{p,q}\subset\mathbb{P}^{k-1}$, $n\ge k>2$, can be
calculated from the Hurwitz formula applied to the map
$\pi_{p,q}^{(k)}\colon\mathcal{C}^{(k)}_{p,q}\to\mathbb{P}^1_\zeta$.  But
some care is needed, since in~general, $\mathcal{C}^{(k)}_{p,q}$~is a
\emph{singular} projective curve, and not being smooth, is not a Riemann
surface.  By resolving singularities, one must first construct a smooth
desingularization
$\tilde{\mathcal{C}}^{(k)}_{p,q}\to\mathcal{C}^{(k)}_{p,q}$.  (Up~to
birational equivalence, $\tilde{\mathcal{C}}^{(k)}_{p,q}$ is unique.)  The
Hurwitz formula can then be applied to the lifted map
$\tilde\pi_{p,q}^{(k)}\colon\tilde{\mathcal{C}}^{(k)}_{p,q}\to\mathbb{P}^1_\zeta$,
which is a degree-$(n-\nobreak k+\nobreak1)_k$ holomorphic map of Riemann
surfaces.

Section \ref{subsec:61} counts the singular points
of~$\mathcal{C}^{(k)}_{p,q}$ and determines their multiplicities.  (See
Table~\ref{tab:2}.)  Section \ref{subsec:62} determines the ramification
structure of~$\tilde\pi_{p,q}^{(k)}$.  Section~\ref{subsec:63} explicitly
parametrizes several plane curves
$\mathcal{C}^{(3)}_{p,q}\subset\mathbb{P}^2$, which according to the
resulting formula for the genus of~$\mathcal{C}^{(k)}_{p,q}$
(Thm.~\ref{thm:genus}) are of genus zero.

\subsection{Desingularization}
\label{subsec:61}

Recall that $\mathcal{C}_{p,q}^{(k)}\subset\mathbb{P}^{k-1}$, $n\ge k\ge2$,
is an algebraic curve including each of the points
$[x_1:\ldots:x_k]\in\mathbb{P}^{k-1}$ in which $x_1,\dots,x_k$ are $k$~of
the $n$~roots (with multiplicity) of some trinomial equation
\begin{equation}
\label{eq:trinomial3}
  x^n - g\,x^p - \beta=0
\end{equation}
(with $n=p+q$, $\gcd(p,q)=1$, and at~most one of $g,\beta$ equaling zero).
Cremona inversion in~$\mathbb{P}^{k-1}\!$, i.e., the substitution
$x_j=1/x_j'$, $j=1,\dots,k$, induces an equivalence
$\mathcal{C}_{p,q}^{(k)}\cong\mathcal{C}_{q,p}^{(k)}$.  It was noted (see
Remark~\ref{rem:51}) that if $k\le p$, $\mathcal{C}^{(k)}_{p,q}$~must
really be defined as a closure in~$\mathbb{P}^{k-1}$, and the taking of the
closure appends a finite number of limit points, which do not come directly
from any trinomial equation of the form~(\ref{eq:trinomial3}).  This will
be elucidated below.

It is a standard fact (see~\S\,\ref{subsec:62}) that if ${\beta\neq0}$,
at~most two roots of~(\ref{eq:trinomial3}) can coincide; and coincidence
occurs if and only if
\begin{equation}
\label{eq:zetadef3}
  \zeta \rsmcoloneq (-)^q \frac{n^n}{p^pq^q}\,\frac{\beta^q}{g^n}
\end{equation}
equals unity.  Also, if a point $[x_1:\nobreak\ldots:\nobreak
  x_k]\in\allowbreak {\mathcal{C}}^{(k)}_{p,q}$ comes directly from an
equation of the form~(\ref{eq:trinomial3}), the coefficients $g,\beta$ are
determined by the point, up to the scaling $(g,\beta)\mapsto
(c^qg,c^n\beta)$ for some $c\neq0$.  (To~see this, consider the system
$x_1^n-\nobreak gx_1^p -\nobreak \beta=0$, $x_2^n-\nobreak gx_2^p -\nobreak
\beta=0$, where by the preceding, $x_1\neq x_2$ can be assumed.  Since
$\gcd(n,p)=1$, if this system has a solution $(g,\beta)$, it has a unique
solution; but of~course $(x_1,x_2)$ can be multiplied by~$c$.) The
formula~(\ref{eq:zetadef3}), which is unaffected by scaling, defines the
map $\pi_{p,q}^{(k)}\colon{\mathcal{C}}^{(k)}_{p,q}\to\mathbb{P}^1_\zeta$.

The points on~$\mathcal{C}^{(k)}_{p,q}$ of special interest here are those
taken by~$\pi_{p,q}^{(k)}$ to $\zeta=0,1,\infty$.  Among points coming from
trinomials, the ``$\zeta=\infty$'' points are ones with $g=0$, and the
``$\zeta=0$'' points are ones with $\beta=0$.  (It will be seen that each
of the limit points present if $k\leq p$ is also a $\zeta=0$ point, by
continuity.)  A~`generic' point on~$\mathcal{C}^{(k)}_{p,q}$ is one that is
taken by~$\pi_{p,q}^{(k)}$ to some $\zeta\not\in\{0,1,\infty\}$.

Which points of $\mathcal{C}^{(k)}_{p,q}\subset \mathbb{P}^{k-1}$ are
non-smooth will now be determined.  First, consider any generic point
$[a]=\allowbreak[a_1:\nobreak\ldots:\nobreak a_k]\in\allowbreak
{\mathcal{C}}^{(k)}_{p,q}$.  That is, consider $a_1,\dots,a_k$, distinct
and nonzero, which are $k$~of the roots (with multiplicity) of a unique
trinomial $x^n-\nobreak g_0x^p-\nobreak \beta_0$ with $g_0,\beta_0\neq0$.
Near~$[a]$, treat $\beta$ as a local parameter
on~${\mathcal{C}}^{(k)}_{p,q}$.  That is, for $\beta$ near~$\beta_0$, let
$[a(\beta)]=\allowbreak [a_1(\beta):\nobreak \ldots:\nobreak
  a_k(\beta)]\in{\mathcal{C}}^{(k)}_{p,q}$ come from the roots of
$x^n-\nobreak g_0x^p-\nobreak \beta$.  For each~$j$,
\begin{equation}
\label{eq:kato1}
  \frac{{\rm d}a_j}{{\rm d}\beta} = \left[\frac{\rm d}{{\rm d}x}
  (x^n-g_0x^p)\bigm|_{x=a_j}\right]^{-1} =
  \left(na_j^{n-1}-pg_0a_j^{p-1}\right)^{-1}\!,
\end{equation}
which (by distinctness of roots) is the reciprocal of a nonzero quantity,
for $\beta$ near~$\beta_0$.  Without loss of generality, assume $[a]$~is in
the affine chart on~$\mathbb{P}^{k-1}$ with coordinate
$(z_1,\dots,z_{k-1})$, where $z_j=x_j/x_k$.  To show that
${\mathcal{C}}^{(k)}_{p,q}$ is smooth at~$[a]$, it suffices to show that
$({\rm d}/{\rm
  d}\beta)(z_1,\dots,z_{k-1})|_{\beta=\beta_0}\neq(0,\dots,0)$, with
$z_j=a_j/a_k$.  Suppose not; in~fact, suppose merely that
\begin{equation}
\label{eq:kato2}
  \frac{\rm d}{{\rm d}\beta}(a_j/a_k) = a_k^{-2}\left(
a_k \frac{{\rm d}a_j}{{\rm d}\beta}
- a_j \frac{{\rm d}a_k}{{\rm d}\beta} \right)
\end{equation}
is zero at $\beta=\beta_0$ for $j$ equal to some~$j^*$, $1\le j^*\le k-1$.
This is equivalent to $a_j^{-1}({\rm d}/{\rm d}\beta)a_j|_{\beta=\beta_0}$
taking equal values when $j=j^*\!,k$, or (by~(\ref{eq:kato1}))
$na_j^n-\nobreak pg_0a_j^p$ taking equal values.  But, equal values are
also taken when $j=j^*\!,k$ by $a_j^n-\nobreak g_0a_j^p=\beta_0$.  Hence
$a_{j_*}^n=a_k^n$ and $a_{j_*}^p=a_k^p$; which since $\gcd(n,p)=1$,
contradicts $a_{j_*}\neq a_k$.

Points $[a]\in{\mathcal{C}}^{(k)}_{p,q}$ coming from trinomials with
$g=g_0=0$ (i.e., ``$\zeta=\infty$'' points) can similarly be shown to be
smooth, by using $g$ rather than~$\beta$ as a local parameter.  Any
point~$[a]$ coming from $x^n-\beta_0=0$ has $a_j=\beta_0^{1/n}
\varepsilon_n^{\pi_j}\!$, where $\pi_1,\dots,\pi_k$ are distinct elements
of $\{0,1,\dots,n-\nobreak1\}$.  For $g$ near $g_0=0$, let
$[a(g)]=[a_1(g):\ldots:a_k(g)]\in{\mathcal{C}}^{(k)}_{p,q}$ come from the
roots of $x^n-gx^p-\beta_0$, and use the affine coordinate
$(z_1,\dots,z_{k-1})$ on~$\mathbb{P}^{k-1}$, where $z_j=x_j/x_k$.  Then
\begin{gather}
  \label{eq:kato3a}
\frac{{\rm d} a_j}{{\rm d} g}\bigm|_{g=g_0=0} = 
\left[
\frac{{\rm d}}{{\rm d}x}(x^q - \beta_0 x^{-p})\bigm|_{x=a_j}
\right]^{-1}
= (n\beta_0)^{-1} a_j^{p+1},\\
\label{eq:kato3b}
\frac{\rm d}{{\rm d}g}(a_j/a_k)\bigm|_{g=g_0=0} = 
a_k^{-2}\left(
a_k \frac{{\rm d}a_j}{{\rm d}g}
- a_j \frac{{\rm d}a_k}{{\rm d}g} \right)\bigm|_{g=g_0=0} = 
(n\beta_0)^{-1}(a_j/a_k)(a_j^p-a_k^p).
\end{gather}
The latter is zero for $j=1,\dots,k-1$ only if $a_1^p,\dots,a_k^p$ are
equal, but they are not.  So, $({\rm d}/{\rm
  d}g)(z_1,\dots,z_{k-1})|_{g=g_0=0}\neq(0,\dots,0)$, and
${\mathcal{C}}^{(k)}_{p,q}$ is smooth at~$[a]$.

The curve ${\mathcal{C}}^{(k)}_{p,q}\subset\mathbb{P}^{k-1}$ can also be
shown to be smooth at any ``$\zeta=1$'' point, i.e., any point coming from
a trinomial with $\beta\neq0$ that has a pair of coincident roots.  Let
$[a]$ be a point on~${\mathcal{C}}^{(k)}_{p,q}$ with $a_{j_1}=a_{j_2}$,
coming from a trinomial $x^n-g_0x^p-\beta_0$ with $\beta_0\neq0$.  On~a
neighborhood of~$[a]$, a~parameter for the curve can be chosen to be
$t=\allowbreak (\beta-\nobreak \beta_0)^{1/2}$, so that
$a_{j_1}(t),a_{j_2}(t)$ are $a_{j_1}(0)\pm\nobreak Ct+\nobreak O(t^2)$ for
some $C\neq0$, and $a_{j}(t)=\allowbreak a_{j}(0)+\nobreak C_jt^2+\nobreak
O(t^3)$ for each $j\neq j_1,j_2$.  The derivative with respect to~$t$ of
any affine coordinate $(z_j)_{j\neq j^*}=\allowbreak (x_j/x_{j^*})_{j\neq
  j^*}$ on~$\mathbb{P}^{k-1}$, where $j^*$~is chosen so that $j^*\neq
j_1,j_2$, is defined and nonzero at~$t=0$.  Hence,
$\mathcal{C}_{p,q}^{(k)}$ is smooth at~$[a]$.

In summary, a point in $\mathcal{C}^{(k)}_{p,q}\subset\mathbb{P}^{k-1}$
coming from a trinomial $x^n-\nobreak gx^p -\nobreak\beta$ can be
non-smooth only if $g\neq0$ and~${\beta=0}$.  That is, it must be one of
the ``$\zeta=0$'' points, which may display high-order coincidences of
roots.  Limit points on~$\mathcal{C}^{(k)}_{p,q}$ coming not directly but
indirectly from trinomials, via the taking of the closure
in~$\mathbb{P}^{k-1}$ (and in~fact, from a $\beta\to0$ limit), will be
considered later.

For ease of understanding, the singular points of~$\mathcal{C}^{(k)}_{p,q}$
will first be determined in the `top' case $k=n$, in which the closure
in~$\mathbb{P}^{k-1}$ need not be taken.  There are $n!/p!q$ points
$\mathfrak{p}=\allowbreak[x_1:\nobreak\ldots:\nobreak x_n]$ in
$\mathcal{C}^{(n)}_{p,q}\subset\mathbb{P}^{n-1}$ that come from trinomials
$x^n-\nobreak gx^p -\nobreak\beta$ with $g\neq0$ and~${\beta=0}$.  One of
them (cf.~(\ref{eq:reduction1})) is~$\mathfrak{p}_0$, defined by
\begin{equation}
\label{eq:reduction1new}
  x_j=\left\{
  \begin{alignedat}{2}
     &\varepsilon_q^{-(j-1)}g^{1/q},&\qquad& j=1,\dots,q,\\
     &0,&\qquad& j=q+1,\dots,n.
  \end{alignedat}
\right.
\end{equation}
Here $g^{1/q}$ signifies one of the $q$'th roots of~$g$, chosen
arbitrarily.  Each point~$\mathfrak{p}$ at which $\mathcal{C}^{(n)}_{p,q}$
can be non-smooth is obtained from~$\mathfrak{p}_0$ by a permutation of
$x_1,\dots,x_n$.

It will now be shown that $\mathfrak{p}_0$~is a \emph{ordinary
  $(p-\nobreak1)!$-fold multiple point} of the
curve~$\mathcal{C}^{(n)}_{p,q}$: there are exactly $(p-\nobreak1)!$
branches (local irreducible components) of~$\mathcal{C}^{(n)}_{p,q}$
at~$\mathfrak{p}_0$, so that the curve is non-smooth at~$\mathfrak{p}_0$ if
and only $p>2$.  Moreover, in a neighborhood of~$\mathfrak{p}_0$, each of
the $(p-\nobreak1)!$ branches at~$\mathfrak{p}_0$ is holomorphically
parametrized by $\xi\rsmcoloneq (-\beta)^{1/p}$.  This is a consequence of
the fact that near $\beta=\beta_0=0$ ($g\neq0$ being held fixed), each of
the $n$~roots of $x^n-\nobreak gx^p-\nobreak \beta$ has a convergent
Puiseux expansion in~$(-\beta)^{1/p}$.

Assume the roots are ordered so that $[x_1:\ldots:x_n]\to\mathfrak{p}_0$ as
$\beta\to0$, i.e., so that for each~$j$, the limit of~$x_j$ agrees
with~(\ref{eq:reduction1new}).  (This assumption fixes the ordering of
$x_1,\dots,x_q$, but leaves unspecified that of $x_{q+1},\dots,x_n$.)  Then
the Puiseux expansions of $x_1,\dots,x_n$ in~$\xi$ will be of the form
\begin{equation}
\label{eq:reduction1newa}
  x_j(\xi)=\left\{
  \begin{alignedat}{2}
     &\varepsilon_q^{-(j-1)}g^{1/q}\left[1+\sum\nolimits_{i=1}^\infty c_{ji}\, \xi^{pi} \right],&\qquad& j=1,\dots,q,\\
     &(\varepsilon_p^{\bar\alpha_j} \xi)\, g^{-1/p}\left[1 + \sum\nolimits_{i=1}^\infty \bar  c_{i}\,(\varepsilon_p^{\bar\alpha_j} \xi)^{qi}\right],&\qquad& j=q+1,\dots,n,
  \end{alignedat}
\right.
\end{equation}
where $(\bar\alpha_{q+1},\dots,\bar\alpha_n)$ is an unspecified permutation
of $(0,1,\dots,p-\nobreak1)$.  Of these two types of expansion, the first
is obtained from $x^n-\nobreak gx^p-\nobreak \beta=0$ by applying the
implicit function theorem to $x^n-\nobreak gx^p+\nobreak \xi^p=0$, and the
second by applying the theorem to $\xi^qy^n -\nobreak gy^p +\nobreak 1 = 0$
(where $y\rsmcoloneq x/\xi$), which is equivalent.

The branch of~$\mathcal{C}^{(n)}_{p,q}$ at~$\mathfrak{p}_0$ that results
from~(\ref{eq:reduction1newa}), locally parametrized by~$\xi$, is affected
by the choice of the permutation $(\bar\alpha_{q+1},\dots,\bar\alpha_n)$.
But cyclically shifting $(\bar\alpha_{q+1},\dots,\bar\alpha_n)$ by
$(1,\dots,1)$ is equivalent to multiplying $\xi$ by~$\varepsilon_p$, which
does not alter the branch.  Hence when counting branches, one must
quotient~out the group~$\mathfrak{C}_p$ of cyclic permutations: there are
not~$p!$ but $(p-\nobreak1)!$ distinct branches, with distinct tangent
lines.


Like~$\mathfrak{p}_0$, each of the $n!/p!q$ points
$\mathfrak{p}\in\mathcal{C}^{(n)}_{p,q}$ coming from a trinomial with
$\beta=0$ is an ordinary $(p-\nobreak1)!$-fold multiple point; therefore
the curve $\mathcal{C}^{(n)}_{p,q}\subset\mathbb{P}^{n-1}$ is singular if
and only if $p>2$.  Each such point~$\mathfrak{p}$,
including~$\mathfrak{p}_0$, lifts to $(p-\nobreak1)!$ points on the
desingularization
$\tilde{\mathcal{C}}^{(n)}_{p,q}\to{\mathcal{C}}^{(n)}_{p,q}$, and each of
the other points on~$\mathcal{C}^{(n)}_{p,q}$ lifts to a single point.
Each point $\tilde{\mathfrak{p}}\in\tilde{\mathcal{C}}^{(n)}_{p,q}$ lifted
from one of the $n!/p!q$ points
${\mathfrak{p}}\in{\mathcal{C}}^{(n)}_{p,q}$ has a neighborhood
in~$\tilde{\mathcal{C}}^{(n)}_{p,q}$ that is biholomorphic to a
neighborhood of the point~$0$ on the complex $\xi$\nobreakdash-line.
That~is, near each of these $(n!/p!q) \times (p-\nobreak1)!=\allowbreak
n!/pq$ lifted points~$\tilde{\mathfrak{p}}$, $\xi$~can be used as a local
parameter on~$\tilde{\mathcal{C}}^{(n)}_{p,q}$.

Composing $\tilde{\mathcal{C}}^{(n)}_{p,q}\to {\mathcal{C}}^{(n)}_{p,q}$
with the map $\pi^{(n)}_{p,q}\colon{\mathcal{C}}^{(n)}_{p,q} \to
\mathbb{P}^{1}_\zeta$, yields the holomorphic lifted map
$\tilde\pi^{(n)}_{p,q}\colon\tilde{\mathcal{C}}^{(n)}_{p,q} \to
\mathbb{P}^{1}_\zeta$, which near any of these
points~$\tilde{\mathfrak{p}}$ is effectively a holomorphic function
$\zeta=\zeta(\xi)$, defined on a neighborhood of ${\xi=0}$.  Since $\zeta$
as defined by~(\ref{eq:zetadef3}) is proportional to~$\beta^q$ and hence
to~$\xi^{pq}$, the composition $\tilde\pi^{(n)}_{p,q}$ takes each of these
points~$\tilde{\mathfrak{p}}$ to $\zeta=0$ with multiplicity~$pq$.

The desingularization of any subsidiary Schwarz curve
$\mathcal{C}^{(k)}_{p,q}$, $n>k\ge 2$, is similar to that of the top curve
$\mathcal{C}^{(n)}_{p,q}$.  With $g\neq0$ fixed, consider any $k$-tuple of
roots $x_1,\dots,x_k$ of the trinomial~(\ref{eq:trinomial3}), and let $\nu$
denote the number of these $k$~roots that tend to nonzero values
as~$\beta\to0$.  Necessarily, $0\le\nu\le k$ with $\nu\le q$ and $k-\nu\le
p$.  These $k$~roots can be expanded as power series in
$\xi\rsmcoloneq(-\beta)^{1/p}$, and up~to a permutation of $x_1,\dots,x_k$,
the expansions will be of the form
\begin{equation}
  \label{eq:expanded}
  x_j(\xi)=\left\{
  \begin{alignedat}{2}
     &\varepsilon_q^{\alpha_j}g^{1/q}\left[1+\sum\nolimits_{i=1}^\infty c_{ji}\, \xi^{pi} \right],&\qquad& j=1,\dots,\nu,\\
     &(\varepsilon_p^{\bar\alpha_j}  \xi)\, g^{-1/p}\left[1 + \sum\nolimits_{i=1}^\infty \bar  c_{i}\,(\varepsilon_p^{\bar\alpha_j}  \xi)^{qi}\right],&\qquad& j=\nu+1,\dots,k,
  \end{alignedat}
\right.
\end{equation}
for certain distinct $\alpha_j\in\{0,\dots,q-1\}$ and certain distinct
$\bar\alpha_j\in\{0,\dots,p-1\}$.  In~(\ref{eq:expanded}), it is the
initial $\nu$~components of $(x_1,\dots,x_k)$ that do not tend to zero
as~$\beta\to0$, and the final $k-\nobreak\nu$ components that~do.

The case $\nu=0$ occurs only when $k\le p$ and is clearly special, since
the limit of $(x_1,\dots,x_k)$ as~$\beta\to0$ in this case is
$(0,\dots,0)$; but the limit of $[x_1:\nobreak\ldots:\nobreak x_k]$ as
$\beta\to0$ exists in~$\mathbb{P}^{k-1}$ as the point
$[\varepsilon_p^{\bar\alpha_1}:\nobreak\ldots:\nobreak\varepsilon_p^{\bar\alpha_k}]$.
This explains Remark~\ref{rem:51}, on the need to take the closure when
defining $\mathcal{C}^{(k)}_{p,q}\subset\mathbb{P}^{k-1}$, if $k\le p$.  On
the fibre of~$\mathcal{C}^{(k)}_{p,q}$ over $\zeta=0$, there are
$(p-\nobreak k+\nobreak 1)_{k-1} =\allowbreak (p-\nobreak k+\nobreak
1)_{k}/p$ distinct limit points $[x_1:\nobreak\ldots:\nobreak x_k]$ of this
$\nu=0$ type, in which $x_1,\dots,x_k$ are distinct $p$'th roots of unity.
The case $\nu=k$ occurs only when $k\le q$, and is also rather special.  On
the fibre over $\zeta=0$, there are $(q-\nobreak k+\nobreak 1)_{k-1}
=\allowbreak (q-\nobreak k+\nobreak 1)_{k}/q$ distinct points of the
$\nu=k$ type, in which $x_1,\dots,x_k$ are distinct $q$'th roots of unity.

\begin{table}
\caption{Branching data for the ordinary multiple points on the algebraic
  curve~$\mathcal{C}^{(k)}_{p,q}\subset\mathbb{P}^{k-1}\!$, $k\ge2$.  They
  are partitioned into types indexed by~$\nu$, with
  $\max(0,k-p)\le\nu\le\min(q,k)$.  If nonzero, $\nu$ counts the~$j$, $1\le
  j\le k$, for which $x_j\neq0$.  There are $N_{p,q}^P(k,\nu)$ points of
  type~$\nu$.}
\begin{center}
{\small
  \begin{tabular}{c||l|l|l}
    \hline
    & $N_{p,q}^P(k,\nu)$ & $N_{p,q}^T(k,\nu)$  & $M_{p,q}(k,\nu)$ \\
    \hline
$\nu=0$ & $(p-k+1)_{k-1}$ & $1$ & $p$\\
$0<\nu<k$ & $\binom{k}{\nu}(q-\nu+1)_{\nu-1}$ &
    $(p-k+\nu+1)_{k-\nu-1}$ & $pq$ \\
$\nu=k$ & $(q-k+1)_{k-1}$ & $1$ & $q$\\
    \hline
  \end{tabular}
}%
\end{center}
\label{tab:2}
\end{table}

A straightforward extension of the preceding treatment of the top curve
$\mathcal{C}^{(n)}_{p,q}$, which in retrospect was a treatment of the case
$(k,\nu)=(n,q)$, yields the branching data in Table~\ref{tab:2}.
$N^P_{p,q}(k,\nu)$ is the number of ordinary multiple points
$\mathfrak{p}\in\mathcal{C}^{(k)}_{p,q}\subset\mathbb{P}^{k-1}$ of the type
indexed by~$\nu$, coming from a trinomial equation with $\beta=0$ (or~in
the $\nu=0$ case, from a $\beta\to0$ limit).  For each such point,
$N^T_{p,q}(k,\nu)$ is the multiplicity: the number of distinct branches
at~$\mathfrak{p}$, or equivalently the number of
points~$\tilde{\mathfrak{p}}$ to which $\mathfrak{p}$~lifts on the
desingularization
$\tilde{\mathcal{C}}^{(k)}_{p,q}\to{\mathcal{C}}^{(k)}_{p,q}$.  For
example, $N^P_{p,q}(n,q)=\allowbreak n!/p!q$ and
$N^T_{p,q}(n,q)=\allowbreak (p-\nobreak1)!$; these are the previously
derived top-curve values.

The formula $\binom{k}{\nu}(q-\nu+1)_{\nu-1}$, given in the table for
$N^P_{p,q}(k,\nu)$ if ${0<\nu<k}$, comes from (i)~choosing $\nu$~of the
roots $x_1,\dots,x_k$ to be the ones with nonzero limits as~$\beta\to0$,
and (ii)~choosing the $\nu$~associated exponents~$\alpha_j$ (i.e., powers
of~$\varepsilon_q$) to be distinct elements of $\{0,1,\dots,q-\nobreak1\}$,
with the group~$\mathfrak{C}_q$ of cyclic permutations quotiented~out.  At
any $\mathfrak{p}\in\mathcal{C}^{(k)}_{p,q}$ specified in this way, there
are $N^T_{p,q}(k,\nu)=\allowbreak (p-\nobreak
k+\nobreak\nu+\nobreak1)_{k-\nu-1}$ branches of~$\mathcal{C}^{(k)}_{p,q}$.
This is because a branch is specified by a choice of $k-\nobreak\nu$
distinct exponents~$\bar\alpha_j$ (i.e., powers of~$\varepsilon_p$) from
$\{0,1,\dots,p-\nobreak1\}$, with the group~$\mathfrak{C}_p$ of cyclic
permutations quotiented~out.

For each $\nu$ with $0<\nu<k$, each of the $N^P_{p,q}(k,\nu)\times
N^T_{p,q}(k,\nu)$ lifted points
$\tilde{\mathfrak{p}}\in\tilde{\mathcal{C}}^{(k)}_{p,q}$ has a neighborhood
that is biholomorphic to a neighborhood of ${\xi=0}$.  The cases $\nu=0$,
resp.~$k$ are special, since the appropriate local parameter for
$[x_1:\nobreak\ldots:\nobreak x_k]\in\mathcal{C}^{(k)}_{p,q}$ is not~$\xi$
but rather $\xi'\rsmcoloneq\xi^q$, resp.\ $\xi''\rsmcoloneq\xi^p$, as is
evident from~(\ref{eq:expanded}).  In both these cases there is only one
branch, i.e., $N^T_{p,q}(k,0)=1$ (each appended limit point is a smooth
point), resp.\ $N^T_{p,q}(k,k) = 1$.

The final quantity $M_{p,q}(k,\nu)$ in the table is the multiplicity with
which each lifted point
$\tilde{\mathfrak{p}}\in\tilde{\mathcal{C}}^{(k)}_{p,q}$ is taken to the
point $\zeta=0$ in~$\mathbb{P}^1_\zeta$ by the lifted
map~$\tilde\pi_{p,q}^{(k)}$, which is locally a holomorphic function
$\zeta=\zeta(\xi)$, or $\zeta=\zeta(\xi')$ resp.\ $\zeta=\zeta(\xi'')$.
This multiplicity is usually~$pq$, as previously seen in the top case, but
in the special case $\nu=0$,~resp.~$k$, it is $p$,~resp.~$q$, because
$\zeta\propto \xi^{pq} = (\xi')^p = (\xi'')^q$.

Note that
\begin{equation}
\label{eq:bincoeffident}
\sum_{\nu=\max(0,k-p)}^{\min(q,k)} N^P_{p,q}(k,\nu)\cdot
N^T_{p,q}(k,\nu)\cdot  M_{p,q}(k,\nu) = (n-k+1)_k,
\end{equation}
since the left side equals the number of points (with multiplicity) on the
fibre of~$\tilde\pi_{p,q}^{(k)}$ above $\zeta=0$, i.e.,
$\deg\tilde\pi_{p,q}^{(k)}=(n-k+1)_k$.  By substituting the data of
Table~\ref{tab:2} into~(\ref{eq:bincoeffident}), one obtains a familiar
binomial coefficient identity.

\subsection{Genera}
\label{subsec:62}

It is now possible to determine, for each $k\ge2$, the ramifications of the
degree-$(n-\nobreak k+\nobreak1)_k$ holomorphic map of Riemann surfaces
$\tilde\pi_{p,q}^{(k)}\colon\tilde{\mathcal{C}}^{(k)}_{p,q}\to\mathbb{P}^1_\zeta$,
where $\tilde{\mathcal{C}}^{(k)}_{p,q}\to{\mathcal{C}}^{(k)}_{p,q}$ is the
desingularization.  Any
$\tilde{\mathfrak{p}}\in\tilde{\mathcal{C}}^{(k)}_{p,q}$ can be mapped with
nontrivial multiplicity to~$\mathbb{P}^1_\zeta$ only if (i)~its image, as
in~\S\,\ref{subsec:61}, is some point
$\mathfrak{p}\in{\mathcal{C}}^{(k)}_{p,q}$ on the fibre
of~$\pi_{p,q}^{(k)}$ over $\zeta=0$; \emph{or}, (ii)~its
image~$\mathfrak{p}$ is taken with nontrivial multiplicity
by~$\pi^{(k)}_{p,q}$ to some point $\zeta\neq0$.  Case~(ii) can occur only
if $\zeta=\pi^{(k)}_{p,q}(\mathfrak{p})$ is a critical value
of~$\pi^{(k)}_{p,q}$, i.e., only if $\pi^{(k)}_{p,q}{}^{-1}\zeta$ comprises
fewer than $(n-\nobreak k+\nobreak 1)_k$ points
on~${\mathcal{C}}^{(k)}_{p,q}$; i.e., only if $\pi^{(n)}_{p,q}{}^{-1}\zeta$
comprises fewer than $n!$~points on~${\mathcal{C}}^{(n)}_{p,q}$.

Case~(ii) can accordingly occur only if two of the $n$~roots (with
multiplicity) of the trinomial $x^n-\nobreak gx^p-\nobreak \beta$ giving
rise to~$\mathfrak{p}$ coincide; \emph{or}, if the roots are distinct but
are proportional to $\{\varepsilon_n^m\}_{m=0}^{n-1}$, the $n$'th roots of
unity, so that there are only $(n-1)!$ distinct ordered $n$-tuples
$[x_1:\ldots:x_n]\in\mathbb{P}^{n-1}$.  By evaluating the compact and
elegant formula for the discriminant of the
trinomial~\cite{Greenfield84,Lefton82}, one finds that when $\beta\neq0$,
the first of these two subcases occurs only if
\begin{equation}
\label{eq:karrin}
(pg/n)^n - (-p\beta/q)^q = 0,
\end{equation}
which by~(\ref{eq:zetadef3}) is equivalent to $\zeta=1$.  In this subcase,
exactly two of the $n$~roots coincide.  (For a description of the monodromy
of the roots around $\zeta=1$, see~\cite[\S\,2]{Mikhalkin2006}.)  The
second subcase occurs only if $g=0$, i.e., $\zeta=\infty$.  One concludes
that
$\tilde\pi_{p,q}^{(k)}\colon\tilde{\mathcal{C}}^{(k)}_{p,q}\to\mathbb{P}^1_\zeta$
can be ramified only over the three points $\zeta=0,1,\infty$.  It is a
so-called \emph{Bely\u\i\ cover}.

If $\zeta=1$, i.e., Eq.~(\ref{eq:karrin}) holds, without loss of generality
one can take $g=n/p$, $\beta=-q/p$, so that the trinomial is proportional
to $px^n-\nobreak nx^p+\nobreak q$, with a single doubled root at $x=1$.
There are $n!/2$ distinct points
$[x_1:\ldots:x_n]\in\mathcal{C}^{(n)}_{p,q}$ on the fibre over~$\zeta=1$:
namely, points in which $x_1,\dots,x_n$ are permutations of the roots of
$px^n-\nobreak nx^p+\nobreak q$, including the single doubled root.

\begin{definition}
For each $p,q\ge1$ with $\gcd(p,q)=1$, a degree-$(p+\nobreak q-\nobreak2)$
polynomial with simple roots $T_{p,q}(x) =\allowbreak \sum_{j=0}^{p+q-2}
t_j\,x^j$, satisfying
\begin{equation*}
p\,x^{p+q} - (p+q)x^p+q = (x-1)^2\, T_{p,q},
\end{equation*}
is defined by
\begin{equation*}
t_j = \left\{
\begin{alignedat}{2}
&(j+1)q, &\qquad& 0\le j \le p-1,\\
&(p+q-1-j)p, &\qquad& p-1\le j \le p+q-2.
\end{alignedat}
\right.
\end{equation*}
Its roots will be denoted $x^*_{p,q;\alpha}$, $1\le \alpha\le p+q-2$.
\end{definition}

\begin{theorem}
\label{thm:pregenus}
  For each\/ $k$, $n\ge k\ge2$, the degree-$(n-k+1)_k$ map of Riemann
surfaces\/
$\tilde\pi_{p,q}^{(k)}\colon\tilde{\mathcal{C}}^{(k)}_{p,q}\to\mathbb{P}^1_\zeta$
is a Bely\u\i\ cover: it is ramified only over\/ $\zeta=0,1,\infty$.
\begin{enumerate}
\item The fibre over\/ $\zeta=0$ contains\/ $(p-k+1)_{k-1}$,
resp.\ $(q-k+1)_{k-1}$ points of multiplicity\/
$p$, resp.\ $q$, all other points being of multiplicity\/ $pq$.
\item The fibre over\/ $\zeta=1$ contains\/ $(n-k-1)_k$ points of unit
  multiplicity, all other points being of multiplicity\/ $2$.
\item The fibre over\/ $\zeta=\infty$ contains\/ $(n-k+1)_{k-1}$ points of
  multiplicity\/ $n$.
\end{enumerate}
\end{theorem}

\begin{remarkaftertheorem}
  The case $k=n$ of this theorem, dealing with the top Schwarz curve
${\mathcal{C}}^{(n)}_{p,q}$, was proved by Kato and Noumi~\cite{Kato2003}.
Note that the fibre over $\zeta=0$ contains no~points of multiplicity~$p$
if $k>p$, and none of multiplicity~$q$ if $k>q$; and that the fibre over
$\zeta=1$ contains no~points of unit multiplicity if $k\ge n-1$.
\end{remarkaftertheorem}

\begin{proof}
  The facts about the fibre over $\zeta=0$ can be read~off from
  Table~\ref{tab:2}.  

  Above $\zeta=1$, the points of ${\mathcal{C}}^{(k)}_{p,q}$ that occur
  with unit multiplicity are the points
  $[x^*_{p,q;\chi(1)}:\ldots:x^*_{p,q;\chi(k)}]$, where the roots
  $x^*_{p,q;\alpha}$, $1\le\alpha\le n-2$ were defined above, and
  $\chi(1),\dots,\chi(k)$ are distinct integers selected from
  $1,\dots,n-2$.  The points of~${\mathcal{C}}^{(k)}_{p,q}$ that occur with
  multiplicity~$2$ are those in which, instead, one or two of the~$x_j$'s
  are equal to the double root~$1$.

  As was already remarked, the points of~${\mathcal{C}}^{(k)}_{p,q}$ above
  $\zeta=\infty$ are the points $[x_1:\ldots:x_k]\in\mathbb{P}^{k-1}$ in
  which $x_1,\dots,x_k$ are distinct $n$'th roots of unity.  There are
  exactly $(n-\nobreak k+\nobreak1)_{k-1}=\allowbreak (n-\nobreak
  k+\nobreak 1)_{k}/n$ such points.
\end{proof}

\begin{theorem}
\label{thm:genus}
  For each\/ $n\ge k\ge 1$, the genus of\/ $\mathcal{C}^{(k)}_{p,q}\subset
  \mathbb{P}^{k-1}$ as an algebraic curve, and the topological genus of the
  Riemann surface\/ $\tilde{\mathcal{C}}^{(k)}_{p,q}$, are equal to
  \begin{sizedisplaymath}{\small}
    1+\left[\frac{(k-1)(2n-k-2)}{4(n-1)} - \frac{n}{2pq}\right]
    (n-k+1)_{k-1}
    -\frac{q-1}{2q}\,    (p-k+1)_{k-1}
    -\frac{p-1}{2p}\,    (q-k+1)_{k-1}.
  \end{sizedisplaymath}
\end{theorem}

\begin{proof}
  Apply the Hurwitz genus formula to the data given in
  Thm.~\ref{thm:pregenus} (which trivially extends to $k=1$).  The genus is
  stable under~$p\leftrightarrow q$, since
  $\mathcal{C}_{p,q}^{(k)}\cong\mathcal{C}_{q,p}^{(k)}$, even though the
  singular points (i.e., their number and type) are not.
\end{proof}

\begin{corollary}
\label{cor:64}
  {\rm(i)} The following Schwarz curves, and only these, are rational,
  i.e., of genus zero.
  \begin{itemize}
  \item The trivial curves\/ $\mathcal{C}^{(1)}_{p,q}\cong\mathbb{P}^1_s$ and\/
  $\mathcal{C}^{(2)}_{p,q}\cong\mathbb{P}^1_t$, for all coprime\/ $p,q\ge1$.
  \item $\mathcal{C}^{(3)}_{p,q}$ for\/ $\{p,q\}=\{1,2\}$ and\/ $\{1,3\}$.
  \item $\mathcal{C}^{(4)}_{p,q}$ for\/ $\{p,q\}=\{1,3\}$.
  \end{itemize}
  In consequence, $\mathcal{C}^{(q)}_{p,q}$ can be rationally parametrized
  only if\/ $q=1$, $q=2$, or\/ $(p,q)=(1,3)${\rm;} and the top curve\/
  $\mathcal{C}^{(n)}_{p,q}$ only if\/ $\{p,q\}=\{1,1\}$, $\{1,2\}$,
  or\/ $\{1,3\}$.  \hfil\break {\rm(ii)} The following Schwarz curves, and
  only these, are elliptic, i.e., of genus~$1$.
  \begin{itemize}
  \item $\mathcal{C}^{(3)}_{p,q}$ for\/ $\{p,q\}=\{1,4\}$.
  \end{itemize}
\end{corollary}
\begin{proof}
  By Thm.~\ref{thm:genus}, these curves are of genus~$0$ and genus~$1$ as
  claimed.  The `only these' statements remain to be proved.

  It follows from Thm.~\ref{thm:pregenus} that for each~$k$ with $n>k>2$,
  $\tilde{\mathcal{C}}^{(k)}_{p,q} \to \tilde{\mathcal{C}}^{(k-1)}_{p,q}$
  is a covering with nontrivial branching.  By the Hurwitz formula, if
  $g(\tilde{\mathcal{C}}^{(k-1)}_{p,q})>0$ then
  $g(\tilde{\mathcal{C}}^{(k)}_{p,q})$ must be strictly greater than
  $g(\tilde{\mathcal{C}}^{(k-1)}_{p,q})$.  Therefore, to determine which
  Schwarz curves with $k\ge3$ are of genus~$0$ or genus~$1$, one should
  focus on the case $k=3$.  Substituting $k=3$ and $n=\allowbreak
  p+\nobreak q$ into Thm.~\ref{thm:genus} yields
  \begin{equation}
    \label{eq:planecurvegenus}
    g(\tilde{\mathcal{C}}^{(3)}_{p,q}) = \left[(p^2+4\,pq+q^2)-9(p+q)+14\right]/\,2.
  \end{equation}
  By examination, this equals~$0$ only if $\{p,q\}=\{1,2\}$ or~$\{1,3\}$,
  and equals~$1$ only if $\{p,q\}=\{1,4\}$.  Moreover, by
  Thm.~\ref{thm:genus}, $g(\tilde{\mathcal{C}}^{(4)}_{1,3})=0$ and
  $g(\tilde{\mathcal{C}}^{(4)}_{1,4})>1$; so the `only these' statements
  are proved.
\end{proof}

\subsection{Projective plane curves}
\label{subsec:63}

The curves $\mathcal{C}^{(3)}_{p,q}$ are projective plane curves and
include the first Schwarz curves of positive genus.
The following theorem makes each
$\mathcal{C}^{(3)}_{p,q}\subset\mathbb{P}^2$ and the associated covering
map $\pi_{p,q}^{(3)}\colon\mathcal{C}_{p,q}^{(3)}\to\mathbb{P}^1_\zeta$
quite concrete.

\begin{theorem}
\label{thm:long3}
  For all coprime\/ $p,q\ge1$ with at least one of\/ $p,q$ greater than\/ $1$
  and\/ $n\rsmcoloneq p+q$, the curve\/
  $\mathcal{C}^{(3)}_{p,q}\subset\mathbb{P}^2$ has defining equation
  \begin{equation}
    \label{eq:defining3}
    \frac{x_1^p(x_2^n-x_3^n) + x_2^p(x_3^n-x_1^n) + x_3^p(x_1^n-x_2^n)}
    {(x_1-x_2)(x_2-x_3)(x_3-x_1)}=0,
  \end{equation}
  where the left side is a symmetric homogeneous polynomial in\/
  $x_1,x_2,x_3$ that is of degree\/ $n+p-3$ and is of degree\/ $n-2$ in any
  single variable.  The curve\/ $\mathcal{C}^{(3)}_{p,q}$ goes through the
  fundamental points of\/ $\mathbb{P}^2$ {\rm(}i.e.,
  $[1:0:0],\allowbreak[0:1:0],\allowbreak[0:0:1]${\rm)} if and only if\/
  $p\ge2$, and is singular if and only if\/ $p\ge3$, in which case each
  fundamental point is an ordinary\/ $(p-1)$-fold multiple point, and
  there are no other singular points.  The\/ {\rm(}partial\,{\rm)} covering
  map\/ $\pi_{p,q}^{(3)}\colon\mathcal{C}_{p,q}^{(3)}\to\mathbb{P}^1_\zeta$
  is given by
  \begin{align*}
    \sigma_n &= \left\{
    \begin{aligned}
      &(-)^{n-1}\,x_1^px_2^px_3^p \,
      \left[
      \frac{x_1(x_2^{q+1}-x_3^{q+1}) + {\rm{cycl.}}}
	   {x_1(x_2^px_3^{p+1}-x_2^{p+1}x_3^p) + {\rm{cycl.}}}
	   \right]
	   ,
      &\qquad&p>1,\\
      &(-)^{n-1}\,x_1^px_2^px_3^p \,
      \left[
      \frac{x_1^q(x_2-x_3) + {\rm{cycl.}}}
	   {x_1^{p+1}(x_2^{p}-x_3^{p}) + {\rm{cycl.}}}
	   \right]
	   ,
      &\qquad& q>1;
    \end{aligned}
    \right.
    \\
    \sigma_q &= \left\{
    \begin{aligned}
      &(-)^{q-1}\,
      \left[
      \frac{x_1(x_2^px_3^{n+1}-x_2^{n+1}x_3^p) + {\rm{cycl.}}}
	   {x_1(x_2^px_3^{p+1}-x_2^{p+1}x_3^p) + {\rm{cycl.}}}
      \right]
      ,
      &\qquad&p>1,\\
      &(-)^{q-1}\,
      \left[
      \frac{x_1^{p+1}(x_2^n-x_3^n)+{\rm{cycl.}}}
	   {x_1^{p+1}(x_2^{p}-x_3^{p})+{\rm{cycl.}}}
      \right]
      ,
      &\qquad& q>1;
    \end{aligned}
    \right.
    \\
    \zeta&=(-)^n\frac{n^n}{p^pq^q}\,\frac{{\sigma_n}^q}{{\sigma_q}^n};
  \end{align*}
  and is of degree\/ $n(n-1)(n-2)$.
\end{theorem}

\begin{proof}
  To get~(\ref{eq:defining3}), substitute the formulas for
  $\beta=(-)^{n-1}\sigma_n$, $g=(-)^{q-1}\sigma_q$ in Lemma~\ref{lem:54}
  into the trinomial equation~(\ref{eq:trinomial3}), and relabel $x$
  as~$x_3$.  To~get the other formulas, apply the elimination procedure
  sketched in~\S\,\ref{sec:curves}.  (Cf.\ Lemma~\ref{lem:54} itself, which
  results from applying the procedure when $k=2$; the details when $k=3$
  are long and are omitted.)  The singular points ($0$~or~$3$ in number)
  come from Table~\ref{tab:2}.  If $p\ge3$, there are singular points on
  the fibre over $\zeta=0$: each fundamental point is a ``$\nu=1$'' one,
  but there are none of types $\nu=0,2,3$.
\end{proof}

\begin{remarkaftertheorem}
  The curves
  $\mathcal{C}^{(3)}_{p,q},\mathcal{C}^{(3)}_{q,p}\subset\mathbb{P}^2$ are
  transforms of each other: the defining equation~(\ref{eq:defining3})
  of~$\mathcal{C}^{(3)}_{p,q}$ is taken by a Cremona inversion
  in~$\mathbb{P}^2$ (i.e., $x_i=1/x_i'$) to itself, altered by
  $p\leftrightarrow q$; as is the formula for $\zeta=\zeta(x_1,x_2,x_3)$.

  The genus of $\mathcal{C}^{(3)}_{p,q}$ was given
  in~(\ref{eq:planecurvegenus}), but comes equally well by substituting the
  data of Thm.~\ref{thm:long3} into the standard genus--degree formula
  \begin{equation}
    g = {\binom{N-1}{2}} - \sum_i {\binom{r_i}2}.
  \end{equation}
  Here $N=n+p-3 = 2p+q-3$ is the degree of $\mathcal{C}^{(3)}_{p,q}$, and
  the sum is over its three singular points, each (if~$p\ge3$) of
  multiplicity $r=p-\nobreak1$.
\end{remarkaftertheorem}

Parametrizing the plane curve $\mathcal{C}^{(3)}_{p,q}$ is a key step
leading to a hypergeometric identity; especially, if either $q$ or
$n=\allowbreak p+\nobreak q$ equals~$3$.  The following examples of
parametrizations are illustrative, and will be exploited
in~\S\,\ref{sec:freeparams}.

\begin{example}
  $\{p,q\}=\{1,2\}$, the simplest case.  The top curves
  $\mathcal{C}^{(3)}_{1,2},\mathcal{C}^{(3)}_{2,1}\subset\mathbb{P}^2$ are
  of genus zero (see Cor.~\ref{cor:64}).  The defining equations are
  $x_1+\nobreak x_2+\nobreak x_3=\nobreak 0$ and $x_1x_2+\nobreak
  x_2x_3+x_3x_1=\nobreak 0$.  They are respectively a line, and a conic
  through the fundamental points of~$\mathbb{P}^2$; and are related by a
  Cremona inversion.

  In~general, $\mathcal{C}_{p,q}^{(n)}\cong\mathcal{C}_{p,q}^{(n-1)}\!$,
  which is reflected in the fact that
  $\phi_{p,q}^{(n)}\colon\mathcal{C}^{(n)}_{p,q}\to\mathcal{C}^{(n-1)}_{p,q}$
  is of degree~$1$.  Hence $t\rsmcoloneq\allowbreak (x_1+\nobreak
  x_2)/\allowbreak (x_1-\nobreak x_2)$, used here to uniformize any
  $\mathcal{C}_{p,q}^{(2)}$, can be used as a rational parameter for
  $\mathcal{C}^{(3)}_{1,2},\mathcal{C}^{(3)}_{2,1}$, yielding
  \begin{equation}
    [x_1:x_2:x_3] = \left\{
      \begin{aligned}[0pt]
        &[t+1:t-1:-2t]{},&\qquad&(p,q)=(1,2);\\
        &[-2t(t+1):-2t(t-1):t^2-1],&\qquad&(p,q)=(2,1).
      \end{aligned}
  \right.
  \end{equation}
  But to respect the $\mathfrak{S}_3$ symmetry it is better to use an
  alternative parameter~$\tilde t$ (related to~$t$ by a M\"obius
  transformation), thus:
\begin{sizeequation}{\small}
\label{eq:68}
  [x_1:x_2:x_3] = \left\{
\begin{aligned}
  &[\omega(1+\tilde t) : \bar\omega(1+\omega \tilde t) : (1+\bar\omega\tilde t)],&\quad&(p,q)=(1,2);\\
  &[\omega(1+\omega\tilde t)(1+\bar\omega\tilde t) : \bar\omega(1+\tilde
     t)(1+\omega\tilde t) : (1+\bar\omega\tilde t)(1+\tilde t)],&\quad&(p,q)=(2,1),
\end{aligned}
\right.
\end{sizeequation}
where $\omega\rsmcoloneq\varepsilon_3$ and
$\bar\omega\rsmcoloneq{\varepsilon_3}^2$. 

According to Thm.~\ref{thm:long3}, the covering
$\pi^{(3)}_{1,2}\colon\mathcal{C}^{(3)}_{1,2}\to\mathbb{P}^1_\zeta$
resp.\ $\pi^{(3)}_{2,1}\colon\mathcal{C}^{(3)}_{2,1}\to\mathbb{P}^1_\zeta$
is performed by a function $\zeta=\zeta(x_1,x_2,x_3)$, namely
\begin{equation}
  \zeta = \left\{
  \begin{aligned}
    & -\,\frac{27}4\, \frac{(x_1x_2x_3)^2}{(x_1x_2+x_2x_3+x_3x_1)^3},  &\qquad&(p,q)=(1,2);\\
    & -\,\frac{27}4\, \frac{x_1x_2x_3}{(x_1+x_2+x_3)^3},  &\qquad&(p,q)=(2,1).
  \end{aligned}
\right.
\end{equation}
Substituting (\ref{eq:68}) yields a degree-$6$ rational map $\tilde
t\mapsto\zeta$, i.e., the Bely\u\i\ map
\begin{equation}
\label{eq:decomp}
  \zeta = \frac{(1+\tilde t^3)^2}{4\tilde t^3}=\left\{
  \begin{aligned}
    &\frac{27}{4}\frac{s^2}{(1-s)^3}\circ \frac{1-\tilde t +  \tilde
      t^2}{(1+\tilde t)^2}&\qquad&(p,q)=(1,2);\\
    &-\frac{27}{4}\frac{s}{(1-s)^3}\circ \frac{(1+\tilde t)^2}{1-\tilde t +  \tilde
      t^2}&\qquad&(p,q)=(2,1),
  \end{aligned}
\right.
\end{equation}
as the covering
$\pi^{(3)}_{1,2}\colon\mathcal{C}^{(3)}_{1,2}\cong\mathcal{C}^{(2)}_{1,2}\to\mathbb{P}^1_\zeta$,
resp.\ $\pi^{(3)}_{2,1}\colon\mathcal{C}^{(3)}_{2,1}\cong\mathcal{C}^{(2)}_{2,1}\to\mathbb{P}^{1}_\zeta$.
Equation~(\ref{eq:decomp}) exhibits the compositions
$\pi^{(2)}_{1,2}=\phi^{(1)}_{1,2}\circ\phi^{(2)}_{1,2}$
and~$\pi^{(2)}_{2,1}=\phi^{(1)}_{2,1}\circ\phi^{(2)}_{2,1}$.  The maps
$\phi^{(2)}_{1,2}$ resp.\ $\phi^{(2)}_{2,1}$ and $\phi^{(1)}_{1,2}$
resp.\ $\phi^{(1)}_{2,1}$, i.e., $\tilde t\mapsto s$ and $s\mapsto\zeta$,
are of degrees $2$ and~$3$.  They come from Thm.~\ref{thm:55}, if one takes
account of the M\"obius transformation relating $t$ and~$\tilde t$.
\end{example}

\begin{example}
\label{ex:13}
  $\{p,q\}=\{1,3\}$.  The curves
$\mathcal{C}^{(3)}_{1,3},\mathcal{C}^{(3)}_{3,1}$ are of genus zero (see
Cor.~\ref{cor:64}).  First consider $(p,q)=(1,3)$.  The defining equation
of~$\mathcal{C}^{(3)}_{1,3}\subset\mathbb{P}^2$ is
  \begin{equation}
    \label{eq:namelessconic}
    x_1^2+x_2^2+x_3^2 + x_1x_2 + x_2x_3 + x_3x_1 = 0,
  \end{equation}
  which follows from Thm.~\ref{thm:long3}, or more directly by
  eliminating~$x_4$ from the equations $\sigma_1=0$, $\sigma_2=0$.  This is
  a conic that does not go through the fundamental points
  of~$\mathbb{P}^2$. It can be parametrized by inspection, with
  parameter~$u\in\mathbb{P}^1$, as
  \begin{equation}
    \label{eq:3param}
    [x_1:x_2:x_3] = [\omega(1-\bar\omega\,u)(1+2\bar\omega\,u):\bar\omega(1-\omega\,u)(1+2\omega\,u):(1-u)(1+2u)],
  \end{equation}
  where $\omega:=\varepsilon_3$.  Substituting (\ref{eq:3param}) into the
  function $\zeta=\zeta(x_1,x_2,x_3)$ of Thm.~\ref{thm:long3} yields a
  degree-$24$ rational map $u\mapsto\zeta$, i.e., the Bely\u\i\ map
\begin{subequations}
  \begin{sizealign}{\small}
    \label{eq:largealtaa}
    \zeta &= -\,\frac{256}{27}\frac{(x_1x_2x_3)^3\, (x_1+x_2+x_3)^3}
    {(x_1+x_2)^4\,(x_2+x_3)^4\,(x_3+x_1)^4}
    \\
    \label{eq:largealta}
    &= -256\,\frac{u^3(1-u^3)^3\,(1+8u^3)^3}{(1-20u^3-8u^6)^4}
    = 1 - \frac{(1+8u^6)^2\,(1+88u^3-8u^6)^2}{(1-20u^3-8u^6)^4}
    \\
    &= -\,\frac{256}{27}\frac{s^3}{(1-s)^4}
    \circ
    \frac{(1-t)(1+3t^2)}{(1+t)^3}
    \circ \left[\left(\frac{\omega-\bar\omega}{3}\right)\,
    \left(\frac{1-2u-2u^2}{1+2u^2}\right)\right],
    \label{eq:largealtb}
  \end{sizealign}
\end{subequations}
  as the covering
  $\pi_{1,3}^{(3)}\colon\mathcal{C}^{(3)}_{1,3}\cong\mathbb{P}^1_u\to\mathbb{P}_\zeta^1$.
  Equation~(\ref{eq:largealtb}) exhibits the composition $\pi_{1,3}^{(3)} =
  \phi_{1,3}^{(1)} \circ \phi_{1,3}^{(2)} \circ \phi_{1,3}^{(3)}$.  The
  maps $\phi_{1,3}^{(3)},\phi_{1,3}^{(2)},\phi_{1,3}^{(1)}$, i.e.,
  $u\mapsto t$, $t\mapsto s$, $s\mapsto\zeta$, are of degrees $2,3,4$.
  They come from Theorem~\ref{thm:55} and the fact that the parameter~$t$
  on any curve $\mathcal{C}^{(2)}_{p,q}$ equals $(x_1+\nobreak
  x_2)/\allowbreak(x_1-\nobreak x_2)$.

  The top Schwarz curve~$\mathcal{C}^{(4)}_{1,3}\subset\mathbb{P}^3$ is
  birationally equivalent to~$\mathcal{C}^{(3)}_{1,3}$, since the map
  $\phi^{(4)}\colon \mathcal{C}^{(4)}_{1,3}\to\mathcal{C}^{(3)}_{1,3}$ is
  of degree~$1$; so it too can be parametrized by~$u$.  Solving the
  equation $\sigma_1=0$ for~$x_4=x_4(u)$ yields $x_4=-3u$, hence
  \begin{multline}
    \label{eq:4param}
    [x_1:x_2:x_3:x_4]\\ = [\omega(1-\bar\omega\,u)(1+2\bar\omega\,u):\bar\omega(1-\omega\,u)(1+2\omega\,u):(1-u)(1+2u):-3u]
  \end{multline}
  is an (asymmetric) rational parametrization of~$\mathcal{C}^{(4)}_{1,3}$.
  Equation (\ref{eq:largealta})~can be viewed as defining the
  degree\nobreakdash-$24$ covering
  $\pi_{1,3}^{(4)}\colon\mathcal{C}^{(4)}_{1,3}\cong\mathbb{P}^1_u\to\mathbb{P}_\zeta^1$.

  The treatment of the case $(p,q)=(3,1)$ is similar.  The genus-zero
  curve $\mathcal{C}^{(3)}_{3,1}\subset\mathbb{P}^2$ is obtained from
  $\mathcal{C}^{(3)}_{1,3}\subset\mathbb{P}^2$ by Cremona inversion
  ($x_i=1/x'_i$), i.e., by the standard quadratic transformation
  ($x_i=x_j'x_k'$).  It has defining equation
  \begin{equation}
    x_1^2x_2^2 + x_2^2x_3^2 + x_3^2x_1^2 + x_1^2 x_2x_3 + x_2^2 x_3x_1 + x_3^2
    x_1x_2 = 0
  \end{equation}
  and is a trinodal quartic; it goes through the fundamental points
  $[1:\nobreak 0:\nobreak 0],\allowbreak[0:\nobreak 1:\nobreak
    0],\allowbreak[0:\nobreak 0:\nobreak 1]$, and each is an ordinary
  double point (`node').  The smooth points
  $[1:\nobreak\omega:\nobreak\bar\omega]$,
  $[1:\nobreak\bar\omega:\nobreak\omega]$ are notable for being appended
  limit points in the sense of~\S\,\ref{subsec:61}: they come indirectly
  rather than directly from trinomial roots.

  By undoing the quadratic transformation, a rational parametrization
  by~$u$ of~$\mathcal{C}^{(3)}_{3,1}$ can be obtained
  from~(\ref{eq:3param}), and one of~$\mathcal{C}^{(4)}_{3,1}$
  from~(\ref{eq:4param}).  Thus the degree\nobreakdash-$24$ rational map
  $u\mapsto\zeta$ given in~(\ref{eq:largealta}) can be used as
  $\pi_{3,1}^{(3)}\colon\mathcal{C}^{(3)}_{3,1}\cong\mathbb{P}^1_u\to\mathbb{P}_\zeta^1$
  and~$\pi_{3,1}^{(4)}\colon\mathcal{C}^{(4)}_{3,1}\cong\mathbb{P}^1_u\to\mathbb{P}_\zeta^1$.
\end{example}

\begin{example}
\label{ex:14}
  $\{p,q\}=\{1,4\}$.  A discussion of the case $(p,q)=(1,4)$ will suffice.
The defining equation of~$\mathcal{C}^{(3)}_{1,4}\subset\mathbb{P}^2$ is
\begin{equation}
\label{eq:trianglecubic}
x_1^3 + x_2^3 + x_3^3 + x_1x_2^2 + x_1x_3^2 + x_2x_3^2 + x_2x_1^2 +
x_3x_1^2 + x_3x_2^2 + x_1x_2x_3 = 0,
\end{equation}
which follows from Thm.~\ref{thm:long3}, or more directly by eliminating
$x_4,x_5$ from the equations $\sigma_1=0$, $\sigma_2=0$, $\sigma_3=0$.
This is a smooth cubic that does not go through the fundamental points
of~$\mathbb{P}^2$.  It is elliptic, i.e., of genus~$1$, with Klein--Weber
invariant $j=-5^2/2$.  It can therefore be uniformized with the aid of
elliptic functions, but it is easier to construct a multivalued
parametrization with radicals.  As~usual, let $t=\allowbreak (x_1+\nobreak
x_2)/\allowbreak (x_1-\nobreak x_2)$, so that $[x_1:\nobreak
  x_2]=\allowbreak[t+\nobreak 1:\nobreak t-\nobreak 1]$, and notice that as
Thm.~\ref{thm:long3} predicts, Eq.~(\ref{eq:trianglecubic})~is of
degree~$n-2=3$ in~$x_3$.  By symmetry, $x_3,x_4,x_5$ are the three roots,
and they are computable in~terms of radicals from~$t$ by Cardano's formula.
It follows that each of
$\mathcal{C}^{(3)}_{1,4},\mathcal{C}^{(4)}_{1,4},\mathcal{C}^{(5)}_{1,4}$
has a multivalued parametrization with radicals in~terms of~$t$.  These
parametrizations are respectively $3$,~$6$, and $6$-valued.
\end{example}

The technique of the last example immediately yields the following theorem.

\begin{theorem}
  For all coprime\/ $p\ge1$, $q\ge2$ with\/ $n\rsmcoloneq p+q\le6$, one can
construct multivalued parametrizations with radicals for the subsidiary
curve\/ $\mathcal{C}^{(q)}_{p,q}\subset\mathbb{P}^{q-1}$ and the top
curve\/ $\mathcal{C}^{(n)}_{p,q}\subset\mathbb{P}^{n-1}$, respectively\/
$(p+1)_{q-2}$-valued and\/ $(n-2)!$-valued.
\end{theorem}

\medskip
\section{Identities with free parameters}
\label{sec:freeparams}

The results of \S\S\,\ref{sec:curves} and~\ref{sec:genera} will now be put
to use, by deriving an interesting collection of parametrized
hypergeometric identities.  The source of many is Thm.~\ref{thm:invrep},
which expressed certain ${}_nF_{n-1}$'s in~terms of algebraic functions.
The expressions in parts (i) and~(ii) of that theorem involve the roots
$x_1,\dots,x_q$, resp.\ $x_1,\dots,x_n$, of the trinomial equation
\begin{equation}
  x^n-g\,x^p-\beta=0,
\end{equation}
with $n:=p+q$ and $\gcd(p,q)=1$.  It follows that any parametrization of
the Schwarz curve $\mathcal{C}^{(q)}_{p,q}$,
resp.~$\mathcal{C}^{(n)}_{p,q}$, will yield a hypergeometric identity.  The
curves $\mathcal{C}^{(q)}_{p,q},\mathcal{C}^{(n)}_{p,q}$ of genus zero were
classified in Cor.~\ref{cor:64}, and the resulting identities are given in
\S\S\,\ref{subsec:71}, \ref{subsec:72},~\ref{subsec:73} below, the
respective parameters used being $s,t,u$.  These are respectively the
rational parameters for any~$\mathcal{C}^{(1)}_{p,q}$,
any~$\mathcal{C}^{(2)}_{p,q}$, and (by Example~\ref{ex:13}) the plane
curves $\mathcal{C}^{(3)}_{1,3},\mathcal{C}^{(3)}_{3,1}$.

Each identity derived from Thm.~\ref{thm:invrep} in this section is really
a family of identities: the represented ${}_nF_{n-1}$ depends on a discrete
parameter $\ell\in\mathbb{Z}$, and the identity involves a rational
function $F_\ell=\allowbreak F_\ell(A,B;y)$.  The functions~$F_\ell$ were
defined in~\S\,\ref{sec:binomial} (see Table~\ref{tab:1} and
Thm.~\ref{thm:iteratable}).  The reader will recall that in particular,
$F_0\equiv1$ and $F_1(A,B;y)=\allowbreak y/[(1-\nobreak B)y+\nobreak B]$.

Each of the ${}_nF_{n-1}$'s also depends on a parameter $a\in\mathbb{C}$.
If $a$~is chosen so that no upper parameter of the corresponding
differential equation~$E_n$ differs by an integer from a lower one, the
monodromy group of the~$E_n$ will be of the imprimitive irreducible type
characterized in Thm.~\ref{thm:BH0}; and if $a\in\mathbb{Q}$, the group
will be finite.  (The case of \emph{equal} upper and lower parameters,
permitting `cancellation,' is possible only for a finite number of choices
of~$a$, such as~$a=\pm1$ when $\ell=0$; it was mentioned in
Thm.~\ref{thm:BH1}.)  One must treat with care the possibility that one of
the lower parameters may be a non-positive integer, in which case
${}_nF_{n-1}$ is undefined (though it may still be possible to interpret
the identity in a limiting sense; cf.\ Lemma~\ref{lem:1}).  \emph{It~is
  assumed for simplicity in this section that\/ $a\in\mathbb{C}$ is chosen
  so that this does not occur, and so that no~division by zero occurs.}

Several of the identities below are rationally parametrized formulas for
${}_{n+1}F_n$'s rather than ${}_{n}F_{n-1}$'s.  They come from
Thm.~\ref{thm:invrepg0} rather than Thm.~\ref{thm:invrep}, and instead of
the rational functions $F_\ell$, $\ell\in\mathbb{Z}$, they involve
interpolating functions $G_0,G_1$ that were defined in
Thm.~\ref{thm:invrepg0} in~terms of~${}_2F_1,{}_3F_2$.  Each of these
identities has an additional free parameter $c\in\mathbb{C}$, and a similar
caveat applies.

\subsection{Parametrizations by {$s$}}
\label{subsec:71}

The rational parametrization of the curve
$\mathcal{C}^{(q)}_{p,q}=\mathcal{C}^{(1)}_{p,1}$ by $s\in\mathbb{P}^1$
given in Eqs.\ (\ref{eq:focusedalignat}) and~(\ref{eq:basecoin}), when
substituted into part~(i) of each of Thms.\ \ref{thm:invrep}
and~\ref{thm:invrepg0}, yields the following.

\begin{theorem}
  For each\/ $p\ge1$, with\/ $n\rsmcoloneq p+1$, define a degree-$n$ Bely\u\i\
  map\/ $\mathbb{P}_s^1\to\mathbb{P}_\zeta^1$ by
  \begin{equation*}
    \zeta=\zeta_{p,1}(s)\rsmcoloneq -\,\frac{n^n}{p^p}\,\frac{s}{(1-s)^n}.
  \end{equation*}
  Then for all\/ $a\in\mathbb{C}$
  and\/ $\ell\in\mathbb{Z}$, one has that near\/ $s=0$,
\begin{equation*}
{}_nF_{n-1}\left(
\begin{array}{l}
\frac{a}n,\dots,\frac{a+(n-1)}n\\
\frac{a-\ell+1}p,\dots,\frac{a-\ell+p}p
\end{array}
\!\!\biggm|\zeta_{p,1}(s)\right)
=(1-s)^a\, F_\ell\!\left(-a,-p;\,(1-s)^{-1}\right).
\end{equation*}
  Moreover, for all\/ $a\in\mathbb{C}$, $c\in\mathbb{C}$ and\/ $\ell=0,1$,
  one has that near\/ $s=0$,
  \begin{multline*}
{}_{n+1}F_{n}\left(
\begin{array}{ll}
\frac{a}n,\dots,\frac{a+(n-1)}n;&\frac{a+c-\ell}p\\
\frac{a-\ell}p,\dots,\frac{a-\ell+(p-1)}p;&\frac{a+c-\ell+p}p
\end{array}
\!\!\biggm|\zeta_{p,1}(s)\right)\\
= (1-s)^a\, G_\ell\!\left(-a,-p,-c;\,(1-s)^{-1}\right).
  \end{multline*}
\label{thm:71}
\end{theorem}

\begin{proof}
  Straightforward, as $g^{-1}x_1 = (1-s)^{-1}$
  by~(\ref{eq:focusedalignat}).
\end{proof}

\begin{remarkaftertheorem}
  For each of $\ell=0,1$, the second identity reduces to the first when
  $c=0$; and as~$c\to\infty$, it reduces to a version of the first in which
  $\ell$~is incremented by~$1$.  Interpolation of this sort is familiar
  from Thm.~\ref{thm:invrepg0} and will be seen repeatedly.
\end{remarkaftertheorem}

Similarly, the parametrization of the top curve
$\mathcal{C}^{(n)}_{p,q}=\mathcal{C}^{(2)}_{1,1}\cong\mathcal{C}^{(1)}_{1,1}$
by $s\in\mathbb{P}^1$, substituted into part~(ii) of each of Thms.\
\ref{thm:invrep} and~\ref{thm:invrepg0}, yields the following.

\begin{theorem}
 As in Theorem~{\rm\ref{thm:71},} let
 \begin{equation*}
   \zeta_{1,1}(s)\rsmcoloneq-\,\frac{4s}{(1-s)^2}.
 \end{equation*}
 Then for all\/ $a\in\mathbb{C}$ and\/ $\ell\in\mathbb{Z}$, and\/
  $\kappa=0,1$, one has that near\/ $s=1$,
\begin{multline*}
{}_2F_{1}\left(
\begin{array}{l}
-a+\ell+\frac\kappa2,\,a+\frac\kappa2\\
\frac12+\kappa
\end{array}
\!\!\biggm|\frac1{\zeta_{1,1}(s)}\right)\\
=\frac{(-)^\kappa(a+\kappa/2)_{1-\ell-\kappa}}{(a)_{1-\ell}}\,\left[-\zeta_{1,1}(s)/4\right]^{\kappa/2}\\
\qquad{}\times\tfrac12\left[s^a\,F_\ell\!\left(-a,\tfrac12;\,s^{-1}\right) 
          +(-)^\kappa s^{-a}\,F_\ell\!\left(-a,\tfrac12;\,s\right)\right].
\end{multline*}
  Moreover, for all\/ $a\in\mathbb{C}$, $c\in\mathbb{C}$ and\/ $\ell=0,1$,
  and\/ $\kappa=0,1$, one has that near\/ $s=1$,
\begin{multline*}
{}_{3}F_{2}\left(
\begin{array}{ll}
-a+\ell+1+\frac\kappa2,\,a+\frac\kappa2;&{-a-c+\ell+\frac\kappa2}\\
\tfrac12+\kappa;&{-a-c+\ell+1+\frac\kappa2}
\end{array}
\!\!\biggm|\frac1{\zeta_{1,1}(s)}\right)\\
=\frac{(-)^\kappa(a+\kappa/2)_{-\ell-\kappa}(a+c-\ell-\kappa/2)}{(a)_{-\ell}(a+c-\ell)}\,\left[-\zeta_{1,1}(s)/4\right]^{\kappa/2}\\
\qquad{}\times\tfrac12\left[s^a\,G_\ell\!\left(-a,\tfrac12,-c;\,s^{-1}\right)
+(-)^\kappa s^{-a}\,G_\ell\!\left(-a,\tfrac12,-c;\,s\right)\right].
\end{multline*}
\label{thm:72}
\end{theorem}

\begin{proof}
  Straightforward, as $\beta^{-1/2}x_1=s^{-1/2}$ and
$\beta^{-1/2}x_2=-s^{1/2}$.
\end{proof}

The cases $\ell=0,1$ of the identities of Thm.~\ref{thm:71} were previously
obtained by Gessel and Stanton using series manipulations~\cite{Gessel82}.
(See their Eqs.\ (5.10)--(5.13),\allowbreak(5.15).)  When $p=2$, the
$\ell=0,1$ cases of the first identity become one-parameter specializations
of Bailey's first cubic transformation of~${}_3F_2$ and its companion.

By comparison, the two identities of Thm.~\ref{thm:72} are elementary.  It
should be possible to derive them, or at~least the first, by using
contiguous function relations or other classical hypergeometric techniques.

\subsection{Parametrizations by {$t$}}
\label{subsec:72}

The rational parametrization of the curve
$\mathcal{C}^{(q)}_{p,q}=\mathcal{C}^{(2)}_{p,2}$ by $t\in\mathbb{P}^1$
given in Eq.~(\ref{eq:basecoin0}) and Lemma~\ref{lem:54}, and the fact that
$[x_1:x_2]=[t+1:t-1]$, when substituted into part~(i) of each of Thms.\
\ref{thm:invrep} and~\ref{thm:invrepg0}, yield the following.

\begin{theorem}
  For each odd\/ $p\ge1$, with\/ $n\rsmcoloneq p+2$, define a map\/
  $\mathbb{P}_t^1\to\mathbb{P}_\zeta^1$ by
\begin{equation*}
  \zeta=\zeta_{p,2}(t)\rsmcoloneq
\frac{4n^n}{p^p}\, \frac{t^2(1-t^2)^{2p}\left[(1+t)^p+(1-t)^p\right]^p}{\left[(1+t)^n+(1-t)^n\right]^n}.
\end{equation*}
Then for all\/ $a\in\mathbb{C}$ and\/ $\ell\in\mathbb{Z}$, and\/
  $\kappa=0,1$, one has that near\/ $t=0$,
\begin{sizealign}{\small}
&{}_nF_{n-1}\left(
\begin{array}{l}
\frac{a}n+\frac\kappa2,\dots,\frac{a+(n-1)}n+\frac\kappa2\\
\frac{a-\ell+1}p+\frac\kappa2,\dots,\frac{a-\ell+p}p+\frac\kappa2;\,\frac12+\kappa
\end{array}
\!\!\biggm|\zeta_{p,2}(t)\right)\notag\\
&= \frac{(-)^\kappa(a+n\kappa/2)_{1-\ell-\kappa}}{(a)_{1-\ell}}
\left[\frac{4p^p}{n^n}\,\zeta_{p,2}(t)\right]^{-\kappa/2}\notag\\
&\qquad{}\times\tfrac12\biggl[ (1+t)^{-2a}\,F_\ell\!\left(-a,-p/2;\,(1+t)^2\left[\frac{(1+t)^p+(1-t)^p}{(1+t)^n+(1-t)^n}\right]\right)\notag\\
&\hphantom{\quad{}\times\tfrac12\bigl[}\qquad{}+(-)^\kappa (1-t)^{-2a}\,F_\ell\!\left(-a,-p/2;\,(1-t)^{2}\left[\frac{(1+t)^p+(1-t)^p}{(1+t)^n+(1-t)^n}\right]\right)\biggr]\notag\\
&\qquad{}\times\left[\frac{(1+t)^n+(1-t)^n}{(1+t)^p+(1-t)^p}\right]^a.\notag
\end{sizealign}
  Moreover, for all\/ $a\in\mathbb{C}$, $c\in\mathbb{C}$ and\/ $\ell=0,1$,
  and\/ $\kappa=0,1$, one has that near\/ $t=0$,
\begin{sizealign}{\small}
&{}_{n+1}F_{n}\left(
\begin{array}{ll}
\frac{a}n+\frac\kappa2,\dots,\frac{a+(n-1)}n+\frac\kappa2;&\frac{a+c-\ell}p+\frac\kappa2\\
\frac{a-\ell}p+\frac\kappa2,\dots,\frac{a-\ell+(p-1)}p+\frac\kappa2;\,\frac12+\kappa;&\frac{a+c-\ell+p}p+\frac\kappa2
\end{array}
\!\!\biggm|\zeta_{p,2}(t)\right)\notag\\
&= \frac{(-)^\kappa(a+n\kappa/2)_{-\ell-\kappa}(a+c-\ell+p\kappa/2)}{(a)_{-\ell}(a+c-\ell)}
\left[\frac{4p^p}{n^n}\,\zeta_{p,2}(t)\right]^{-\kappa/2}\notag\\
&\qquad{}\times\tfrac12\biggl[ (1+t)^{-2a}\,G_\ell\!\left(-a,-p/2,-c;\,(1+t)^2\left[\frac{(1+t)^p+(1-t)^p}{(1+t)^n+(1-t)^n}\right]\right)\notag\\
&\hphantom{\qquad{}\times\tfrac12\bigl[}\quad{}+(-)^\kappa (1-t)^{-2a}\,G_\ell\!\left(-a,-p/2,-c;\,(1-t)^{2}\left[\frac{(1+t)^p+(1-t)^p}{(1+t)^n+(1-t)^n}\right]\right)
\biggr]\notag\\
&\qquad{}\times\left[\frac{(1+t)^n+(1-t)^n}{(1+t)^p+(1-t)^p}\right]^a.\notag
\end{sizealign}
\label{thm:73}
\end{theorem}

\begin{remarkaftertheorem}
  The even function $\zeta=\zeta_{p,2}(t)$ is a degree-$[n(n-1)]$
  Bely\u\i\ map.  The case $\ell=0$, $\kappa=0$ of the first identity
  appeared in~\S\,\ref{sec:intro} as the sample result~(\ref{eq:sample}).
  Note that if $\kappa=0$, there is simplification: the right-hand
  prefactor becomes unity.
\end{remarkaftertheorem}

Similarly, the rational parametrization of the top curve
$\mathcal{C}^{(n)}_{p,q}=\mathcal{C}^{(3)}_{1,2}\cong\mathcal{C}^{(2)}_{1,2}$
given in~(\ref{eq:68}) by $\tilde t\in\mathbb{P}^1$ (related to~$t$ by a
M\"obius transformation), substituted into part~(ii) of each of Thms.\
\ref{thm:invrep} and~\ref{thm:invrepg0}, yields the following.

\begin{theorem}
  Define a degree-\/$6$ Bely\u\i\ map\/ $\mathbb{P}^1_{\tilde
  t}\to\mathbb{P}^1_\zeta$ by
  \begin{equation*}
  \zeta=\zeta_{1,2}(\tilde t)\rsmcoloneq \frac{(1+\tilde t^3)^2}{4\tilde t^3}
  = 1 + \frac{(1-\tilde t^3)^2}{4\tilde t^3}.
\end{equation*}
Then for all\/ $a\in\mathbb{C}$ and\/ $\ell\in\mathbb{Z}$, and\/
  $\kappa=0,1,2$, one has that near\/ $\tilde t=0$,
\begin{sizealign}{\footnotesize}
&{}_3F_{2}\left(
\begin{array}{l}
-a+\ell+\frac\kappa3;\,\frac{a}{2}+\frac\kappa3,\,\frac{a+1}2+\frac\kappa3\\
b_1(\kappa),\,b_2(\kappa)
\end{array}
\!\!\biggm|\frac1{\zeta_{1,2}(\tilde t)}\right)\notag\\
&=
\frac{(-)^\kappa(1)_\kappa(a+2\kappa/3)_{1-\ell-\kappa}}{(a)_{1-\ell}}\,\left[\frac{4}{27}\,\zeta_{1,2}(\tilde{t})\right]^{\kappa/3}
\notag\\
&\quad{}\times
\frac13\Biggl\{
\left[\frac{(1+\tilde t)^3}{1+\tilde t^3}\right]^{-a}F_\ell\left(-a,\tfrac13;\,\frac{(1+\tilde  t)^3}{1+\tilde t^3}\right)
+\bar\omega^\kappa
\left[\frac{(1+\omega\tilde t)^3}{1+\tilde t^3}\right]^{-a}F_\ell\left(-a,\tfrac13;\,\frac{(1+\omega\tilde t)^3}{1+\tilde t^3}\right)\notag\\
&\hphantom{
\quad{}\times
\frac13\Biggl\{
\left[\frac{(1+\tilde t)^3}{1+\tilde t^3}\right]^{-a}F_\ell\left(-a,\tfrac13;\,\frac{(1+\tilde  t)^3}{1+\tilde t^3}\right)
}
{}+\omega^{\kappa}
\left[\frac{(1+\bar\omega\tilde t)^3}{1+\tilde t^3}\right]^{-a}F_\ell\left(-a,\tfrac13;\,\frac{(1+\bar\omega\tilde t)^3}{1+\tilde t^3}\right)
\Biggr\}.\notag
\end{sizealign}
  Moreover, for all\/ $a\in\mathbb{C}$, $c\in\mathbb{C}$ and\/ $\ell=0,1,2$,
  and\/ $\kappa=0,1$, one has that\hfil\break near\/ $\tilde t=0$,
\begin{sizealign}{\footnotesize}
&{}_4F_{3}\left(
\begin{array}{ll}
-a+\ell+1+\frac\kappa3;\,\frac{a}{2}+\frac\kappa3,\,\frac{a+1}2+\frac\kappa3;&-a-c+\ell+\frac\kappa3\\
b_1(\kappa),\,b_2(\kappa);&-a-c+\ell+1+\frac\kappa3
\end{array}
\!\!\biggm|\frac1{\zeta_{1,2}(\tilde t)}\right)\notag\\
&=
\frac{(-)^\kappa(1)_\kappa(a+2\kappa/3)_{-\ell-\kappa}(a+c-\ell-\kappa/3)}{(a)_{-\ell}(a+c-\ell)}\,\left[\frac{4}{27}\,\zeta_{1,2}(\tilde t)\right]^{\kappa/3}
\notag\\
&\quad{}\times
\frac13\Biggl\{
\left[\frac{(1+\tilde t)^3}{1+\tilde t^3}\right]^{-a}G_\ell\left(-a,\tfrac13,-c;\,\frac{(1+\tilde t)^3}{1+\tilde t^3}\right)
+\bar\omega^\kappa
\left[\frac{(1+\omega\tilde t)^3}{1+\tilde t^3}\right]^{-a}G_\ell\left(-a,\tfrac13,-c;\,\frac{(1+\omega\tilde t)^3}{1+\tilde t^3}\right)\notag\\
&\hphantom{
\quad{}\times
\frac13\Biggl\{
\left[\frac{(1+\tilde t)^3}{1+\tilde t^3}\right]^{-a}G_\ell\left(-a,\tfrac13,-c;\,\frac{(1+\tilde t)^3}{1+\tilde t^3}\right)
}
{}+\omega^{\kappa}
\left[\frac{(1+\bar\omega\tilde t)^3}{1+\tilde t^3}\right]^{-a}G_\ell\left(-a,\tfrac13,-c;\,\frac{(1+\bar\omega\tilde t)^3}{1+\tilde t^3}\right)
\Biggr\}.\notag
\end{sizealign}
In these two identities, the lower parameters\/ $b_1,b_2$ are\/
$\tfrac13,\tfrac23${\rm;} or\/ $\tfrac23,\tfrac43${\rm;} or\/
$\tfrac43,\tfrac53${\rm;} depending on whether\/ $\kappa=0$, $1$, or\/ $2$.
Also, $\omega\rsmcoloneq\varepsilon_3 = \exp(2\pi{\rm i}/3)$ and
$\bar\omega\rsmcoloneq{\varepsilon_3}^2$.
\end{theorem}

\begin{proof}
  Straightforward, as $\beta=\sigma_3=x_1x_2x_3=1+\tilde t^3$
  by~(\ref{eq:68}).
\end{proof}


One can also derive identities from the parametrization by $\tilde
t\in\mathbb{P}^1$ of the birationally equivalent top curve
$\mathcal{C}^{(3)}_{2,1}\cong\mathcal{C}^{(2)}_{2,1}$, given along with the
parametrization of $\mathcal{C}^{(3)}_{1,2}\cong\mathcal{C}^{(2)}_{1,2}$
in~(\ref{eq:68}).  Details are left to the reader.

\subsection{Parametrizations by {$u$}}
\label{subsec:73}

The rational parametrization of the plane curve
$\mathcal{C}^{(q)}_{p,q}=\mathcal{C}^{(3)}_{1,3}$ by $u\in\mathbb{P}^1$
given in~(\ref{eq:3param}), in Example~\ref{ex:13}, when substituted into
part~(i) of each of Thms.\ \ref{thm:invrep} and~\ref{thm:invrepg0}, yields
the following.

\begin{theorem}
  Define a degree-\/$24$ Bely\u\i\ map\/
  $\mathbb{P}_u^1\to\mathbb{P}^1_\zeta$ by
\begin{sizedisplaymath}{\small}
  \zeta=\zeta_{1,3}(u)\rsmcoloneq
-256\,\frac{u^3(1-u^3)^3(1+8u^3)^3}{(1-20u^3-8u^6)^4}
= 1 - \frac{(1+8u^6)^2(1+88u^3-8u^6)^2}{(1-20u^3-8u^6)^4}.
\end{sizedisplaymath}
Then for all\/ $a\in\mathbb{C}$ and\/ $\ell\in\mathbb{Z}$, and\/
  $\kappa=0,1,2$, one has that near\/ $u=0$,
\begin{sizealign}{\small}
&{}_4F_{3}\left(
\begin{array}{l}
\frac{a}4+\frac\kappa3,\dots,\frac{a+3}4+\frac\kappa3\\
{a-\ell+1}+\frac\kappa3;\,b_1(\kappa),\,b_2(\kappa)
\end{array}
\!\!\biggm|\zeta_{1,3}(u)\right)\notag\\
&= \frac{(-)^\kappa(1)_\kappa(a+4\kappa/3)_{1-\ell-\kappa}}{(a)_{1-\ell}}
\left[-\,\frac{27}{256}\,\zeta_{1,3}(u)\right]^{-\kappa/3}\notag\\
&\quad{}\times
\frac13\Biggl\{
\left[\frac{(1-u)^3(1+2u)^3}{1-20u^3-8u^6}\right]^{-a}F_\ell\left(-a,-\tfrac13;\,\frac{(1-u)^3(1+2u)^3}{1-20u^3-8u^6}\right)
\notag\\
&\quad\qquad\qquad{}+\bar\omega^{\kappa}
\left[\frac{(1-\omega u)^3(1+2\omega u)^3}{1-20u^3-8u^6}\right]^{-a}F_\ell\left(-a,-\tfrac13;\,\frac{(1-\omega u)^3(1+2\omega u)^3}{1-20u^3-8u^6}\right)
\notag\\
&\quad\qquad\qquad{}+\omega^{\kappa}
\left[\frac{(1-\bar\omega u)^3(1+2\bar\omega u)^3}{1-20u^3-8u^6}\right]^{-a}F_\ell\left(-a,-\tfrac13;\,\frac{(1-\bar\omega u)^3(1+2\bar\omega u)^3}{1-20u^3-8u^6}\right)
\Biggr\}.\notag
\end{sizealign}
Moreover, for all\/ $a\in\mathbb{C}$, $c\in\mathbb{C}$ and\/ $\ell=0,1$,
and\/ $\kappa=0,1,2$, one has that\hfil\break near\/ $u=0$,
\begin{sizealign}{\small}
&{}_5F_{4}\left(
\begin{array}{ll}
\frac{a}4+\frac\kappa3,\dots,\frac{a+3}4+\frac\kappa3;&a+c-\ell+\frac\kappa3\\
{a-\ell}+\frac\kappa3;\,b_1(\kappa),\,b_2(\kappa);&a+c-\ell+1+\frac\kappa3
\end{array}
\!\!\biggm|\zeta_{1,3}(u)\right)\notag\\
&= \frac{(-)^\kappa(1)_\kappa(a+4\kappa/3)_{-\ell-\kappa}(a+c-\ell+\kappa/3)}{(a)_{-\ell}(a+c-\ell)}
\left[-\,\frac{27}{256}\,\zeta_{1,3}(u)\right]^{-\kappa/3}\notag\\
&\quad{}\times
\frac13\Biggl\{
\left[\frac{(1-u)^3(1+2u)^3}{1-20u^3-8u^6}\right]^{-a}G_\ell\left(-a,-\tfrac13,-c;\,\frac{(1-u)^3(1+2u)^3}{1-20u^3-8u^6}\right)
\notag\\
&\quad\qquad\qquad{}+\bar\omega^{\kappa}
\left[\frac{(1-\omega u)^3(1+2\omega u)^3}{1-20u^3-8u^6}\right]^{-a}G_\ell\left(-a,-\tfrac13,-c;\,\frac{(1-\omega u)^3(1+2\omega u)^3}{1-20u^3-8u^6}\right)
\notag\\
&\quad\qquad\qquad{}+\omega^{\kappa}
\left[\frac{(1-\bar\omega u)^3(1+2\bar\omega u)^3}{1-20u^3-8u^6}\right]^{-a}G_\ell\left(-a,-\tfrac13,-c;\,\frac{(1-\bar\omega u)^3(1+2\bar\omega u)^3}{1-20u^3-8u^6}\right)
\Biggr\}.\notag
\end{sizealign}
In both identities, the lower parameters\/ $b_1,b_2$ are defined as in the
previous theorem.
\end{theorem}

\begin{proof}
  Straightforward, as $g=\sigma_3=x_1x_2x_3=1-20u^3-8u^6$ by~(\ref{eq:3param}).
\end{proof}

One can similarly derive identities from the rational parametrizations of
the pair of top curves
$\mathcal{C}^{(4)}_{1,3}\cong\mathcal{C}^{(3)}_{1,3}$ and
$\mathcal{C}^{(4)}_{3,1}\cong\mathcal{C}^{(3)}_{3,1}$, also given in
Example~\ref{ex:13}, by substituting them into Thms.\ \ref{thm:invrep}(ii)
and~\ref{thm:invrepg0}(ii).  The resulting identities involve
${}_4F_3,{}_5F_4$, like the preceding.

\section{Identities without free parameters}
\label{sec:nofreeparams}

In this final section an alternative approach to the parametrizing of
certain ${}_nF_{n-1}$'s, based on computation in rings of symmetric
polynomials, is sketched.  The need for a strengthened approach is
indicated by the nature of the preceding results.  Each identity
in~\S\,\ref{sec:freeparams} was derived from either Thm.~\ref{thm:invrep}
or Thm.~\ref{thm:invrepg0}: it had a free parameter~$a\in\mathbb{C}$ and
was based on a (rational) parametrization of a Schwarz curve, either
$\mathcal{C}_{p,q}^{(q)}$ or~$\mathcal{C}_{p,q}^{(n)}$.  But not many such
curves are of zero or low genus, making it difficult to derive large
numbers of hypergeometric identities, even if multivalued parametrizations
with radicals are allowed.

In~\S\,\ref{subsec:81}, it is shown that if $a\in\mathbb{Z}$, the relevant
curve is a \emph{quotiented} Schwarz curve $\mathcal{C}_{p,q}^{(q;{\rm
    symm.})}$ or~$\mathcal{C}_{p,q}^{(n;{\rm symm.})}$; and if
$qa\in\mathbb{Z}$, resp.\ $na\in\mathbb{Z}$, it is another such curve
$\mathcal{C}_{p,q}^{(q;{\rm cycl.})}$, resp.\ $\mathcal{C}_{p,q}^{(n;{\rm
    cycl.})}$\!.  It is actions of the symmetric group $\mathfrak{S}_q$
resp.~$\mathfrak{S}_n$ and the cyclic group $\mathfrak{C}_q$
resp.~$\mathfrak{C}_n$ that are quotiented~out, and quotienting may lower
the genus; which facilitates explicit parametrization, rational or
otherwise.  The several interesting examples worked~out in
\S\,\ref{subsec:82} and~\S\,\ref{subsec:83} illustrate this, though they
provide representations of algebraic ${}_nF_{n-1}$'s arising from $E_n$'s
that are not of the imprimitive irreducible type characterized in
Thm.~\ref{thm:BH0}, which was previously the focus.  (The difference in
type is due to reducible monodromy and/or a lower hypergeometric parameter
being a non-positive integer.)

\subsection{Theory}
\label{subsec:81}

Theorem~\ref{thm:inta} below is a restatement of the `Birkeland' $\ell=0$
case of Thm.~\ref{thm:invrep}, which is the simplest because when $\ell=0$,
the right-hand function~$F_\ell$ degenerates to unity.  The restatement
emphasizes the role played by symmetric functions of the trinomial roots
$x_1,\dots,x_n$.  As usual, $\sigma_l=\sigma_l(x_1,\dots,x_n)$ denotes the
$l$'th elementary symmetric polynomial, and
$\mathcal{C}_{p,q}^{(n)}\subset\mathbb{P}^{n-1}$ comprises all points
$[x_1:\ldots:x_n]\in\mathbb{P}^{n-1}$ such that
$\sigma_1,\dots,\sigma_{q-1}$ and $\sigma_{q+1},\dots,\sigma_{n-1}$ equal
zero.

\begin{definition}
  For each $p,q\ge1$ with $\gcd(p,q)=1$ and $n:=p+q$, let
  \begin{sizedisplaymath}{\small}
    {\mathcal{F}}_{p,q}^{(q;\kappa)}(a;\,\zeta)\rsmcoloneq
    {}_nF_{n-1}\left(
    \begin{array}{l}
      \frac{a}n+\frac{\kappa}q,\dots,\frac{a+(n-1)}n+\frac{\kappa}q\\
      \frac{a+1}p+\frac{\kappa}q,\dots,
      \frac{a+p}p+\frac{\kappa}q;\,
      \frac1q+\frac{\kappa}q,\dots,\widehat{\frac{q-\kappa}q + \frac{\kappa}q},\dots,\frac{q}q+\frac{\kappa}q
    \end{array}
    \!\!\biggm|\zeta\right)
  \end{sizedisplaymath}
  for $\kappa=0,\dots,q-1$, and 
  \begin{sizedisplaymath}{\small}
    {\mathcal{F}}_{p,q}^{(n;\kappa)}(a;\,\zeta)\rsmcoloneq
    {}_nF_{n-1}\left(
    \begin{array}{l}
      \frac{-a}p+\frac{\kappa}n,\dots,
      \frac{-a+(p-1)}p+\frac{\kappa}n;\,
      \frac{a}q+\frac{\kappa}n,\dots,\frac{a+(q-1)}q+\frac{\kappa}n\\
      \frac1n+\frac{\kappa}n,\dots,\widehat{\frac{n-\kappa}n + \frac{\kappa}n},\dots,\frac{n}n+\frac{\kappa}n
    \end{array}
    \!\!\biggm|\zeta\right)
  \end{sizedisplaymath}
  for $\kappa=0,\dots,n-1$.  For all $a\in\mathbb{C}$, each function is
  defined and analytic in a neighborhood of $\zeta=0$, unless (in the
  first) a lower parameter is a non-positive integer.  If exactly one lower
  parameter equals such an integer, $-m$, and an upper parameter equals an
  integer $-m'$ with $-m\le-m'\le0$, then the function can be defined as a
  limit, with the aid of Lemma~\ref{lem:1}.
\label{def:lastdef}
\end{definition}

  By a congruential argument, the differential equation~$E_n$ corresponding
  to each of these ${}_nF_{n-1}$'s will have reducible monodromy, i.e.,
  have an upper and a lower parameter that differ by an integer, iff
  $qa\in\mathbb{Z}$, resp.\ $na\in\mathbb{Z}$.  In~fact, the~$E_n$ will
  have $\zeta=0$, resp.\ $\zeta=\infty$ as a logarithmic point, i.e., have
  a pair of lower, resp.\ upper parameters that differ by an integer, iff
  $qa\in\mathbb{Z}$, resp.\ $na\in\mathbb{Z}$.

  \begin{theorem}
    \label{thm:inta}
    {\rm (i)} Near the ``\/$\zeta=0$'' point\/
    $[x_1:\ldots:x_n]\in\mathcal{C}_{p,q}^{(n)}$ with
    \begin{equation*}
      x_j=\left\{
      \begin{alignedat}{2}
	&\varepsilon_q^{-(j-1)},&\qquad& j=1,\dots,q,\\
	&0,&\qquad& j=q+1,\dots,n,
      \end{alignedat}
      \right.
    \end{equation*}
    for each\/ $\kappa=0,\dots,q-1$ one has
    \begin{sizemultline}{\small}
      \mathcal{F}_{p,q}^{{(q;\,\kappa)}}(a;\,\zeta) 
      =
      \frac{(-)^\kappa(1)_\kappa(a+n\kappa/q)_{1-\kappa}}{(a)_1} \notag\\
      {}\times\left[(-)^q\frac{p^pq^q}{n^n}\,\zeta\right]^{-\kappa/q}
      \left[(-)^{q-1}\sigma_q\right]^a
      \left[\frac1q\sum_{j=1}^q (\varepsilon_q^{-n})^{(j-1)\kappa} \left(\varepsilon_q^{(j-1)}x_j\right)^{-qa}\right].
      \notag
    \end{sizemultline}
    {\rm (ii)} Near the ``\/$\zeta=\infty$'' point\/
    $[x_1:\ldots:x_n]\in\mathcal{C}_{p,q}^{(n)}$ with
    \begin{equation*}
      x_j=\varepsilon_n^{-(j-1)},\qquad j=1,\dots,n,
    \end{equation*}
    for each\/ $\kappa=0,\dots,n-1$ one has
    \begin{sizemultline}{\small}
      \mathcal{F}_{p,q}^{{(n;\,\kappa)}}(a;\,\zeta^{-1})
      =
      \frac{(-)^\kappa(1)_\kappa(a+q\kappa/n)_{1-\kappa}}{(a)_1} \notag\\
      {}\times\left[(-)^q\frac{p^pq^q}{n^n}\,\zeta\right]^{\kappa/n}
      \left[(-)^{n-1}\sigma_n\right]^a
      \left[\frac1n\sum_{j=1}^n (\varepsilon_n^{-q})^{(j-1)\kappa} \left(\varepsilon_n^{(j-1)}x_j\right)^{-na}\right].
      \notag
    \end{sizemultline}
  In both\/ {\rm(i)} and\/ {\rm(ii),}
  \begin{equation*}
  \zeta \rsmcoloneq (-)^n\frac{n^n}{p^pq^q}\,\frac{{\sigma_n}^q}{{\sigma_q}^n},
  \qquad\quad
  (-)^q\frac{p^pq^q}{n^n}\,\zeta = (-)^p\,\frac{{\sigma_n}^q}{{\sigma_q}^n},
  \end{equation*}
  as usual; and it is assumed that\/ $a\in\mathbb{C}$ is such that the
  left-hand function is defined.  When\/ $\kappa>0$, branches must be
  chosen appropriately.
  \label{thm:82}
  \end{theorem}

In parts (i) and~(ii) of this theorem, the trinomial roots appear in the
form $x_j^{-qa}$, resp.~$x_j^{-na}$.  In consequence, the case when
$a\in\mathbb{Z}$, in particular when $a$~is a negative integer, is
especially nice.  If $\kappa=0$, which will be assumed henceforth, this is
seen as follows.  The identities of parts (i) and~(ii) specialize if
$\kappa=0$ to
\begin{subequations}
\begin{align}
\mathcal{F}_{p,q}^{{(q;\,0)}}(a;\,\zeta)&= \left[(-)^{q-1}\sigma_q\right]^a
\cdot\frac1q\sum_{j=1}^q \bigl(\varepsilon_q^{(j-1)}x_j\bigr)^{-qa}, \label{eq:81pa} \\
\mathcal{F}_{p,q}^{{(n;\,0)}}(a;\,\zeta^{-1})&=
\left[(-)^{n-1}\sigma_n\right]^a \cdot\frac1n\sum_{j=1}^n \bigl(\varepsilon_n^{(j-1)}x_j\bigr)^{-na}.
\label{eq:81pb}\end{align}\label{eq:81p}\end{subequations}When
$a\in\mathbb{Z}$, these reduce to
\begin{subequations}
\begin{align}
\mathcal{F}_{p,q}^{{(q;\,0)}}(a;\,\zeta)&= \left[(-)^{q-1}\sigma_q\right]^a
\cdot\frac1q\sum_{j=1}^q x_j^{-qa}, \label{eq:81a} \\
\mathcal{F}_{p,q}^{{(n;\,0)}}(a;\,\zeta^{-1})&=
\left[(-)^{n-1}\sigma_n\right]^a \cdot\frac1n\sum_{j=1}^n x_j^{-na}.
\label{eq:81b}\end{align}\label{eq:81}\end{subequations}The right side 
of (\ref{eq:81a}), resp.~(\ref{eq:81b}), is \emph{single-valued} on
$\mathcal{C}_{p,q}^{(q)}\subset\mathbb{P}^{q-1}$,
resp.\ $\mathcal{C}_{p,q}^{(n)}\subset\mathbb{P}^{n-1}$.  Therefore, when
$a\in\mathbb{Z}$ the algebraic functions
$\zeta\mapsto\mathcal{F}_{p,q}^{{(q;\,0)}}(a;\,\zeta)$ and
$\zeta\mapsto\mathcal{F}_{p,q}^{{(n;\,0)}}(a;\,\zeta^{-1})$ are
respectively \emph{uniformized} by
$\mathcal{C}_{p,q}^{(q)},\mathcal{C}_{p,q}^{(n)}$.

In fact, a stronger result is true.  There is an obvious permutation
symmetry in~(\ref{eq:81}ab) under $\mathfrak{S}_q$, resp.~$\mathfrak{S}_n$,
which can be quotiented~out, as will now be explained.  Consider the action
of~$\mathfrak{S}_k$ on the curve
$\mathcal{C}_{p,q}^{(k)}\subset\mathbb{P}^{k-1}$, for $k=2,\dots,n$.  It
was shown in~\S\,\ref{sec:curves} that $\mathcal{C}_{p,q}^{(k)}$~is defined
by a system of $k-\nobreak2$ equations, each of which sets to zero a
homogeneous polynomial in the invariants $\bar\sigma_1,\dots,\bar\sigma_k$,
the elementary symmetric polynomials in~$x_1,\dots,x_k$.  The function
field of~$\mathcal{C}_{p,q}^{(k)}$ is the field of rational functions
in~$x_1,\dots,x_k$ that are homogeneous of degree zero, with these defining
equations quotiented~out.

The function field of~$\mathcal{C}_{p,q}^{(k)}$ contains a subfield of
$\mathfrak{S}_k$-stable functions, comprising all rational functions
in~$\bar\sigma_1,\dots,\bar\sigma_k$ that are homogeneous
(in~$x_1,\dots,x_k$) of degree zero; again, with the defining equations
quotiented~out.  One element of this subfield is the function~$\zeta$,
which gives the degree-$(n-\nobreak k+\nobreak 1)_k$ covering
$\pi^{(k)}_{p,q}\colon \mathcal{C}_{p,q}^{(k)}\to\mathbb{P}^1_\zeta$.
Up~to birational equivalence, define a quotiented curve
$\mathcal{C}_{p,q}^{(k;{\rm symm.})}  \rsmcoloneq
\mathcal{C}_{p,q}^{(k)}/\mathfrak{S}_k$ that has this subfield of
$\mathfrak{S}_k$-stable functions as its function field, so that
\begin{equation}
\mathcal{C}_{p,q}^{(k)}\longrightarrow\mathcal{C}_{p,q}^{(k;{\rm
symm.})}\longrightarrow\mathbb{P}^1_\zeta
\end{equation}
is a decomposition of~$\pi^{(k)}_{p,q}$ into maps of respective degrees
$k!,\binom{n}{k}$.  

The genus-zero case $k=2$ is illustrative.  It is the case that
$\mathcal{C}_{p,q}^{(2)}\cong\mathbb{P}_t^1$, where by convention
$t=\allowbreak(x_1+\nobreak x_2)/\allowbreak(x_1-\nobreak x_2)$ is the
rational parameter; and $\mathcal{C}_{p,q}^{(2;{\rm
    symm.})}\cong\mathbb{P}_v^1$, where $v\rsmcoloneq t^2$.  This is
because the only nontrivial element of~$\mathfrak{S}_2$ is the involution
$x_1\leftrightarrow x_2$, which on~account of the parametrization
$[x_1:x_2]=[t+1:t-1]$ is the involution $t\mapsto-t$; so the function field
of~$\mathcal{C}_{p,q}^{(2;{\rm symm.})}$ comprises only functions even
in~$t$.  When $k=1$ or~$0$, $\mathcal{C}_{p,q}^{(k;{\rm symm.})}$~will be
defined as follows.  By convention, $\mathcal{C}_{p,q}^{(1;{\rm
    symm.})}\cong \mathcal{C}_{p,q}^{(1)}\cong\mathbb{P}_s^1$ and
$\mathcal{C}_{p,q}^{(0;{\rm symm.})}\cong
\mathcal{C}_{p,q}^{(0)}\cong\mathbb{P}_\zeta^1$.

\begin{lemma}
$\mathcal{C}_{p,q}^{(k;\rm{symm.})}  \cong
\mathcal{C}_{p,q}^{(n-k;\rm{symm.})}$ for\/ $k=0,\dots,n$, due to the
respective function fields being isomorphic.  Here\/ $n\rsmcoloneq p+q$, as
usual.
\label{lem:ttov}
\end{lemma}

\begin{proof}
  To prove the case $k=n$ (or~$k=0$), observe that the function field of
  ${\mathcal{C}}_{p,q}^{(n;{\rm symm.})}\!$ comprises all rational
  functions of ${\sigma_n}^q/{\sigma_q}^n\propto\zeta$.  (Also,
  $\mathfrak{S}_n$~is the covering group of $\pi^{(n)}_{p,q}\colon
  \mathcal{C}^{(n)}_{p,q}\to\mathbb{P}_\zeta^1$.)  If $0<k<n$, use the fact
  that any element of~${\mathcal{C}}_{p,q}^{(n-k;{\rm symm.})}$ can be
  viewed as a quotient of two symmetric polynomials in~$x_{k+1},\dots,x_n$
  that are homogeneous and of the same degree.  Any homogeneous symmetric
  polynomial in~$x_{k+1},\dots,x_n$ can be expressed as a quotient of two
  such polynomials in~$x_1,\dots,x_k$ by the elimination procedure
  of~\S\,\ref{sec:curves}; cf.~Lemma~\ref{lem:dual}.
\end{proof}

\begin{theorem}
For each\/ $a\in\mathbb{Z}$, 
  \begin{itemize}
    \item[{\rm(i)}] the algebraic function\/
    $\zeta\mapsto\mathcal{F}_{p,q}^{{(q;\,0)}}(a;\,\zeta)$ is
    uniformized by\/ $\mathcal{C}_{p,q}^{(q;{\rm symm.})}$\!.
    \item[{\rm(ii)}] the algebraic function\/
    $\zeta\mapsto\mathcal{F}_{p,q}^{{(n;\,0)}}(a;\,\zeta^{-1})$ is
    uniformized by\/ $\mathcal{C}_{p,q}^{(n;{\rm symm.})}$\!.
  \end{itemize}
\label{thm:strengthened1}
\end{theorem}

\begin{proof}
  Immediate, by the permutation symmetry of formulas
  (\ref{eq:81a}),(\ref{eq:81b}).
\end{proof}

The following is a specialization of Thm.~\ref{thm:strengthened1}, with a
constructive proof.

\begin{theorem}
  For each\/ $a\in\mathbb{Z}$,
  \begin{itemize}
    \item[{\rm(i)}] if\/ $p$ or\/ $q$ equals\/ $1$, resp.\ $2$, the algebraic
    function\/ $\zeta\mapsto\mathcal{F}_{p,q}^{{(q;\,0)}}(a;\zeta)$ can be
    rationally parametrized by\/ $s$, resp.\ $v:=t^2$, where\/ $s,t$ are
    the rational parameters of\/
    $\mathcal{C}_{p,q}^{(1)}\cong\mathbb{P}^1_s$ and\/
    $\mathcal{C}_{p,q}^{(2)}\cong\mathbb{P}^2_t$.  {\rm(}This is because
    $\mathcal{C}_{p,q}^{(q;{\rm symm.})}\cong\mathbb{P}^1_s$, resp.\
    $\mathbb{P}^1_v$.{\rm)}
    \item[{\rm(ii)}] the algebraic function\/
    $\zeta\mapsto\mathcal{F}_{p,q}^{{(n;\,0)}}(a;\zeta^{-1})$ is in fact
    rational,
    i.e., is uniformized by\/
    $\mathcal{C}_{p,q}^{(0)}\rsmcoloneq\mathbb{P}^1_\zeta$. 
  \end{itemize}
In part\/ {\rm(i),} the rational parametrization of the \emph{argument} of
the algebraic function, whether\/ $\zeta=\pi^{(1)}_{p,q}(s)$ or\/
$\zeta=\pi^{(2)}_{p,q}(t)$, was already given in Theorem\/
{\rm\ref{thm:55}}; the latter is even in\/ $t$, and hence, single-valued
in\/ $v$.
\label{thm:83}
\end{theorem}

\begin{proof}
Part~(ii) is trivial, since $\mathcal{F}_{p,q}^{(n;0)}(a;\zeta)$ will be a
polynomial in~$\zeta$ if $a\in\mathbb{Z}$, because one of the upper
hypergeometric parameters will then be a nonpositive integer (see
Defn.~\ref{def:lastdef}).  Also, the cases $q=1,2$ of part~(i) are familiar
from~\S\,\ref{sec:curves}.  If $q=1$ then the right side of~(\ref{eq:81a})
equals $(1-\nobreak s)^a$, and if $q=2$ then one can use the
parametrization $[x_1:x_2]=[t+1:t-1]$ and the formula
for~$\sigma_q=\sigma_q(x_1,x_2)$ given in Lemma~\ref{lem:54} to express the
right side of~(\ref{eq:81a}) as a rational function of the parameter~$t$;
in~fact, of $v=t^2$.

The cases $p=n-q=1,2$ of part~(i) are more interesting.  They can be viewed
as consequences of $\mathcal{C}_{p,q}^{(q;{\rm
    symm.})}\cong\mathcal{C}_{p,q}^{(p;{\rm symm.})}\!$, which comes from
Lemma~\ref{lem:ttov}.  But a constructive rather than an abstract proof
will be given, in the elimination-theoretic spirit of Lemma~\ref{lem:dual}.
Suppose the integer $\gamma=-qa$ is positive.  Write
\begin{equation}
\label{eq:84}
  \sum_{j=1}^q x_j^\gamma = 
\left\{
\begin{alignedat}{2}
  &p_\gamma - x_n^\gamma, &\qquad& p=1,\\
  &p_\gamma - x_n^\gamma - x_{n-1}^\gamma, &\qquad& p=2,
\end{alignedat}
\right.
\end{equation}
where $p_\gamma\rsmcoloneq \sum_{j=1}^nx_j^\gamma$ is a symmetric
polynomial in $x_1,\dots,x_n$, the so-called $\gamma$'th power-sum
symmetric polynomial.  The elementary symmetric polynomials
$\{\sigma_l\}_{l=1}^n$ in $x_1,\dots,x_n$ form an algebraic basis for the
ring of symmetric ones, so the $\{p_\gamma\}_{\gamma=1}^\infty$ can be
expressed as polynomials in the $\{\sigma_l\}_{l=1}^n$; e.g., by inverting
the Newton--Girard formula~\cite[\S\,I.2]{Macdonald95}
\begin{equation}
\label{eq:newtongirard}
\sigma_l = l^{-1}\sum_{\gamma=1}^l (-)^{\gamma-1} p_\gamma\,\sigma_{l-\gamma},
\end{equation}
which holds for all $l\ge1$, it being understood that $\sigma_l=0$ if
$l>n$.  On the curve $\mathcal{C}^{(n)}_{p,q}$, each $\sigma_l$ other than
$\sigma_0,\sigma_q,\sigma_n$ equals zero, which introduces simplifications.

Now, introduce an alternative sequence of subsidiary Schwarz curves
\begin{equation}
  \mathcal{C}^{(n)}_{p,q}\longrightarrow\mathcal{C}^{(n-1)}_{p,q}{}^\prime\longrightarrow\cdots\longrightarrow\mathcal{C}^{(2)}_{p,q}{}^\prime\longrightarrow\mathcal{C}^{(1)}_{p,q}{}^\prime\longrightarrow\mathbb{P}^1_\zeta,
\label{eq:83}
\end{equation}
based on the successive elimination of $x_1,\dots,x_n$ rather than of
$x_n,\dots,x_1$, so that
$\mathcal{C}_{p,q}^{(1)}{}^\prime\cong\mathbb{P}_{s'}^1$, where
\begin{equation}
\label{eq:primedsigmas}
\sigma_q = (-)^{q-1}x_n^q\cdot(1-s'), \qquad\quad  \sigma_n=(-)^{n-1}x_n^n\cdot s',
\end{equation}
(cf.~(\ref{eq:focusedalignat})), and
$\mathcal{C}_{p,q}^{(2)}{}^\prime\cong\mathbb{P}_{t'}^1$ is parametrized by
$t'\rsmcoloneq (x_n+x_{n-1})/(x_n-x_{n-1})$, so that
$[x_n:x_{n-1}]=[t'+1:t'-1]$.  By using
(\ref{eq:84}),(\ref{eq:newtongirard}) and~(\ref{eq:primedsigmas}), one can
write the right side of~(\ref{eq:81a}) as a rational function of $s'$,
resp.\ of~$t'$ (in~fact of $v':={t'}^2$, by invariance under
$x_n\leftrightarrow x_{n-1}$).  Moreover, the expressions for $\zeta$
in~terms of $s',t'$ are the same as those for $\zeta$ in~terms of~$s,t$.

The proof when $\gamma=-qa<0$ is an easy modification of the preceding.
\end{proof}

\begin{remarkaftertheorem}
Illustrations of Thm.~\ref{thm:83}(i), showing how the algorithm embedded
in its proof is implemented, are given in~\S\,\ref{subsec:82} below.
\end{remarkaftertheorem}

Now consider what the two identities of Thm.~\ref{thm:inta} amount~to, when
$a$~is \emph{not} an integer, but nonetheless $qa\in\mathbb{Z}$,
resp.\ $na\in\mathbb{Z}$.  If $a=\frac{m}q$, resp.\ $a=\frac{m}n$, with
$m\in\mathbb{Z}$, the $\kappa=0$ specializations (\ref{eq:81p}ab) reduce to
\begin{subequations}
\begin{align}
\mathcal{F}_{p,q}^{{(q;\,0)}}(\tfrac{m}q;\,\zeta)&=
\left[(-)^{q-1}\sigma_q\right]^{m/q}
\cdot\frac1q\sum_{j=1}^q \varepsilon_q^{-(j-1)m}  x_j^{-m},  \label{eq:86a}
\\
\mathcal{F}_{p,q}^{{(n;\,0)}}(\tfrac{m}n;\,\zeta^{-1})&=
\left[(-)^{n-1}\sigma_n\right]^{m/n}
\cdot\frac1n\sum_{j=1}^n \varepsilon_n^{-(j-1)m} x_j^{-m}.  \label{eq:86b}
\end{align}
\label{eq:86}
\end{subequations}

\noindent
It follows that the power
$\left[\mathcal{F}_{p,q}^{{(q;\,0)}}(\frac{m}q;\,\zeta)\right]^q$, resp.\
$\left[\mathcal{F}_{p,q}^{{(n;\,0)}}(\frac{m}n;\,\zeta^{-1})\right]^n$, is
a rational function of $x_1,\dots,x_q$, resp.\ $x_1,\dots,x_n$, i.e., is in
the function field of $\mathcal{C}_{p,q}^{(q)}$,
resp.~$\mathcal{C}_{p,q}^{(n)}$.  Equivalently, these powers of algebraic
functions are uniformized by
$\mathcal{C}_{p,q}^{(q)},\mathcal{C}_{p,q}^{(n)}$.  But because of the
presence of the coefficients
$\varepsilon_q^{-(j-1)m}\!,\varepsilon_n^{-(j-1)m}\!$, they are not stable
under the action of $\mathfrak{S}_q$ on~$x_1,\dots,x_q$, resp.\
$\mathfrak{S}_n$ on~$x_1,\dots,x_n$.  Rather, they are stable under the
associated subgroups of \emph{cyclic} permutations.  The group of cyclic
permutations of $x_1,\dots,x_k$ will be denoted by~$\mathfrak{C}_k$, and
the group of dihedral ones (if~$k>2$) by~$\mathfrak{D}_k$.  The relevant
order and subgroup indices are $\left|\mathfrak{C}_k\right|=k$ and
$(\mathfrak{D}_k:\mathfrak{C}_k)=2$,
$(\mathfrak{S}_k:\mathfrak{D}_k)=(k-1)!/2$.

\begin{definition}
  The quotient curves $\mathcal{C}_{p,q}^{(k;\text{cycl.})}\!$,
$\mathcal{C}_{p,q}^{(k;\text{dihedr.})}\!$,
$\mathcal{C}_{p,q}^{(k;\text{symm.})}$ are defined to be
$\mathcal{C}_{p,q}^{(k)}/\,\Gamma$ with
$\Gamma=\mathfrak{C}_k,\mathfrak{D}_k,\mathfrak{S}_k$, so that if~$k>2$
one has the sequence of coverings
  \begin{equation*}
    \mathcal{C}_{p,q}^{(k)}
    \longrightarrow
    \mathcal{C}_{p,q}^{(k;\text{cycl.})}
    \longrightarrow
    \mathcal{C}_{p,q}^{(k;\text{dihedr.})}
    \longrightarrow
    \mathcal{C}_{p,q}^{(k;\text{symm.})}
    \longrightarrow
    \mathbb{P}^1_\zeta,
  \end{equation*}
which have respective degrees $k$, $2$, $(k-1)!/2$, and~$\binom{n}{k}$, as
a decomposition of the degree-$(n-\nobreak k+\nobreak1)_k$ covering
map~$\pi_{p,q}^{(k)}\colon\mathcal{C}_{p,q}^{(k)}\to\mathbb{P}^1_\zeta$.
If $k=1,2$ then there is no~dihedral curve, but
$\mathcal{C}_{p,q}^{(k;\text{cycl.})}\to
\mathcal{C}_{p,q}^{(k;\text{symm.})}$ has degree $(k-1)!$ in all cases.
\end{definition}

\begin{theorem}
For each\/ $m\in\mathbb{Z}$,
  \begin{itemize}
    \item[{\rm(i)}] the algebraic function\/
    $\zeta\mapsto\left[\mathcal{F}_{p,q}^{{(q;\,0)}}(\tfrac{m}q;\,\zeta)\right]^q$
    is uniformized by\/ $\mathcal{C}_{p,q}^{(q;{\rm
    cycl.})}\!$.
    \item[{\rm(ii)}] the algebraic function\/
    $\zeta\mapsto\left[\mathcal{F}_{p,q}^{{(n;\,0)}}(\tfrac{m}n;\,\zeta^{-1})\right]^n$
    is uniformized by\/ $\mathcal{C}_{p,q}^{(n;{\rm
    cycl.})}$\!.
  \end{itemize}
\label{thm:strengthened2}
\end{theorem}

\begin{proof}
  Immediate, by the cyclic permutation symmetry of
  (\ref{eq:86a}),(\ref{eq:86b}).
\end{proof}

The following is a specialization of Theorem~\ref{thm:strengthened2}(i),
with a constructive proof.

\begin{theorem}
  For each\/ $m=-1,-2,\dots$,
\begin{itemize}
    \item[{\rm(i)}] for each\/ $q\ge1$, the algebraic function\/
      $\zeta\mapsto\left[\mathcal{F}_{1,q}^{{(q;\,0)}}(\tfrac{m}q;\,\zeta)\right]^q$
      has a\/ ${(q-1)!}$-valued parametrization by\/ $s$, the rational
      parameter of\/ $\mathcal{C}_{1,q}^{(q;{\rm
          symm.})}\cong\mathcal{C}_{1,q}^{(1)}\cong\mathbb{P}^1_s$.  That
      is, it satisfies a\/ degree-$(q-\nobreak1)!$ polynomial equation with
      coefficients in\/ $\mathbb{Z}[s]$.
    \item[{\rm(ii)}] for each odd\/ $q\ge1$, the algebraic function\/
      $\zeta\mapsto\left[\mathcal{F}_{2,q}^{{(q;\,0)}}(\tfrac{m}q;\,\zeta)\right]^q$
      has a\/ ${(q-1)!}$-valued parametrization by\/ $v=t^2$, the rational
      parameter of\/ $\mathcal{C}_{2,q}^{(q;{\rm
          symm.})}\cong\mathbb{P}^1_v$, i.e., by the square of the rational
      parameter of\/ $\mathcal{C}_{1,q}^{(2)}\cong\mathbb{P}^1_t.$ That is,
      it satisfies a\/ degree-$(q-\nobreak1)!$ polynomial equation with
      coefficients in\/ $\mathbb{Z}[v]$.
\end{itemize}
The rational parametrizations $\zeta=\pi^{(1)}_{1,q}(s)$ and
$\zeta=\pi^{(2)}_{2,q}(t)$ of the \emph{argument} of the algebraic
function, for\/ {\rm(i),(ii),} were already given in Theorem\/
{\rm\ref{thm:55};} the latter is even in\/ $t$, and hence, single-valued
in\/ $v$.
\label{thm:87}
\end{theorem}

\begin{proof}
  The rings of symmetric polynomials in $x_1,\dots,x_q$ and
  in~$x_1,\dots,x_n$ have algebraic bases $\{\sigma_l^{(q)}\}_{l=1}^q$
  or~$\{p_\gamma^{(q)}\}_{\gamma=1}^q$, resp.\ $\{\sigma_l^{(n)}\}_{l=1}^n$
  or~$\{p_\gamma^{(n)}\}_{\gamma=1}^n$, of elementary or power-sum
  symmetric polynomials.  Consider
  \begin{equation}
    G_{q,-m}(y) \rsmcoloneq
    \prod_\chi\left\{
    y - \left[
      \sum_{j=1}^q \varepsilon_q^{-(j-1)m}x_{\chi(j)}^{-m}
      \right]^q
    \right\},
  \end{equation}
  the product being taken over one representative~$\chi$ of each of the
  $(q-1)!$ cosets of $\mathfrak{C}_q$ in~$\mathfrak{S}_q$, acting
  on~$x_1,\dots,x_q$.  By Eq.~(\ref{eq:86a}),
  \begin{equation}
    \left[\mathcal{F}_{p,q}^{(q;0)}(\tfrac{m}q;\,\zeta)\right]^q
    = q^{-q}\left[(-)^{q-1}\sigma_q\right]^m\, y_0,
  \end{equation}
  where $y_0$~is a certain root of $G_{q,-m}(y)$.  Each coefficient of the
  degree-$(q-\nobreak1)!$ polynomial $G_{q,-m}(y)$ is stable under the
  action of~$\mathfrak{S}_q$ on~$x_1,\dots,x_q$, and can therefore be
  expressed as a polynomial in the~$\{p_\gamma^{(q)}\}_{\gamma=1}^q$.  But,
\begin{equation}
 p_\gamma^{(q)} =
  \sum_{j=1}^q x_j^\gamma = 
\left\{
\begin{alignedat}{2}
  &p_\gamma^{(n)} - x_n^\gamma, &\qquad& p=1,\\
  &p_\gamma^{(n)} - x_n^\gamma - x_{n-1}^\gamma, &\qquad& p=2.
\end{alignedat}
\right.
\end{equation}
Hence, each coefficient of $G_{q,-m}(y)$ is expressible in~terms of the
$\{\sigma_l^{(n)}\}_{l=1}^n$ and~$x_{n}$, resp.\ $x_{n},x_{n-1}$.  But
on~$\mathcal{C}_{p,q}^{(n)}$, each of the~$\{\sigma_l^{(n)}\}_{l=1}^n$
other than $\sigma_q^{(n)}\!,\sigma_n^{(n)}$\!  equals zero.  By
introducing the alternative sequence~(\ref{eq:83}) of subsidiary Schwarz
curves, one can express $\sigma_q^{(n)}\!,\sigma_n^{(n)}$ and~$x_n$, resp.\
$x_{n},x_{n-1}$, in~terms of $s'$~resp.~$t'$, the rational parameter of
$\mathcal{C}_{1,p}^{(1)}{}^\prime$
resp.~$\mathcal{C}_{2,p}^{(2)}{}^\prime$.  And the expressions for $\zeta$
in~terms of~$s',t'$ are the same as those for $\zeta$ in~terms of~$s,t$.
Each rational function of~$t'$ encountered is a function of
$v'\rsmcoloneq(t')^2$, by invariance under $x_{n}\leftrightarrow x_{n-1}$.
\end{proof}

\begin{remarkaftertheorem}
Illustrations of Thm.~\ref{thm:87}, showing how the algorithm embedded in
its proof is implemented, are given in~\S\,\ref{subsec:83} below.
\end{remarkaftertheorem}

\subsection{Parametrizations when {$a\in\mathbf{Z}$}}
\label{subsec:82}

Theorems\ \ref{thm:dominici} and (especially) \ref{thm:810} below are
sample applications of the algorithm embedded in the proof of
Thm.~\ref{thm:83}(i).  They show how it is possible rationally to
parametrize many algebraic ${}_nF_{n-1}$'s with $a\in\mathbb{Z}$, due~to
their being uniformized by quotient curves $\mathcal{C}_{p,q}^{(q;{\rm
    symm.})}$ that are of genus zero.  These examples also illustrate how
lower parameters that are non-positive integers can be handled by taking a
limit, and applying the auxiliary Lemma~\ref{lem:1}.

The first theorem is relatively simple, since the governing curve is
$\mathcal{C}_{p,q}^{(q)}=\mathcal{C}_{n-1,1}^{(1)}$,\! which is of genus
zero without quotienting; though Lemma~\ref{lem:1} is used.

\begin{theorem}
  For each\/ $n\ge2$, 
  \begin{equation*}
{}_{n-1} F_{n-2}
\left(
\begin{array}{c}
-\frac1n;\frac1n,\dots,\frac{n-2}n\\
\frac1{n-1},\dots,\frac{n-2}{n-1}
\end{array}
    \!\!\biggm|
-\,\frac{n^n}{(n-1)^{n-1}}\,\frac{s}{(1-s)^n}
\right)
= \frac{(n-1)+s}{(n-1)(1-s)}
  \end{equation*}
in a neighborhood of\/ $s=0$.
\label{thm:dominici}
\end{theorem}

\begin{proof}
  As mentioned in the proof of Thm.~\ref{thm:83}(i), if $q=1$ then
  $\mathcal{F}_{p,q}^{{(q;0)}}(a;\zeta)$ equals $(1-\nobreak s)^a$, where
  \begin{equation}
    \zeta = (-)^q\frac{n^n}{p^pq^q}\,\frac{s^q}{(1-s)^n}
  \end{equation}
  is the usual map
  $\phi^{(1)}\colon\mathcal{C}^{(1)}_{p,q}\cong\mathbb{P}^1_s\to\mathbb{P}^1_\zeta$.
  Hence, $\mathcal{F}_{n-1,1}^{{(1;0)}}(-1;\zeta)$ equals $(1-\nobreak
  s)^{-1}$.  Formally,
  \begin{equation}
    \mathcal{F}_{n-1,1}^{{(1;0)}}(-1;\,\zeta) = 
    {}_nF_{n-1}\left(
    \begin{array}{c}
      \frac{-1}n,\frac{0}n,\frac{1}n\dots,\frac{n-2}n\\
      \frac{0}{n-1},\frac{1}{n-1}\dots,
      \frac{n-2}{n-1}
    \end{array}
    \!\!\biggm|\zeta\right),
    \end{equation}
  but the right side is undefined, on account of the presence of a
  non-positive integer $-m=0$ as a lower parameter (accompanied,
  fortunately, by a matching upper parameter, $-m'=0$).  One must interpret
  this in a limiting sense, i.e.,
  \begin{equation}
  \mathcal{F}_{n-1,1}^{{(1;0)}}(-1;\,\zeta) = \lim_{a\to-1}\,
  {}_nF_{n-1}\left(
    \begin{array}{c}
      \frac{a}n,\frac{a+1}n,\frac{a+2}n\dots,\frac{a+(n-1)}n\\
      \frac{a+1}{n-1},\frac{a+2}{n-1}\dots,
      \frac{a+(n-1)}{n-1}
    \end{array}
    \!\!\biggm|\zeta\right).
  \end{equation}
  By Lemma~\ref{lem:1}, this limit equals
  \begin{equation}
    \frac1n + \frac{n-1}n \cdot
    {}_{n-1} F_{n-2}
\left(
\begin{array}{c}
-\frac1n;\frac1n,\dots,\frac{n-2}n\\
\frac1{n-1},\dots,\frac{n-2}{n-1}
\end{array}
    \!\!\biggm|
    \zeta\right),
  \end{equation}
  and the theorem now follows by a bit of algebra.
\end{proof}

\begin{remarkaftertheorem}
  Theorem~\ref{thm:dominici} can also be viewed as a degenerate case of
  Thm.~\ref{thm:71}.  Explicit representations for this algebraic ${}_{n-1}
  F_{n-2}$, for low~$n$, were recently given by
  Dominici~\cite{Dominici2008}.  The differential equation~$E_{n-1}$ of
  which this ${}_{n-1} F_{n-2}$ is a solution has monodromy group $H<{\it
  GL}_{n-1}(\mathbb{C})$ isomorphic to~$\mathfrak{S}_n$, by
  Thm.~\ref{thm:BH1}.
\end{remarkaftertheorem}

The following formula, of Girard type, facilitates computation in the
function fields of any top Schwarz curve and its subsidiaries.

\begin{lemma}
\label{lem:inversenewtongirard}
Consider the subset of\/ $\mathbb{C}^n$, coordinatized by\/
$x_1,\dots,x_n$, on which the elementary symmetric polynomials\/
$\sigma_1,\dots,\sigma_{q-1}$ and $\sigma_{q+1},\dots,\sigma_{n-1}$ equal
zero.  {\rm(}The image of this subset under\/ 
$\mathbb{C}^n\to\mathbb{P}^{n-1}$ is the top curve\/ 
$\mathcal{C}_{p,q}^{(n)}$.{\rm)} On this subset, one can express any
power-sum symmetric polynomial\/ $p_\gamma=\sum_{j=1}^n {x_j}^\gamma$,
$\gamma\ge1$, in terms of\/ $\sigma_q$ and\/ $\sigma_n$ by
  \begin{equation*}
    p_\gamma = \sum c_{m_q,m_n}\,{\sigma_q}^{m_q} {\sigma_n}^{m_n},
  \end{equation*}
  where the sum is over all\/ $m_q,m_n\ge 0$ with\/ $m_qq+m_nn=\gamma$, and
  \begin{equation*}
        c_{m_q,m_n} = (-)^{\chi(q)m_q+\chi(n)m_n}
    \left[
      q\binom{m_q+m_n-1}{m_n} + n\binom{m_q+m_n-1}{m_q}
      \right],
  \end{equation*}
  with\/ $\chi(l)\rsmcoloneq1,0$ if\/ $l\equiv0,1$ {\rm(mod $2$).}
\end{lemma}

\begin{proof}
  This comes, e.g., from the determinantal formula
  \begin{equation}
    p_\gamma =\left|
    \begin{array}{ccccc}
      \sigma_1 & 1 & 0 & \dots & 0 \\
      2\sigma_2 & \sigma_1 & 1 & \dots & 0 \\
      \vdots & \vdots & \vdots & & \vdots \\
      \gamma\sigma_\gamma & \sigma_{\gamma-1} & \sigma_{\gamma-2} & \dots & \sigma_1
    \end{array}
    \right|,
  \end{equation}
with $\sigma_l=0$ if $l>n$.  It is inverse to the Newton--Girard
formula~\cite[\S\,I.2]{Macdonald95}.
\end{proof}

Each hypergeometric function in the following theorem, with
$a\in\mathbb{Z}$, is algebraic with governing curve $\mathcal{C}^{(q;{\rm
    symm.})}_{p,q}=\mathcal{C}^{(3;{\rm symm.})}_{2,3}$\!.  This quotient
curve is of genus zero, although the unquotiented curve
$\mathcal{C}^{(q)}_{p,q}=\mathcal{C}^{(3)}_{2,3}$, which sextuply covers
it, is of genus~$3$ by the formula of Thm.~\ref{thm:genus}.  The point of
the theorem is that if $a\in\mathbb{Z}$, a~uniformization by rational
functions becomes possible.

\begin{theorem}
For each integer\/ $a\le-1$, the function\/
$\mathcal{F}_{2,3}^{{(3;0)}}(a;\zeta)$ satisfies
\begin{equation}
\label{eq:identity23}
  \mathcal{F}_{2,3}^{{(3;0)}}(a;\,\zeta(t))
  = \tfrac13\left[\mathfrak{s}_3(t)\right]^a \left[
P_a\left(\mathfrak{s}_3(t),\mathfrak{s}_5(t)\right) - (t+1)^{-3a} - (t-1)^{-3a}
\right]
\end{equation}
in a neighborhood of\/ $t=0$, where
\begin{gather*}
  {\mathfrak{s}}_3(t)=\frac{1+10t^2+5t^4}{2\,t},\qquad\quad {\mathfrak{s}}_5(t)=-\,\frac{(1-t^2)^2(1+3t^2)}{2\,t},\\
  \zeta(t) = -\,\frac{5^5}{2^23^3}\,\frac{{{\mathfrak{s}}_5}^3(t)}{{{\mathfrak{s}}_3}^5(t)}=\frac{5^5}{3^3}\,\frac{t^2(1-t^2)^6(1+3t^2)^3}{(1+10t^2+5t^4)^5},
\end{gather*}
and\/ $P_a\in\mathbb{Z}[\mathfrak{s}_3,\mathfrak{s}_5]$ is defined by, e.g.,
\begin{equation*}
P_a(\mathfrak{s}_3,\mathfrak{s}_5)=
\left\{
\begin{alignedat}{2}
&3\,\mathfrak{s}_3^{-a},&\qquad & a=-1,-2,-3,-4;\\  
&3\,{\mathfrak{s}_3}^{5}+5\,{\mathfrak{s}_5}^{3},&\qquad & a=-5.
\end{alignedat}
\right.
\end{equation*}
For each\/ $a$, the special function\/
$\mathcal{F}_{2,3}^{{(3;0)}}(a;\zeta)$ can be expressed in terms of
nondegenerate hypergeometric functions by\/ Lemma\/ {\rm\ref{lem:1};} e.g.,
\begin{align*}
 \mathcal{F}_{2,3}^{{(3;0)}}(-1;\,\zeta) &= 
 1\,-\,\frac{2^23^3}{5^5}\,\zeta \cdot
 {}_5F_4\left(
 \begin{array}{c}
   \frac{4}5;\,\frac{6}5,\frac{7}5,\frac{8}5;\,1\\
   \frac32;\,\frac{4}3,\frac{5}3;\,2
 \end{array}
    \!\!\biggm|
    \zeta\right),\\
                                         &= 
    \frac35 + \frac25\cdot{}_4F_3\left(
 \begin{array}{c}
   -\frac15;\,\frac15,\frac25,\frac35\\
   \frac12;\,\frac13,\frac23
 \end{array}
    \!\!\biggm|
    \zeta\right),\\
 \mathcal{F}_{2,3}^{{(3;0)}}(-5;\,\zeta) &= 
 1\,-\,\frac{2^23^2}{5^4}\,\zeta\,-\,\frac{2^63^9}{5^{14}}\,\zeta^3 \cdot
 {}_5F_4\left(
 \begin{array}{c}
   \frac{11}5,\frac{12}5,\frac{13}5,\frac{14}5;\,2\\
   \frac32;\,\frac{10}3,\frac{11}3;\,4
 \end{array}
    \!\!\biggm|
    \zeta\right).
\end{align*}
\label{thm:810}
\end{theorem}

\begin{proof}
  This is the $(p,q)=(2,3)$ case of Thm.~\ref{thm:83}(i), made explicit,
  and the identity~(\ref{eq:identity23}) is a rational parametrization of
  the formula~(\ref{eq:81a}) for~$\mathcal{F}_{p,q}^{{(q;0)}}(a,\zeta)$
  when $a\in\mathbb{Z}$.  The functions
  $\mathfrak{s}_3(t),\mathfrak{s}_5(t)$ are
  $\sigma_3/x_5^3,\sigma_5/x_5^5$, i.e., $\sigma_q/x_n^q,\sigma_n/x_n^n$,
  expressed in~terms of~$t$, the rational parameter of the alternative
  Schwarz curve $\mathcal{C}_{p,q}^{(p)}{}^\prime =
  \mathcal{C}_{2,3}^{(2)}{}^\prime$.  The given expressions and the formula
  for $\zeta=\pi_{p,q}^{(2)}(t)$ come from Lemma~\ref{lem:54} and
  Thm.~\ref{thm:55}; though $x_1,x_2$ are to be replaced by~$x_5,x_4$, and
  `$t$'~is to be understood as $(x_5+\nobreak x_4)/\allowbreak
  (x_5-\nobreak x_4)$, not $(x_1+\nobreak x_2)/\allowbreak (x_1-\nobreak
  x_2)$, so that $[x_5:x_4]=[t+1:t-1]$.

  For each integer $a\le-1$, the quantity
  $P_a=P_a(\mathfrak{s}_3,\mathfrak{s}_5)$ is a normalized version of the
  power-sum symmetric polynomial~$p_\gamma$, i.e.,
  \begin{equation}
    p_\gamma =\sum_{j=1}^n x_j^\gamma = \sum_{j=1}^5 x_j^\gamma,\qquad
    \gamma\rsmcoloneq-qa=-3a,
  \end{equation}
  expressed as a polynomial in~$\sigma_3,\sigma_5$ by the Girard formula of
  Lemma~\ref{lem:inversenewtongirard}.  The subtracted terms
  $(t+\nobreak1)^{-3a}$, $(t-\nobreak1)^{-3a}$ in~(\ref{eq:identity23})
  come from subtracting $x_4^\gamma,x_5^\gamma$ from~$p_\gamma$ to obtain
  $x_1^\gamma+\nobreak x_2^\gamma +\nobreak x_3^\gamma$, as
  in~(\ref{eq:84}).

  The removal from $\mathcal{F}_{2,3}^{(3,0)}(a;\zeta)$, which is a
  ${}_5F_4(\zeta)$, of lower hypergeometric parameters that are
  non-positive integers, is a straightforward application of
  Lemma~\ref{lem:1}.  E.g., for $a=-1,-2,\allowbreak-3,\allowbreak-4,-5$,
  the relevant lower/upper parameters $(-m,-m')$ are respectively $(0,0)$,
  $(0,0)$, $(-1,0)$, $(-1,0)$, $(-2,-1)$.  Only if $a=-1$ or~$-2$ does
  $-m=-m'$, permitting the ${}_5F_4$ to be reduced to a~${}_4F_3$.
\end{proof}

\begin{remarkaftertheorem}
  In the two cases $a=-1,-2$ in which an upper and a lower parameter of
  $\mathcal{F}_{2,3}^{(3,0)}(a;\zeta)$ can be cancelled, reducing it to
  a~${}_4F_3(\zeta)$, the~$E_4$ of which this~${}_4F_3$ is a solution has
  monodromy group $H<{\it GL}_4(\mathbb{C})$ isomorphic
  to~$\mathfrak{S}_5$, by Thm.~\ref{thm:BH1}.

  If $a\le-3$ then $\mathcal{F}_{2,3}^{(3,0)}(a;\zeta)$ is essentially
  a~${}_5F_4(\zeta)$, which is the solution of an~$E_5$ that has reducible
  monodromy, due~to an upper and a lower parameter differing by an integer.
  (See the typical case $a=-5$, below.)
\end{remarkaftertheorem}

\begin{corollary}
\label{cor:811}
  Define a degree-\/$10$ Bely\u\i\ map
  $\mathbb{P}_v^1\to\mathbb{P}_\zeta^1$ by
  \begin{multline*}
    \zeta = \zeta_{2,3}(v)\rsmcoloneq 
    \frac{5^5}{3^3}\,\frac{v(1-v)^6(1+3v)^3}{(1+10v+5v^2)^5}\\
    = 1 - 
    \frac{(1-35v-125v^2-225v^3)^2(27+115v+25v^2+25v^3)}{27\,(1+10v+5v^2)^5}.
  \end{multline*}
  Then in a neighborhood of\/ $v=0$,
  \begin{equation*}
{}_4F_3\left(
 \begin{array}{c}
   -\frac15;\,\frac15,\frac25,\frac35\\
   \frac12;\,\frac13,\frac23
 \end{array}
    \!\!\biggm|
    \zeta_{2,3}(v)
    \right)\\
    = \frac{3+5v^2}{3(1+10v+5v^2)},
  \end{equation*}
  \begin{multline*}
{}_5F_4\left(
 \begin{array}{c}
   \frac{11}5,\frac{12}5,\frac{13}5,\frac{14}5;\,2\\
   \frac32;\,\frac{10}3,\frac{11}3;\,4
 \end{array}
    \!\!\biggm|
    \zeta_{2,3}(v)\right)\\
    =
    \frac{(3+v)(5+10v+v^2)(1+28v+134v^2+92v^3+v^4)(1+10v+5v^2)^{10}}{15\,(1-v)^{18}(1+3v)^9}.
  \end{multline*}
\end{corollary}
\begin{proof}
  These are the rational parametrizations of the nondegenerate
  hypergeometric functions obtained in the theorem (when $a=-1$,
  resp.~$-5$).  They are even in~$t$ and expressible in~terms of
  $v\rsmcoloneq t^2$, as expected.
\end{proof}

\begin{remarkaftertheorem}
  The uniformizing parameter~$v$ is a rational parameter for the genus-$0$
  quotient curve $\mathcal{C}^{(q;{\rm symm.})}_{p,q} =
  \mathcal{C}^{(3;{\rm symm.})}_{2,3}$\!.  The reason why the above
  ${}_4F_3(\zeta)$ and ${}_5F_4(\zeta)$ are $10$-branched functions
  of~$\zeta$ is that the maps
  \begin{equation}
    \mathcal{C}^{(3)}_{2,3} \to
    \mathcal{C}^{(3;{\rm symm})}_{2,3}\cong\mathbb{P}^1_v \to \mathbb{P}^1_\zeta
  \end{equation}
  have respective degrees $6,10$.  The $10$-sheetedness of the latter map
  comes from $\binom{n}q=\binom{p+q}q=\binom{5}{3}=10$.

  The subsidiary curve $\mathcal{C}^{(3)}_{2,3}\subset\mathbb{P}^2$ is of
  genus~$3$ and is a smooth quartic through the fundamental points
  of~$\mathbb{P}^2$, with defining equation
  \begin{multline}
    \label{eq:lastminute}
    (x_1+x_2+x_3)(x_1x_2x_3) + (x_1x_2 + x_2x_3 + x_3x_1)^2 \\
    {}- (x_1+x_2+x_3)^2 (x_1x_2 + x_2x_3 + x_3x_1) = 0,
  \end{multline}
  as follows from Thm.~\ref{thm:long3}.  A formula for the degree-$6$
  quotient map $\mathcal{C}^{(3)}_{2,3}\to\mathcal{C}^{(3;{\rm
      symm.})}_{2,3}\!\cong\mathbb{P}^1_v$ can be derived by eliminating
  $x_4,x_5$ from the defining equations $\sigma_1=0$, $\sigma_2=0$,
  $\sigma_4=0$ of the top curve
  $\mathcal{C}^{(5)}_{2,3}\subset\mathbb{P}^4$, using $v={t'}^2$ with
  $t'=\allowbreak (x_5+\nobreak x_4)/\allowbreak (x_5-\nobreak x_4)$.  It
  is
  \begin{equation}
    \label{eq:degree6}
    v=\frac{(x_1+x_2+x_3)^2}{
      4(x_1x_2+x_2x_3+x_3x_1)
      - 
      3(x_1+x_2+x_3)^2
    }.
  \end{equation}
\end{remarkaftertheorem}

\subsection{Parametrizations when {$qa\in\mathbf{Z}$}}
\label{subsec:83}

Theorems\ \ref{thm:812} and~\ref{thm:fullstop} below are sample
applications of the algorithm embedded in the proof of Thm.~\ref{thm:87},
to the cases $(p,q)=(2,3)$ and~$(1,4)$ with $a=-\frac1q$.  They show how
one can construct parametrizations of many algebraic ${}_nF_{n-1}$'s with
$qa\in\mathbb{Z}$, due~to powers of the ${}_nF_{n-1}$'s being uniformized
by quotient curves $\mathcal{C}_{p,q}^{(q;{\rm cycl.})}\!$ that if not
rational, are at~least low-degree covers of rational ones.  In these two
theorems the rational lower curves will respectively be
$\mathcal{C}_{p,q}^{(q;{\rm symm.})}$ and $\mathcal{C}_{p,q}^{(q;{\rm
    dihedr.})}\!$, and the parametrizations of the ${}_nF_{n-1}$'s will for
the first time involve radicals.

\begin{theorem}
\label{thm:812}
  Define a degree-\/$10$ Bely\u\i\ map\/
  $\mathbb{P}^1_v\to\mathbb{P}^1_\zeta$ by the rational formula for\/
  $\zeta=\zeta_{2,3}(v)$ given in Corollary\/ {\rm\ref{cor:811}.}  Then in
  a neighborhood of\/ $v=0$,
  \begin{multline}
    (1+10v+5v^2)^{1/3}
    \:
    {}_4F_3\left(\!
    \begin{array}{c}
      -\frac1{15},\frac2{15},\frac8{15},\frac{11}{15}\\
      \frac13,\frac23,\frac56
    \end{array}
    \!\!\biggm|
    \zeta_{2,3}(v)
    \right)
    \\
    = \left\{\tfrac12 + \tfrac53v - \tfrac5{54}v^2 + \tfrac12\left[(1+3v)
      \left(1+\tfrac{115}{27}v+\tfrac{25}{27}v^2+\tfrac{25}{27}v^3\right)\right]^{1/2}\right\}^{1/3}\!\!.
    \label{eq:identity23a}
  \end{multline}
\end{theorem}

\begin{proof}
  This is the $(p,q)=(2,3)$ case of Thm.~\ref{thm:87}(ii), made explicit,
  and the identity~(\ref{eq:identity23a}) is a rational parametrization
  of~$\bigl[\mathcal{F}_{p,q}^{{(q;0)}}(a,\zeta)\bigr]^q=\bigl[\mathcal{F}_{2,3}^{{(3;0)}}(a,\zeta)\bigr]^3$
  when $m=-1$, i.e., $a=-1/3$ and $qa=-1$.  

  As in Thm.~\ref{thm:810}, the relevant top and subsidiary Schwarz curves
  are $\mathcal{C}_{p,q}^{(n)}=\mathcal{C}_{2,3}^{(5)}$ and
  $\mathcal{C}_{p,q}^{(q)}=\mathcal{C}_{2,3}^{(3)}$, and the complementary
  subsidiary curve is
  $\mathcal{C}_{p,q}^{(p)}{}^\prime=\mathcal{C}_{2,3}^{(2)}{}^\prime$.  The
  homogeneous coordinates on $\mathcal{C}_{2,3}^{(3)}\subset\mathbb{P}^2$
  are $x_1,x_2,x_3$ and those on
  $\mathcal{C}_{2,3}^{(2)}{}^\prime\cong\mathbb{P}^1_t$ are $x_4,x_5$.  The
  rational parameter on the latter curve is $t\rsmcoloneq\allowbreak
  (x_5+\nobreak x_4)/\allowbreak (x_5-\nobreak x_4)$, and if
  $\mathfrak{s}_3 \rsmcoloneq \sigma_3/x_5^3$, $\mathfrak{s}_5 \rsmcoloneq
  \sigma_5/x_5^5$, then
  \begin{sizegather}{\small}
    {\mathfrak{s}_3}=\frac1{x_5^3}\,\frac{x_5^5 - x_4^5}{x_5^2 - x_4^2}=\frac{1+10t^2+5t^4}{2\,t},
    \qquad
    {\mathfrak{s}_5}= -\,\frac{x_4^2}{x_5^3}\,\frac{x_5^3 -
    x_4^3}{x_5^2-x_4^2}=-\,\frac{(1-t^2)^2(1+3t^2)}{2\,t},  \notag\\
    \zeta=-\,\frac{5^5}{2^2 3^3}\,\frac{\mathfrak{s}_5^3}{\mathfrak{s}_3^5} =
    \frac{5^5}{3^3}\,\frac{t^2(1-t^2)^6(1+3t^2)^3}{(1+10t^2+5t^4)^5}
    \notag
  \end{sizegather}
  follow from Lemma~\ref{lem:54} (with $x_1,x_2$ relabelled
  as~$x_5,x_4$).  This was the origin of the formula $\zeta=\zeta_{2,3}(v)$
  given in Corollary~{\rm\ref{cor:811}}, with $v\rsmcoloneq t^2$.

  By the formula in Eq.~(\ref{eq:86a}), 
  \begin{sizeequation}{\small}
    \mathcal{F}_{2,3}^{(3;0)}(-\tfrac{1}3;\,\zeta) \rsmcoloneq
	    {}_5F_4\left(\!
	    \begin{array}{c}
	      -\frac1{15},\frac2{15},\frac13,\frac8{15},\frac{11}{15}\\
	      \frac13,\frac56;\,\frac13,\frac23
	    \end{array}
	    \!\!\biggm|
	    \zeta
	    \right)
	    = 
	    {}_4F_3\left(\!
	    \begin{array}{c}
	      -\frac1{15},\frac2{15},\frac8{15},\frac{11}{15}\\
	      \frac13,\frac23,\frac56
	    \end{array}
	    \!\!\biggm|
	    \zeta
	    \right)  
  \end{sizeequation}
  (an upper and a lower parameter being cancelled) has the representation
  \begin{subequations}
  \begin{align}
    \mathcal{F}_{2,3}^{(3;0)}(-\tfrac{1}3;\,\zeta) &= 
    \tfrac13\,{\sigma_3}^{-1/3}\left[x_1 + \varepsilon_3\,x_2 +
    \varepsilon_3^2\,x_3\right], \label{eq:817a}\\
    \left[\mathcal{F}_{2,3}^{(3;0)}(-\tfrac{1}3;\,\zeta)\right]^3 &= 
    \tfrac1{3^3}\,{\sigma_3}^{-1}\left[x_1 + \varepsilon_3\,x_2 +
    \varepsilon_3^2\,x_3\right]^3.
    \label{eq:finalsub}
  \end{align}
  \end{subequations}
  Define a polynomial $G_{3,1}(y)$, symmetric in $x_1,x_2,x_3$, one of the
  roots of which is the final factor in~(\ref{eq:finalsub}), as a product
  over the two cosets of $\mathfrak{C}_3$ in~$\mathfrak{S}_3$, i.e.,
  \begin{equation}
    \begin{split}
    G_{3,1}(y) &= \left\{y - \left[x_1 + \varepsilon_3\,x_2 +
    \varepsilon_3^2\,x_3\right]^3\right\} \left\{y - \left[x_1 + \varepsilon_3^2\,x_2 + \varepsilon_3\,x_3\right]^3\right\}\\
    &= y^2 + \left[-2\hat\sigma_1^3 + 9\hat\sigma_1\hat\sigma_2 -
    27\hat\sigma_3\right]y + \left[\hat\sigma_1^2-3\hat\sigma_2\right]^3,
    \end{split}
  \end{equation}
  as one finds by a bit of computation.  Here, $\{\hat\sigma_l\}_{l=1}^3$
  are the elementary symmetric polynomials in $x_1,x_2,x_3$ alone; the ones
  in~$x_4,x_5$ will be denoted by~$\{\bar\sigma_l\}_{l=1}^2$.

  Each coefficient in $G_{3,1}$ can be expressed rationally and
  symmetrically in terms of~$x_4,x_5$, as the function fields of
  $\mathcal{C}_{2,3}^{(3;{\rm symm.})}\!,\mathcal{C}_{2,3}^{(2;{\rm
      symm.})}$ are isomorphic.  One way to do this is to exploit the
  formula given in Lemma~\ref{lem:dual}.  Another uses the structure of the
  ring of symmetric polynomials in~$x_1,\dots,x_5$, and was sketched in the
  proof of Thm.~\ref{thm:87}.  One can express $\{\hat\sigma_l\}_{l=1}^3$
  in~terms of the power-sum symmetric polynomials $\{\hat
  p_\gamma\}_{\gamma=1}^3$, which can be expressed in~terms of the overall
  power-sum symmetric polynomials $\{p_\gamma\}_{\gamma=1}^5$
  in~$x_1,x_2,x_3,x_4,x_5$, along with $x_4,x_5$.  But the
  $\{p_\gamma\}_{\gamma=1}^5$ can be expressed in~terms of
  $\{\sigma_l\}_{l=1}^5$.  Of~these, the only nonzero ones are
  $\sigma_3,\sigma_5$, which were expressed above in~terms of~$x_4,x_5$.

  Regardless of which technique one uses, one finds that
  \begin{equation}
    \label{eq:819}
    G_{3,1}(y) = y^2 + \left[\frac{-7\bar\sigma_1^4 + 36\bar\sigma_1^2\bar\sigma_2-27\bar\sigma_2^2}{\bar\sigma_1}\right]y + \left[-2\bar\sigma_1^2 + 3\bar\sigma_2\right]^3.
  \end{equation}
  Let $F_3$ denote the left side of Eq.~(\ref{eq:finalsub}), so that
  $F_3=y/3^3\sigma_3$, where $y$~is a root of~$G_{3,1}$.  It follows
  from~(\ref{eq:819}) and the formula for $\mathfrak{s}_3 = \sigma_3/x_5^3$
  in~terms of~$t$, and $[x_5:x_4]=[t+1:t-1]$, that $F_3$~satisfies
  \begin{equation}
    {F_3}^2 - \frac{27+90v-5v^2}{27\,(1+10v+5v^2)}\,F_3 -\frac{4v(3+5v)^3}{729\,(1+10v+5v^2)^2} = 0,
  \end{equation}
  where $v=t^2$.  The theorem now follows from the quadratic formula.
\end{proof}

\begin{remarkaftertheorem}
  The cube of $\mathcal{F}_{2,3}^{(3;0)}(-\tfrac13;\zeta)$, i.e., of the
  ${}_4F_3$ in Thm.~(\ref{thm:812}), is uniformized by
  $\mathcal{C}_{2,3}^{(3;{\rm cycl.})}\!$, which doubly covers
  $\mathcal{C}_{2,3}^{(3;{\rm symm.})}\cong\mathbb{P}_v^1$.  (This is clear
  from~(\ref{eq:finalsub}).)  In fact, the statement of the theorem
  implicitly contains a plane model of $\mathcal{C}_{2,3}^{(3;{\rm
      cycl.})}$ as a double cover of $\mathcal{C}_{2,3}^{(3;{\rm
      symm.})}\cong\mathbb{P}_v^1$, namely
  \begin{equation}
  w^2 = (1+3v)
      \left(1+\tfrac{115}{27}v+\tfrac{25}{27}v^2+\tfrac{25}{27}v^3\right).
  \end{equation}
  This affine quartic $\mathcal{C}_{2,3}^{(3;{\rm symm.})}\ni(v,w)$ is
  elliptic (of genus~$1$), with Klein--Weber invariant $j=-2^{12}5^2/3$.
  It~is triply covered by the unquotiented Schwarz curve
  $\mathcal{C}_{2,3}^{(3)}\subset\mathbb{P}^2$, which is of genus~$3$ and
  has defining equation~(\ref{eq:lastminute}).

  In summary, the maps
  \begin{equation}
    \mathcal{C}_{2,3}^{(3)} 
    \longrightarrow
    \mathcal{C}_{2,3}^{(3;{\rm cycl.})}\!
    \longrightarrow
    \mathcal{C}_{2,3}^{(3;{\rm symm.})} \!\cong \mathbb{P}^1_v
    \longrightarrow
    \mathbb{P}^1_\zeta
    \label{eq:slippedin}
  \end{equation}
  have respective degrees $3,2,10$.  This diagram contrasts with
  \begin{equation}
    \mathcal{C}_{2,3}^{(3)} 
    \longrightarrow
    \mathcal{C}_{2,3}^{(2)} \cong\mathbb{P}^1_t
    \longrightarrow
    \mathcal{C}_{2,3}^{(1)} \cong\mathbb{P}^1_s
    \longrightarrow
    \mathbb{P}^1_\zeta,
  \end{equation}
  in which the maps have respective degrees $3,4,5$.

  An $\mathfrak{S}_3$\nobreakdash-invariant formula $v=v(x_1,x_2,x_3)$ for
  the degree-$6$ quotient map
  $\mathcal{C}_{2,3}^{(3)}\to\mathcal{C}_{2,3}^{(3;{\rm
      symm.})}\cong\mathbb{P}^1_v$, which is the composition of the first
  two maps in~(\ref{eq:slippedin}), was given in Eq.~(\ref{eq:degree6}).
  The first of the two, the degree\nobreakdash-$3$ map
  $\mathcal{C}_{2,3}^{(3)}\to\mathcal{C}_{2,3}^{(3;{\rm cycl.})}$\!, is of
  the form $(x_1,x_2,x_3)\mapsto(v,w)$.  A~rational formula for
  $w=w(x_1,x_2,x_3)$, ${\mathfrak{C}}_3$\nobreakdash-invariant rather than
  ${\mathfrak{S}}_3$\nobreakdash-invariant, can also be worked~out.
\end{remarkaftertheorem}

\begin{theorem}
  Define a degree-\/$15$ Bely\u\i\ map\/
  $\mathbb{P}^1_x\to\mathbb{P}^1_\zeta$ by
  \begin{sizealign}{\small}
    \zeta = \zeta_{1,4}(x) :&=
    \frac14\,\frac{x(1-5x)^4 (5+6x+5x^2)^4}{(1-x)^5(1+10x+5x^2)^5}
    =\frac{5^5}{4^4}\frac{s^4}{(1-s)^5}\circ \frac{-(1-5x)(5+6x+5x^2)}{64\,x}
    \notag\\
    &= 1 -
    \frac{(4-5x-10x^2-5x^3)(1-55x-5x^2-5x^3)^2(1+5x^2+10x^3)^2}{4(1-x)^5(1+10x+5x^2)^5}.\notag
  \end{sizealign}
  Then in a neighborhood of\/ $x=0$,
  \begin{multline}
    (1-x)^{1/4}(1+10x+5x^2)^{1/4}
    \:
    {}_4F_3\left(\!
    \begin{array}{c}
      -\frac1{20},\frac3{20},\frac7{20},\frac{11}{20}\\
      \frac14,\frac12,\frac34
    \end{array}
    \!\!\biggm|
    \zeta_{1,4}(x)
    \right)
    \\
    = \left\{\tfrac12 - \tfrac5{16}x - \tfrac5{8}x^2 - \tfrac5{16}x^3 +
    \tfrac12 \left[
      1-\tfrac54x-\tfrac{5}{2}x^2-\tfrac{5}{4}x^3\right]^{1/2}\right\}^{1/4}\!\!.
    \label{eq:identity14}
  \end{multline}
\label{thm:fullstop}
\end{theorem}

\begin{proof}
  This is the $(p,q)=(1,4)$ case of Thm.~\ref{thm:87}(i), made explicit,
  and the identity~(\ref{eq:identity14}) is a rational parametrization
  of~$\bigl[\mathcal{F}_{p,q}^{{(q;0)}}(a,\zeta)\bigr]^q=\bigl[\mathcal{F}_{1,4}^{{(4;0)}}(a,\zeta)\bigr]^4$
  when $m=-1$, i.e., $a=-1/4$ and $qa=-1$.  The relevant top and subsidiary
  Schwarz curves are $\mathcal{C}_{p,q}^{(n)}=\mathcal{C}_{1,4}^{(5)}$ and
  $\mathcal{C}_{p,q}^{(q)}=\mathcal{C}_{1,4}^{(4)}$, and the complementary
  subsidiary curve is
  $\mathcal{C}_{p,q}^{(p)}{}^\prime=\mathcal{C}_{1,4}^{(1)}{}^\prime\cong\mathbb{P}^1_{s'}$.

  The computations resemble those in the proof of Thm.~\ref{thm:812} and
  will only be sketched.  One defines a degree\nobreakdash-$6$ polynomial
  $G_{4,1}(y)$, symmetric in~$x_1,x_2,x_3,x_4$, as a product over the six
  cosets of~$\mathfrak{C}_4$ in~$\mathfrak{S}_4$.  Each coefficient
  in~$G_{4,1}$ can be expressed rationally in $x_5$ and the rational
  parameter $s'=\sigma_5/x_5^5$.  This leads to a degree\nobreakdash-$6$
  polynomial equation for~$F_4$, defined as the fourth power of
  $\mathcal{F}_{1,4}^{(4;0)}(-\tfrac14;\zeta)$, i.e., of the ${}_4F_3$ in
  the theorem.  The coefficients of the degree\nobreakdash-$6$ polynomial
  are polynomial in~$s'$.  By examination, if one substitutes
  \begin{equation}
    s'=s'(x) = \frac{-(1-5x)(5+6x+5x^2)}{64\,x},
  \end{equation}
  this polynomial will factor.  The only relevant factor is a quadratic
  one, namely
  \begin{sizeequation}{\small}
    \label{eq:f4eqn}
    {F_4}^2 - \frac{8-5x-10x^2-5x^3}{8(1-x)(1+10x+5x^2)}\,F_4 + \frac{25\,x^2(1+x)^4}{256(1-x)^2(1+10x+5x^2)^2}=0.
  \end{sizeequation}
  The theorem follows from Eq.~(\ref{eq:f4eqn}) by the quadratic formula.
\end{proof}

\begin{remarkaftertheorem}
  The fourth power of $\mathcal{F}_{1,4}^{(4;0)}(-\tfrac14;\zeta)$, i.e.,
  of the ${}_4F_3$ in Thm.~(\ref{thm:fullstop}), is uniformized by
  $\mathcal{C}_{1,4}^{(4;{\rm cycl.})}\!$, which doubly covers
  $\mathcal{C}_{1,4}^{(4;{\rm dihedr.})}\!$, which in~turn, triply covers
  $\mathcal{C}_{1,4}^{(4;{\rm symm.})}\cong\mathcal{C}_{1,4}^{(1;{\rm
      symm.})}\cong\mathbb{P}_{s'}^1$.  Hence, $x$~can be identified as a
  rational parameter for the genus-$0$ curve $\mathcal{C}_{1,4}^{(4;{\rm
      dihedr.})}$\!.  The statement of the theorem implicitly contains a
  plane model of $\mathcal{C}_{1,4}^{(4;{\rm cycl.})}$ as a double cover of
  $\mathcal{C}_{1,4}^{(4;{\rm dihedr.})}\cong\mathbb{P}_x^1$, namely
  \begin{equation}
  w^2 = 1-\tfrac54 x - \tfrac52 x^2 - \tfrac54 x^3.
  \end{equation}
  This affine cubic $\mathcal{C}_{1,4}^{(4;{\rm cycl.})}\ni(x,w)$ is
  elliptic (of genus~$1$), with Klein--Weber invariant $j=-5^2/2$,
  like~$\mathcal{C}_{1,4}^{(3)}$.  (See Example~\ref{ex:14}; the equality
  of the $j$\nobreakdash-invariants remains to be investigated.)  It~is
  quadruply covered by the top curve
  $\mathcal{C}_{1,4}^{(5)}\cong\mathcal{C}_{1,4}^{(4)}\subset\mathbb{P}^3$,
  which is of genus~$4$.

  In summary, the maps
  \begin{equation}
    \label{eq:lastcite}
    \mathcal{C}_{1,4}^{(5)} \cong
    \mathcal{C}_{1,4}^{(4)} 
    \longrightarrow
    \mathcal{C}_{1,4}^{(4;{\rm cycl.})}\!
    \longrightarrow
    \mathcal{C}_{1,4}^{(4;{\rm dihedr.})} \!\cong \mathbb{P}^1_x
    \longrightarrow
    \mathcal{C}_{1,4}^{(4;{\rm symm.})} \!\cong \mathbb{P}^1_{s'}
    \longrightarrow
    \mathbb{P}^1_\zeta
  \end{equation}
  have respective degrees $4,2,3,5$.  This diagram contrasts with
  \begin{equation}
    \mathcal{C}_{1,4}^{(5)} \cong
    \mathcal{C}_{1,4}^{(4)} 
    \longrightarrow
    \mathcal{C}_{1,4}^{(3)} 
    \longrightarrow
    \mathcal{C}_{1,4}^{(2)} \cong\mathbb{P}^1_t
    \longrightarrow
    \mathcal{C}_{1,4}^{(1)} \cong\mathbb{P}^1_s
    \longrightarrow
    \mathbb{P}^1_\zeta,
  \end{equation}
  in which the maps have respective degrees $2,3,4,5$.

  An $\mathfrak{S}_4$\nobreakdash-invariant formula for the degree-$24$
  quotient map $\mathcal{C}^{(4)}_{1,4}\to\mathcal{C}^{(4;{\rm
      symm.})}_{1,4}\!\cong\mathbb{P}^1_{s'}$, which is the composition of
  the first three of the four maps in~(\ref{eq:lastcite}), can be derived
  by eliminating $x_5$ from the defining equations $\sigma_1=0$,
  $\sigma_2=0$, $\sigma_3=0$ of the top curve
  $\mathcal{C}^{(5)}_{1,4}\subset\mathbb{P}^4$, using $s'=\sigma_5/x_5^5$.
  It is
  \begin{equation}
    s'=\frac{x_1x_2x_3x_4}{(x_1+x_2+x_3+x_4)^4}.
  \end{equation}
\end{remarkaftertheorem}

\begin{acknowledgement}
  We thank the anonymous referee for many suggestions that led to
  improvements in the exposition.
\end{acknowledgement}


\begin{thebibliography}{29}
\expandafter\ifx\csname natexlab\endcsname\relax\def\natexlab#1{#1}\fi
\providecommand{\url}[1]{\texttt{#1}}
\providecommand{\href}[2]{#2}
\providecommand{\path}[1]{#1}
\providecommand{\DOIprefix}{doi:}
\providecommand{\ArXivprefix}{arXiv:}
\providecommand{\URLprefix}{URL: }
\providecommand{\Pubmedprefix}{pmid:}
\providecommand{\doi}[1]{\href{http://dx.doi.org/#1}{\path{#1}}}
\providecommand{\Pubmed}[1]{\href{pmid:#1}{\path{#1}}}
\providecommand{\bibinfo}[2]{#2}
\ifx\xfnm\relax \def\xfnm[#1]{\unskip,\space#1}\fi
\bibitem[{Belardinelli(1959)}]{Belardinelli59}
\bibinfo{author}{G.~Belardinelli}, \bibinfo{title}{Risoluzione analitica delle
  equazioni algebriche generali}, \bibinfo{journal}{Rend. Sem. Mat. Fis.
  Milano} \bibinfo{volume}{29} (\bibinfo{year}{1959}) \bibinfo{pages}{13--45}.
\bibitem[{Berndt(1985)}]{BerndtI}
\bibinfo{author}{B.C. Berndt}, \bibinfo{title}{Ramanujan's Notebooks, {P}art
  {I}}, \bibinfo{publisher}{Springer-Verlag}, \bibinfo{address}{New
  York/Berlin}, \bibinfo{year}{1985}.
\bibitem[{Beukers and Heckman(1989)}]{Beukers89}
\bibinfo{author}{F.~Beukers}, \bibinfo{author}{G.~Heckman},
  \bibinfo{title}{Monodromy for the hypergeometric function {$\sb nF\sb
  {n-1}$}}, \bibinfo{journal}{Invent. Math.} \bibinfo{volume}{95}
  (\bibinfo{year}{1989}) \bibinfo{pages}{325--354}.
\bibitem[{Birkeland(1927)}]{Birkeland27}
\bibinfo{author}{R.~Birkeland}, \bibinfo{title}{{\"U}ber die {A}ufl{\"o}sung
  algebraischer {G}leichungen durch hypergeometrische {F}unktionen},
  \bibinfo{journal}{Math.~Z.} \bibinfo{volume}{26} (\bibinfo{year}{1927})
  \bibinfo{pages}{566--578}.
\bibitem[{Charalambides(1977)}]{Charalambides77}
\bibinfo{author}{C.A. Charalambides}, \bibinfo{title}{A new kind of numbers
  appearing in the {$n$}-fold convolution of truncated binomial and negative
  binomial distributions}, \bibinfo{journal}{{SIAM} J.~Appl. Math.}
  \bibinfo{volume}{33} (\bibinfo{year}{1977}) \bibinfo{pages}{279--288}.
\bibitem[{Chu(1995)}]{Chu95}
\bibinfo{author}{Wenchang Chu}, \bibinfo{title}{Binomial convolutions and
  hypergeometric identities}, \bibinfo{journal}{Rend. Circ. Mat. Palermo~(2)}
  \bibinfo{volume}{43} (\bibinfo{year}{1995}) \bibinfo{pages}{333--360}.
\bibitem[{Comtet(1974)}]{Comtet74}
\bibinfo{author}{L.~Comtet}, \bibinfo{title}{Advanced Combinatorics},
  \bibinfo{publisher}{Reidel}, \bibinfo{address}{Boston/Dordrecht},
  \bibinfo{year}{1974}.
\bibitem[{Dominici(2008)}]{Dominici2008}
\bibinfo{author}{D.~Dominici}, \bibinfo{title}{Asymptotic analysis of
  generalized {H}ermite polynomials}, \bibinfo{journal}{Analysis (Munich)}
  \bibinfo{volume}{28} (\bibinfo{year}{2008}) \bibinfo{pages}{239--261},
  available on-line as arXiv:math/0606324.
\bibitem[{Eagle(1939)}]{Eagle39}
\bibinfo{author}{A.~Eagle}, \bibinfo{title}{Series for all the roots of a
  trinomial equation}, \bibinfo{journal}{Amer. Math. Monthly}
  \bibinfo{volume}{46} (\bibinfo{year}{1939}) \bibinfo{pages}{422--425}.
\bibitem[{Erd{\'e}lyi et~al.(  55)Erd{\'e}lyi, Magnus, Oberhettinger and
  Tricomi}]{Erdelyi53}
\bibinfo{editor}{A.~Erd{\'e}lyi}, \bibinfo{editor}{W.~Magnus},
  \bibinfo{editor}{F.~Oberhettinger}, \bibinfo{editor}{F.G. Tricomi} (Eds.),
  \bibinfo{title}{Higher Transcendental Functions},
  \bibinfo{publisher}{McGraw--Hill}, \bibinfo{address}{New York},
  \bibinfo{year}{1953--55}. \bibinfo{note}{Also known as {\em The Bateman
  Manuscript Project\/}}.
\bibitem[{Fell(1980)}]{Fell80}
\bibinfo{author}{H.~Fell}, \bibinfo{title}{The geometry of zeros of trinomial
  equations}, \bibinfo{journal}{Rend. Circ. Mat. Palermo~(2)}
  \bibinfo{volume}{29} (\bibinfo{year}{1980}) \bibinfo{pages}{303--336}.
\bibitem[{Gessel and Stanton(1982)}]{Gessel82}
\bibinfo{author}{I.~Gessel}, \bibinfo{author}{D.~Stanton},
  \bibinfo{title}{Strange evaluations of hypergeometric series},
  \bibinfo{journal}{{SIAM} J.~Math. Anal.} \bibinfo{volume}{13}
  (\bibinfo{year}{1982}) \bibinfo{pages}{295--308}.
\bibitem[{Glasser(2000)}]{Glasser2000}
\bibinfo{author}{M.L. Glasser}, \bibinfo{title}{Hypergeometric functions and
  the trinomial equation}, \bibinfo{journal}{J.~Comput. Appl. Math.}
  \bibinfo{volume}{118} (\bibinfo{year}{2000}) \bibinfo{pages}{169--173}.
\bibitem[{Gould(1956)}]{Gould56}
\bibinfo{author}{H.W. Gould}, \bibinfo{title}{Some generalizations of
  {V}andermonde's convolution}, \bibinfo{journal}{Amer. Math. Monthly}
  \bibinfo{volume}{63} (\bibinfo{year}{1956}) \bibinfo{pages}{84--91}.
\bibitem[{Gould(1961)}]{Gould61}
\bibinfo{author}{H.W. Gould}, \bibinfo{title}{A series transformation for
  finding convolution identities}, \bibinfo{journal}{Duke Math.~J.}
  \bibinfo{volume}{28} (\bibinfo{year}{1961}) \bibinfo{pages}{193--202}.
\bibitem[{Gould(1972)}]{Gould72}
\bibinfo{author}{H.W. Gould}, \bibinfo{title}{Combinatorial Identities: A
  Standardized Set of Tables Listing 500 Binomial Coefficient Summations},
  \bibinfo{address}{Morgantown, WV}, \bibinfo{year}{1972}.
\bibitem[{Greenfield and Drucker(1984)}]{Greenfield84}
\bibinfo{author}{G.R. Greenfield}, \bibinfo{author}{D.~Drucker},
  \bibinfo{title}{On the discriminant of a trinomial}, \bibinfo{journal}{Linear
  Algebra Appl.} \bibinfo{volume}{62} (\bibinfo{year}{1984})
  \bibinfo{pages}{105--112}.
\bibitem[{Hall(1941)}]{Hall41}
\bibinfo{author}{N.A. Hall}, \bibinfo{title}{The solution of the trinomial
  equation in infinite series by the method of iteration},
  \bibinfo{journal}{Natl. Math. Mag.} \bibinfo{volume}{15}
  (\bibinfo{year}{1941}) \bibinfo{pages}{219--229}.
\bibitem[{Jordan(1965)}]{Jordan65}
\bibinfo{author}{C.~Jordan}, \bibinfo{title}{Calculus of Finite Differences},
  \bibinfo{edition}{3rd} ed., \bibinfo{publisher}{Chelsea},
  \bibinfo{address}{New York}, \bibinfo{year}{1965}.
\bibitem[{Kato(2006)}]{Kato2006}
\bibinfo{author}{M.~Kato}, \bibinfo{title}{Minimal {S}chwarz maps of
  {${}_3F_2$} with finite irreducible monodromy groups},
  \bibinfo{journal}{Kyushu J.~Math.} \bibinfo{volume}{60}
  (\bibinfo{year}{2006}) \bibinfo{pages}{27--46}.
\bibitem[{Kato and Noumi(2003)}]{Kato2003}
\bibinfo{author}{M.~Kato}, \bibinfo{author}{M.~Noumi},
  \bibinfo{title}{Monodromy groups of hypergeometric functions satisfying
  algebraic equations}, \bibinfo{journal}{Tohoku Math.~J.~(2)}
  \bibinfo{volume}{55} (\bibinfo{year}{2003}) \bibinfo{pages}{189--205}.
\bibitem[{Kimura and Shima(1991)}]{Kimura91}
\bibinfo{author}{T.~Kimura}, \bibinfo{author}{K.~Shima}, \bibinfo{title}{A note
  on the monodromy of the hypergeometric differential equation},
  \bibinfo{journal}{Japan. J.~Math. (N.S.)} \bibinfo{volume}{17}
  (\bibinfo{year}{1991}) \bibinfo{pages}{137--163}.
\bibitem[{Lefton(1982)}]{Lefton82}
\bibinfo{author}{P.~Lefton}, \bibinfo{title}{A trinomial discriminant formula},
  \bibinfo{journal}{Fibonacci Quart.} \bibinfo{volume}{20}
  (\bibinfo{year}{1982}) \bibinfo{pages}{363--365}.
\bibitem[{Macdonald(1995)}]{Macdonald95}
\bibinfo{author}{I.G. Macdonald}, \bibinfo{title}{Symmetric Functions and
  {H}all Polynomials}, \bibinfo{edition}{2nd} ed., \bibinfo{publisher}{Oxford
  Univ. Press}, \bibinfo{year}{1995}. \bibinfo{note}{With contributions by
  A.~Zelevinsky}.
\bibitem[{Merlini et~al.(2006)Merlini, Sprugnoli and Verri}]{Merlini2006}
\bibinfo{author}{D.~Merlini}, \bibinfo{author}{R.~Sprugnoli},
  \bibinfo{author}{M.C. Verri}, \bibinfo{title}{Lagrange inversion: When and
  how}, \bibinfo{journal}{Acta Appl. Math.} \bibinfo{volume}{94}
  (\bibinfo{year}{2006}) \bibinfo{pages}{233--249}.
\bibitem[{Mikhalkin(2006)}]{Mikhalkin2006}
\bibinfo{author}{E.N. Mikhalkin}, \bibinfo{title}{On solutions of general
  algebraic equations by means of integrals of elementary functions},
  \bibinfo{journal}{Siberian Math.~J.} \bibinfo{volume}{47}
  (\bibinfo{year}{2006}) \bibinfo{pages}{301--306}. \bibinfo{note}{{R}ussian
  original in {{Sibirsk. Mat. Zh.}} 47 (2006) 365--371}.
\bibitem[{Passare and Tsikh(2004)}]{Passare2004}
\bibinfo{author}{M.~Passare}, \bibinfo{author}{A.~Tsikh},
  \bibinfo{title}{Algebraic equations and hypergeometric series}, in:
  \bibinfo{editor}{O.A. Laudal}, \bibinfo{editor}{R.~Piene} (Eds.),
  \bibinfo{booktitle}{The Legacy of {N}iels {H}enrik {A}bel},
  \bibinfo{publisher}{Springer-Verlag}, \bibinfo{address}{New York/Berlin},
  \bibinfo{year}{2004}, pp. \bibinfo{pages}{653--672}.
\bibitem[{P{\'o}lya(  22)}]{Polya21}
\bibinfo{author}{G.~P{\'o}lya}, \bibinfo{title}{Sur les s{\'e}ries
  enti{\`e}res, dont la somme est une fonction alg{\'e}brique},
  \bibinfo{journal}{Enseignement Math.} \bibinfo{volume}{22}
  (\bibinfo{year}{1921--22}) \bibinfo{pages}{38--47}.
\bibitem[{Zeitlin(1970)}]{Zeitlin70}
\bibinfo{author}{D.~Zeitlin}, \bibinfo{title}{A new class of generating
  functions for hypergeometric polynomials}, \bibinfo{journal}{Proc. Amer.
  Math. Soc.} \bibinfo{volume}{25} (\bibinfo{year}{1970})
  \bibinfo{pages}{405--412}.

\end{thebibliography}

\end{document}